\documentclass[a4paper,11pt]{article}
\title{
\LARGE{Tessellation and Lyubich-Minsky laminations\\
 associated with quadratic maps II:} 
 \\ 
 \Large{Topological structures of 3-laminations
}
}
\author{
Tomoki Kawahira
\\ \small{Graduate School of Mathematics, Nagoya University}
}

\usepackage{enumerate}
\usepackage[dvipdfm,pdftex]{graphicx}
\usepackage{amsmath}
\usepackage{amscd}
\usepackage{amssymb}
\usepackage{amsfonts}
\usepackage{theorem}
\newtheorem{thm}{Theorem}[section]
\newtheorem{prop}[thm]{Proposition}
\newtheorem{lem}[thm]{Lemma}
\newtheorem{cor}[thm]{Corollary}
\theorembodyfont{\rmfamily}

\newtheorem{pf}{Proof.}

\newcommand{\thmref}[1]{Theorem \ref{#1}}
\newcommand{\propref}[1]{Proposition \ref{#1}}
\newcommand{\corref}[1]{Corollary \ref{#1}}
\newcommand{\lemref}[1]{Lemma \ref{#1}}

\newcommand{\figref}[1]{Figure \ref{#1}}
\newcommand{\secref}[1]{Section \ref{#1}}
\newcommand{\A}{\mathbb{A}}
\newcommand{\C}{\mathbb{C}}
\newcommand{\Cstar}{\mathbb{C}^\ast}
\newcommand{\Chat}{\hat{\mathbb{C}}}
\newcommand{\Cbar}{\bar{\mathbb{C}}}
\newcommand{\Dbar}{\bar{\mathbb{D}}}
\newcommand{\R}{\mathbb{R}}
\newcommand{\D}{\mathbb{D}}

\newcommand{\Hyp}{\mathbb{H}}
\newcommand{\Z}{\mathbb{Z}}
\newcommand{\V}{\mathbb{V}}
\newcommand{\N}{\mathbb{N}}
\newcommand{\T}{\mathbb{T}}

\newcommand{\AAA}{\mathcal{A}}%
\newcommand{\BB}{\mathcal{B}}

\newcommand{\DD}{\mathcal{D}}

\newcommand{\FF}{\mathcal{F}}%

\newcommand{\HH}{\mathcal{H}}%
\newcommand{\II}{\mathcal{I}}%
\newcommand{\JJ}{\mathcal{J}}%
\newcommand{\KK}{\mathcal{K}}%
\newcommand{\LL}{\mathcal{L}}
\newcommand{\MM}{\mathcal{M}}%
\newcommand{\NN}{\mathcal{N}}%
\newcommand{\OO}{\mathcal{O}}
\newcommand{\PP}{\mathcal{P}}
\newcommand{\QQ}{\mathcal{Q}}

\newcommand{\SSS}{\mathcal{S}}
\newcommand{\TT}{\mathcal{T}}
\newcommand{\UU}{\mathcal{U}}%
\newcommand{\VV}{\mathcal{V}}

\newcommand{\XX}{\mathcal{X}}
\newcommand{\YY}{\mathcal{Y}}
\newcommand{\ZZ}{\mathcal{Z}}
\newcommand{\norm}[1]{{\left\| #1 \right\|}}

\newcommand{\kakko}[1]{{\left( #1 \right)}}
\newcommand{\skakko}[1]{{\left\{ #1 \right\}}}
\newcommand{\kk}[1]{{\langle #1 \rangle}}

\newcommand{\llist}[1]{{{#1}_1, \ldots, {#1}_l}}

\newcommand{\bo}[1]{{ (#1_0, #1_{-1}, \ldots) }}
\newcommand{\bon}[1]{(#1_{-n})_{n \ge 0}}
\newcommand{\Tate}[2]{\begin{pmatrix} 
                        {#1} \\ 
                        {#2} 
                        \end{pmatrix}}
\newcommand{\QED}{\hfill $\blacksquare$}
\newcommand{\ee}{~=~}
\newcommand{\dee}{~:=~}

\renewcommand{\Bar}{\overline}
\newcommand{\bs}[1]{\boldsymbol{#1}}
\newcommand{\parag}[1]{
\medskip
\noindent {\bfseries #1}
}
\newcommand{\card}{\mathrm{card}}

\newcommand{\rp}{\mathrm{Re}\,}
\newcommand{\ip}{\mathrm{Im}\,}

\newcommand{\dist}{\mathrm{dist}}

\newcommand{\Rat}{\mathrm{Rat}}
\newcommand{\Tess}{\mathrm{Tess}}

\newcommand{\id}{\mathrm{id}}
\newcommand{\h}{\mathrm{h}}
\newcommand{\n}{\mathrm{n}}
\renewcommand{\u}{\mathrm{u}}
\renewcommand{\v}{\mathrm{v}}
\renewcommand{\a}{\mathrm{a}}

\newcommand{\Aff}{\mathrm{Aff}}
\newcommand{\pr}{\mathrm{pr}}
\newcommand{\q}{\mathrm{q}}

\newcommand{\al}{\alpha}
\newcommand{\gam}{\gamma}
\newcommand{\lam}{\lambda}
\newcommand{\Lam}{\varLambda}
\newcommand{\s}{\sigma}
\newcommand{\e}{\epsilon}
\newcommand{\cc}{\circ}
\newcommand{\fe}{f_\epsilon}
\newcommand{\gep}{g_\epsilon}

\newcommand{\he}{h_\epsilon}

\newcommand{\He}{H_\epsilon}
\newcommand{\Kfc}{K_f^\circ}
\newcommand{\Kgc}{K_g^\circ}
\newcommand{\Rf}{\mathcal{R}_f}

\newcommand{\Af}{\AAA_f}
\newcommand{\Ag}{\AAA_g}
\newcommand{\Afn}{\AAA_f^\n}
\newcommand{\Agn}{\AAA_g^\n}
\newcommand{\Afl}{\AAA_f^{\ell}}

\newcommand{\Kf}{\KK_f}

\newcommand{\Hf}{\HH_f}
\newcommand{\Hg}{\HH_g}
\newcommand{\Mf}{\MM_f}
\newcommand{\Mg}{\MM_g}
\newcommand{\Nf}{\NN_f}
\newcommand{\Ng}{\NN_g}
\newcommand{\Jf}{\JJ_f}
\newcommand{\Jg}{\JJ_g}
\newcommand{\Ff}{\FF_f}
\newcommand{\Fg}{\FF_g}

\newcommand{\Pf}{\PP_f}
\newcommand{\Pg}{\PP_g}
\newcommand{\Qf}{\QQ_f}
\newcommand{\Qg}{\QQ_g}
\newcommand{\QPf}{\QQ\PP_f}
\newcommand{\QPg}{\QQ\PP_g}
\newcommand{\Bf}{\BB_f}
\newcommand{\Bg}{\BB_g}
\newcommand{\Df}{\DD_f}
\newcommand{\Dg}{\DD_g}
\newcommand{\Sf}{\SSS_f}
\newcommand{\Sg}{\SSS_g}
\newcommand{\IO}{{\II\OO}}

\newcommand{\QP}{{\QQ\PP}}
\newcommand{\zhat}{\hat{z}}
\newcommand{\what}{\hat{w}}
\newcommand{\fhat}{\hat{f}}
\newcommand{\ghat}{\hat{g}}
\newcommand{\hhat}{\hat{h}}
\newcommand{\hhatA}{\hat{h}_\a}
\newcommand{\hhatH}{\hat{h}_\h}
\newcommand{\hhatM}{\hat{h}_\q}
\newcommand{\psihat}{\hat{\psi}}
\newcommand{\phihat}{\hat{\phi}}
\newcommand{\kaphat}{\hat{\kappa}}
\newcommand{\alhat}{\hat{\alpha}}
\newcommand{\betahat}{\hat{\beta}}

\newcommand{\zetahat}{\hat{\zeta}}
\newcommand{\gamhat}{\hat{\gamma}}
\newcommand{\delhat}{\hat{\delta}}
\newcommand{\omehat}{\hat{\omega}}
\newcommand{\Ihat}{\hat{I}}
\newcommand{\Ohat}{\hat{O}}

\newcommand{\Thetahat}{\hat{\Theta}}
\newcommand{\thetahat}{\hat{\theta}}

\newcommand{\Lhat}{\hat{L}}

\newcommand{\TThat}{\hat{\mathbb{T}}}
\newcommand{\UUhat}{\hat{\mathcal{U}}}
\newcommand{\Kfhat}{\hat{\KK}_f}
\newcommand{\Kghat}{\hat{\KK}_g}
\newcommand{\xhat}{\hat{x}}
\newcommand{\yhat}{\hat{y}}
\newcommand{\xtil}{\tilde{x}}
\newcommand{\ytil}{\tilde{y}}
\newcommand{\lbar}{{\bar{l}}}
\newcommand{\qbar}{{\bar{q}}}
\newcommand{\pihat}{\hat{\pi}}
\newcommand{\Pif}{\Pi_f}
\newcommand{\Pig}{\Pi_g}
\newcommand{\Tf}{T_f}
\newcommand{\Tg}{T_g}
\newcommand{\vz}{\boldsymbol{z}}

\newcommand{\vx}{\boldsymbol{x}}
\newcommand{\vy}{\boldsymbol{y}}
\newcommand{\vv}{\boldsymbol{v}}

\newcommand{\Mfbar}{\overline{\Mf}}
\newcommand{\Mgbar}{\overline{\Mg}}

\begin{document}

\maketitle

\begin{abstract}
According to an analogy to quasi-Fuchsian groups, we investigate topological and combinatorial structures of Lyubich and Minsky's affine and hyperbolic 3-laminations associated with the hyperbolic and parabolic quadratic maps.

We begin by showing that hyperbolic rational maps in the same hyperbolic component have quasi-isometrically the same 3-laminations. 
This gives a good reason to regard the main cardioid of the Mandelbrot set as an analogue of the Bers slices in the quasi-Fuchsian space. 
Then we describe the topological and combinatorial changes of laminations associated with hyperbolic-to-parabolic degenerations (and parabolic-to-hyperbolic bifurcations) of quadratic maps. 
For example, the differences between the structures of the quotient 3-laminations of Douady's rabbit, the Cauliflower, and $z \mapsto z^2$ are described.

The descriptions employ a new method of \textit{tessellation} inside the filled Julia set introduced in Part I \cite{Ka3} that works like external rays outside the Julia set.

\noindent
{\bf The 2000 Mathematics Subject Classification:} 37F45, 37F99.
\end{abstract}

\setcounter{tocdepth}{1}
\tableofcontents

\section{Introduction}\label{sec_01}
As an analogue of the hyperbolic 3-orbifolds associated with Kleinian groups, Lyubich and Minsky \cite{LM} introduced hyperbolic orbifold 3-laminations associated with rational maps. 
For a rational map $f:\Cbar \to \Cbar$, there exist three kinds of laminations $\Af, ~\Hf$, and $\Mf$ as following:
\begin{itemize}
\item The \textit{affine lamination} $\Af$, a Riemann surface lamination with leaves isomorphic to $\C$ or its quotient orbifold.
\item The {$\Hyp^3$-lamination} $\Hf$, a hyperbolic 3-lamination with leaves isomorphic to $\Hyp^3$ or its quotient orbifold. (This is a 3D-extension of the affine lamination of $f$.)
\end{itemize}
As a Kleinian group acts on $\Cbar$ and $\Hyp^3$, there is a natural homeomorphic and leafwise-isomorphic action $\fhat$ on $\Af$ and $\Hf$ induced from that of $f$. In particular, the action is properly discontinuous on $\Hf$. The third lamination is:
\begin{itemize}
\item The \textit{quotient 3-lamination} $\Mf:=\Hf/\fhat$ , whose leaves are hyperbolic 3-orbifolds.  
\end{itemize}
In addition, we have the Fatou-Julia decomposition $\Ff \sqcup \Jf = \Af$ of the affine lamination. As a Kleinian group acts on the complement of the limit set properly discontinuously, $\fhat$ acts on $\Ff$ properly discontinuously. Hence the quotient $\Ff/\fhat$ forms a Riemann surface lamination and we regard it as the \textit{conformal boundary} of $\Mf$. In this paper we say the quotient $\Mfbar = (\Hf \cup \Ff)/\fhat$ the \textit{Kleinian lamination}.

Lyubich and Minsky applied analogous arguments on rigidity theorems of hyperbolic 3-orbifolds to the hyperbolic 3-laminations, and showed a rigidity result of rational maps that have no recurrent critical points or parabolic points \cite[Theorem 9.1]{LM}. 
However, even for hyperbolic quadratic maps like $z^2-1$ or Douady's rabbit, the precise structure of their laminations have not been investigated. (See the questions in \cite[\S 10]{LM}.)
Moreover, even for the simplest parabolic quadratic map $z^2+1/4$ (``the Cauliflower") whose 3-lamination has a cuspidal part, its structure had not been precisely known. 
The aim of this paper is to give a method to describe the topological and combinatorial changes of laminations associated with the motion of parameter $c$ of $z^2+c$ from one hyperbolic component to another via parabolic parameters. 
In particular, our investigation of the Kleinian laminations is explained in terms of the analogy between the Bers slices of the quasi-Fuchsian deformation and the Mandelbrot set. 

\parag{Outline of the paper.}
\secref{sec_02} is devoted for the constructions of the Lyubich-Minsky laminations according to \cite{LM}.  Then we prove that the affine lamination $\Af$ of an expansive (parabolic) rational map $f$ (that is, the Julia set contains no critical point) is actually identified as a proper subset of the \textit{natural extension} (projective limit) $\Nf$ of $f$.

In \secref{sec_03}, we consider perturbation of the rational maps with superattracting cycles. We will show that we can relax such cycles to attracting cycles without changing the geometry of $\Afn$. (In our terminology, attracting cycles are \textit{not} superattracting.) As a corollary, all laminations associated with hyperbolic rational maps in the same hyperbolic component are quasiconformally (or quasi-isometrically) equivalent.

The analogy between the Bers slices of the quasi-Fuchsian deformation space and the Mandelbrot set is explained in \secref{sec_04} in the lamination context. 
Our stance is that the quadratic map $f_c(z)=z^2+c$ is a deformation of $f_0(z)=z^2$, whose Kleinian lamination is a product of Sullivan's solenoidal Riemann surface lamination $\SSS_0$ and a closed interval $[0,1]$. 
By using the results in the \secref{sec_03}, we will show that all $f_c$ with $c$ in the main cardioid of the Mandelbrot set have quasi-isometrically the same Kleinian laminations. 
We also define the lower and upper ends of the quotient 3-laminations  which correspond to the two ends of the Kleinian manifolds of quasi-Fuchsian groups.

Then we go back to the settings of Part I \cite{Ka3} and study degeneration and bifurcation processes of quadratic laminations.

In \secref{sec_05}, we summerize the results and notation given in Part I \cite{Ka3}. 
We formalize hyperbolic-to-parabolic degeneration of quadratic maps in terms of \textit{degeneration pairs} $(f \to g)$ and we recall some properties of \textit{tessellation} of the interior of the filled Julia sets associated with $(f \to g)$. 

From this section, we often use the following standard examples to explain the ideas: 
\begin{itemize}
\item 
{\bf The Cauliflowers}: The family $\{z^2 + c \}$ with $0<c \le 1/4$, which gives the simplest degeneration process. 
\item 
{\bf The Rabbits}: The family $\{z^2 + c \}$ with parameter $c$ moving from $0$ to the center of the $1/3$-limb (\textit{Douady's rabbit}, with the value denoted by $c_{\mathrm{rab}}$) of the Mandelbrot set. (See Figure 1 of Part I \cite{Ka3}.) 
This is a typical degeneration and bifurcation process. 
\item 
{\bf The Airplanes}: 
The family $\{z^2 + c\}$ with $c_{\mathrm{air}} < c \le -1.75$ and $c_{\mathrm{air}}=-1.7548 \cdots$, which gives a real superattracting cycle of period three. This is a typical degeneration process to the root of a primitive copy of the Mandelbrot set.
\end{itemize}

In \secref{sec_06}, we consider hyperbolic-to-parabolic degeneration of the quadratic natural extensions and affine laminations. For hyperbolic $f$ and parabolic $g$ with $f \to g$,  we first lift the pinching semiconjugacy $h$ on the sphere constructed in the first part of this work (\cite{Ka3}) to their natural extensions; $\hhat: \Nf \to \Ng$. Then we describe the structure of $\Ag$ as a degeneration of $\Af$ by using $\hhat$. 
The conclusion here is that as $f \to g$ the topology of $\Af$ is leafwise preserved to $\Ag$ except on finitely many leaves, called the {\it principal leaves}. 


In \secref{sec_07}, we extend the semiconjugacy $\hhat$ on the affine lamination to the $\Hyp^3$-lamination. Though it is impossible to extend it to the whole $\Hf$, we can find the best possible extension which is enough to describe the structure of $\Hg$. 
We will show that the leafwise topology is also preserved as $f \to g$ except on the 3D extension of the principal leaves.

In \secref{sec_08}, we construct a pinching map on the quotient 3-laminations by using the semiconjugacy above. 
This pinching map will give a topological picture of $\Mg$ based on that of $\Mf$.
We will see that the topology of each leaf of $\Mf$ is preserved as $f \to g$, but the geometry changes at the quotient principal leaf. 

The lower ends in the conformal boundaries better reflect this geometric change. 
In \secref{sec_09}, we investigate the structures of the lower ends of the Kleinian laminations $\Bar{\Mf}$ and $\Bar{\Mg}$ by checking how tessellations are lifted to the conformal boundaries. 
Then we will conclude that $\Bar{\Mf}$ cannot be homeomorphic to $\Bar{\Mg}$. 
As an application, we describe the Kleinian laminations of the Cauliflowers (\S \ref{subsec_cauli}) in detail.

In \secref{sec_10}, we consider the topological change of the lower ends associated with a bifurcation. More precisely, we take two distinct degeneration pairs $(f_1 \to g_1)$ and $(f_2 \to g_2)$ with $g_1=g_2$ and we describe the combinatorial difference between the lower ends of $\Bar{\MM_{f_1}}$ and $\Bar{\MM_{f_2}}$ by means of a Dehn-twist-like operation.
For instance, we will conclude:
\begin{thm}\label{thm_introduction_1}
The Kleinian laminations of $f_0:z \mapsto z^2$, the Cauliflower $f_{1/4}:z \mapsto z^2 + 1/4$, and Douady's Rabbit are not homeomorphic each other. 
\end{thm}
See \S \ref{subsec_rabbits} for more detailed description of the  Kleinian lamination of Douady's rabbit. 
We will also give a result which supports a geometric analogy between the Bers slices in the quasi-Fuchsian deformation spaces and the main cardioid of the Mandelbrot set:

\begin{thm}\label{thm_introduction_2}
For any $f_c(z) = z^2 + c$, its Kleinian lamination $\Bar{\MM_{f_c}}$ has a boundary component that is conformally the same as $\SSS_0$. It also has a product structure $\approx \SSS_0 \times [0,1]$ if and only if $c$ is in the main cardioid of the Mandelbrot set. 
\end{thm}
See \S 4.2 and \thmref{thm_prod_str_vs_main_cardioid}.

\parag{H\'enon mappings with small Jacobian.}
There is another viewpoint of the results in \secref{sec_06}. In \cite{HO}, Hubbard and Oberste-Vorth showed a strong connection between the projective limits of hyperbolic polynomials and the Julia sets of H\'enon mappings with small Jacobian. Since Lyubich-Minsky laminations are based on the projective limits of rational maps (the \textit{natural extension} in their/our terminology), our method automatically gives topological description of the Julia sets of some quadratic H\'enon mappings with small Jacobian.
The last section (Appendix \ref{sec_11}) of this paper is a brief survey on relations between the Julia sets of H\'enon mappings and the affine laminations.

\parag{Note.}
Most of pictures in this paper are colored. The most recent version of this paper and author's other articles are available at: 
\verb|http://www.math.nagoya-u.ac.jp/~kawahira|

\parag{Acknowledgments.}
I would like to thank C. Cabrera, K. Ito, M. Lyubich, and Y. Tanaka for stimulating discussions, and the referee for comments. M. Lyubich also gave me opportunities to visit SUNY at Stony Brook, University of Toronto, and the Fields Institute where part of this work was being prepared. I also would like to thank the Fields Institute and the Institute des Hautes \'Etudes Scientifiques for their hospitality while this work was being revised.

This research is partially supported by JSPS Research Fellowships for Young Scientists, JSPS Grant-in-Aid for Young Scientists, Nagoya University, the Circle for the Promotion of Science and Engineering, Inamori Foundation, and the IH\'ES. I sincerely appreciate their supports.

\section{Lyubich-Minsky laminations}\label{sec_02}

In this section we sketch Lyubich and Minsky's universal construction of the affine, $\Hyp^3$-, and quotient laminations associated with general rational maps in \cite[\S 6 and \S 7]{LM}. 
(See \cite{L} for a brief survey by Lyubich, or \cite[\S 3]{KL} would be helpful also. For the general theory of laminations, see \cite{CC}. See also \cite[\S 6]{MS} and \cite{S}.)

Furthermore, for any expansive ($=$ parabolic) rational map $f$, we will show that its affine lamination $\Af$ is identified as its affine part $\Afn$, which is an analytically well-behaved part of the natural extension. 
In particular, for hyperbolic and parabolic quadratic maps $f$ and $g$, we may regard $\Af$ and $\Ag$ as $\Afn$ and $\Agn$.

\subsection{The universal orbifold laminations}
Let $\UU$ be the set of all non-constant meromorphic functions on the complex plane $\C$, with topology induced by the uniform convergence on compact sets. Let us consider the product space $\hat{\UU}:=\UU \times \UU \times \cdots$. We consider each element $\psihat=\bo{\psi} \in \hat{\UU}$ as a map from $\C$ to $\Cbar \times \Cbar \times \cdots$:
\begin{align*}
\psihat:&~\C~\to~\Cbar \times \Cbar \times \cdots \\
    &~w~\mapsto~ (\psi_0(w),\psi_{-1}(w), \ldots ).
\end{align*} 
By taking the value at $w=0$ of this map, we have the following continuous map:
\begin{align*}
\wp:&~\UUhat~\to~\Cbar \times \Cbar \times \cdots \\
    &~\psihat~\mapsto~ (\psi_0(0),\psi_{-1}(0), \ldots ).
\end{align*} 

\parag{Right actions and quotient universal laminations.}
Let $\mathrm{Aff}$ be the set of complex affine maps on $\C$. Then we may identify $\Aff$ as $\A:=\C \times \Cstar$ by regarding $\delta(w)=a+b w$ as $(a,b) \in \A$. 

For $\delta \in \Aff$ and $\psihat \in \UUhat$, $\psihat \cc \delta(w)$ is also an element of $\UUhat$. Thus $\delta \in \Aff$ acts on the right and the orbits $\{\psihat \cc \Aff: \psihat \in \UUhat\}$ form a foliation of $\UUhat$ with leaves isomorphic to $\A$ or its quotient manifold (\cite[Lemma 7.1]{LM}).

We have natural projections from $\A$ over $\Hyp^3 = \C \times \R^+$ and $\C$ given by $(a,b) \mapsto (a, |b|)$ and $(a,b) \mapsto a$. Correspondingly, we consider the following two equivalent relations in $\UUhat$:
\begin{enumerate}
\item $\psihat \sim_\h \psihat'$ if there exists a $t \in \R$ such that $\psihat(w) =\psihat'( e ^{2 \pi i t} w)$ for all $w \in \C$.
\item $\psihat \sim_\a \psihat'$ if there exists a $b \in \Cstar$ such that $\psihat(w) =\psihat'(bw)$ for all $w \in \C$.
\end{enumerate}
Set $\UUhat^\h:= \UUhat/\!\!\sim_\h$ and $\UUhat^\a:=\UUhat/\!\!\sim_\a$. It is known that they are 3- and 2-orbifold foliations \cite[Corollary 7.3]{LM}. In fact, for each $\psihat \in \UUhat$, a natural map 
$$
\A ~\ni ~(a,b) ~\longmapsto~ \psihat(a+bw)~  \in~\psihat \cc \Aff~ \subset~ \UUhat
$$
gives a branched chart for the orbit $\psihat \cc \Aff$ in general. For example, consider the orbit $\psihat \cc \Aff$ of $\psihat(w)=(e^{w^2},e^{w^2},\ldots)$. Then $\psihat(a+bw)$ and $\psihat(-a-bw)$ determine the same point in $\UUhat$ and the chart have singularities over $a=0$. The quotients $(\psihat \cc \Aff)/\!\!\sim_\a$ and $(\psihat \cc \Aff)/\!\!\sim_\h$ may have more singularities in general and thus they form 3- and 2-orbifolds. They also have natural finitely branched charts from $\Hyp^3$ and $\C$ respectively, which is given by projections from $\A$ and its consistent equivalent relations above. We call such charts the \textit{orbifold charts with respect to $\psihat$}.

There is a natural vertical projection $\pr: \UUhat^\h \to \UUhat^\a$ which sends $[\psihat(a+|b|w)]_\h$ to $[\psihat(a+w)]_\a$ for any $\psihat \in \UUhat$ and $(a,b) \in \A$. It is consistent with the vertical projection of orbifold charts that projects $(a,|b|) \in \Hyp^3$ to $a \in \C$. 

In this paper, we will mainly work with non-singular cases: That is, the quotients of the orbit $\psihat \cc \Aff$ are isomorphic to $\Hyp^3$ or $\C$. For $(a,b) \in \A$, we say \textit{$[\psihat(a+bw)]_\h \in \UUhat^\h$ (resp. $[\psihat(a+bw)]_\a \in \UUhat^\a$) has coordinate $(a,|b|) \in \Hyp^3$ (resp. $a \in \C$) with respect to $\psihat$}. 

Note that the continuous map $\wp$ still makes sense from $\UUhat^\h$ or $\UUhat^\a$ to $\Cbar \times \Cbar \times \cdots$.

\parag{Left actions and the global attractor.}
For a rational map $f$ and $\psi \in \UU$, $f \cc \psi$ is also an element of $\UU$. Thus $f$ acts on $\UU$ from the left and we denote this action $f_\u: \psi \mapsto f \cc \psi$. Note that this action is not generally injective. For example, if $f(z)=z^2$ then $f_\u(e^w)=f_\u(e^{w+ \pi i})$. 

Let us consider a closed invariant set $\Kf$ of the semigroup $\kk{f_\u}$ defined by $\bigcap_{n \ge 0} f^n_\u  (\UU)$, which we call the \textit{global attractor} \cite[\S\S 7.3]{LM}. To obtain a homeomorphic action over $\Kf$, we consider the natural extension
$$
\Kfhat:=\skakko{\psihat=\bo{\psi} \in \UUhat: \psi_0\in\Kf, f_\u(\psi_{-n})=\psi_{-n+1}}
$$
with topology induced from $\UUhat$. On this set $f_\u$ is lifted to a homeomorphic action
$$
\fhat_\u~: \Kfhat ~\ni~ \bo{\psi} ~\mapsto~ 
(f_\u(\psi_0),\psi_0, \psi_{-1}, \ldots ) ~\in~ \Kfhat
$$
which is a leafwise isomorphism. It is also known that $\Kfhat^\h:=\Kfhat/\!\!\sim_\h$ and $\Kfhat^\a:=\Kfhat/\!\!\sim_\a$ are 3- or 2-orbifold laminations \cite[Lemma 7.4]{LM}. (In \cite{KL}, the set $\Kfhat^\a$ is denoted by $\boldsymbol{\AAA}_f$ and called the \textit{universal affine lamination}.)
Moreover, the action of $\fhat_\u$ descends to well-defined actions $\fhat_\a$ and $\fhat_\h$ on $\Kfhat^\h$ and $\Kfhat^\a$ as following:
$$
\fhat_\h([\psihat]_\h):=[\fhat_\u(\psihat)]_\h ~~~\text{and}~~~
\fhat_\a([\psihat]_\a):=[\fhat_\u(\psihat)]_\a.
$$
One can easily check that they send a leaf to a leaf isomorphically.

\subsection{The affine and hyperbolic laminations}
\parag{Natural extension.}
For the dynamics of rational $f$, the set of all possible backward orbits
$$
\NN_{f}:=\skakko{\zhat=(z_0,z_{-1}, \ldots)
:~z_0 \in \Cbar,~f(z_{-n-1})=z_{-n}}
$$
is called the \textit{natural extension} of $f$, which is equipped with a topology from $\Cbar \times \Cbar \times \cdots$ \cite[\S 3]{LM}. For $\zhat=\bo{z} \in \Nf$, the lift $\fhat_\n$ of $f$ and the natural projection $\pi_f$ are defined by
\begin{align*} 
\fhat_\n(\zhat)&:=(f(z_0), f(z_{-1}), \ldots)=(f(z_0), z_0, z_{-1}, \ldots) ~~\text{and}\\
\pi_{f}(\zhat) & :=z_0.
\end{align*}
It is clear that $\fhat_\n$ is a homeomorphism, and satisfies $\pi_f \circ \fhat_\n=f \circ \pi_f$. 

\parag{Remarks on notation.}
\begin{itemize}
\item For simplicity, we will abuse $\fhat$ to denote $\fhat_\u$, $\fhat_\h$, $\fhat_\a$, and $\fhat_\n$ when it does not lead confusion.

\item We will use $\xtil$ or $\ytil$ to denote points in $\Kfhat^\h$, and $\xhat$ or $\yhat$ to denote points in $\Kfhat^\a$. The points in the natural extension $\Nf$ are denoted by $\zhat$ or $\zetahat$, etc. 
\end{itemize}

\parag{Regular leaf space and its affine part.}
An element $\zhat=(z_0,z_{-1}, \ldots) \in \mathcal{N}_f$ is \textit{regular} if there exists a neighborhood $U_0$ of $z_0$ such that its pull-backs 
$$
z_0 \in U_0~\leftarrow~ U_{-1} ~\leftarrow~ U_{-2} ~\leftarrow~  \cdots
$$
along the backward orbit $\zhat$ are eventually univalent. For example, backward orbits remaining in an attracting or parabolic cycle are \textit{not} regular, since at least one critical orbit accumulates on such a cycle. 
We say that such backward orbits form an attracting or parabolic cycle in $\Nf$.

Let $\Rf$ denote the set of all regular points in $\Nf$. We call $\Rf$ the \textit{regular leaf space} of $f$ \cite[\S\S 3.2]{LM}. A \textit{leaf} of $\Rf$ is a path-connected component of $\Rf$. Note that $\Rf$ is invariant under $\fhat$. It is known that each leaf of $\Rf$ has a natural complex structure of non-compact Riemann surface induced from the dynamics downstairs. We take the union of all the leaves isomorphic to $\C$, and call it the \textit{affine part} $\Afn$ of $\Rf$ \cite[\S\S 4.2]{LM}. 
By Lemma 3.1, Proposition 4.5, Lemma 4.8 of \cite{LM}, we have

\begin{prop}[Affine part]\label{prop_LM_4.5}
The affine part $\Afn$ satisfies the following properties:
\begin{enumerate}[\rm (1)]
\item The map $\fhat$ acts on $\Afn$ homeomorphically. Moreover, for any leaf $L$ in $\Afn$, $\fhat:L \to \fhat(L)$ is an isomorphism.
\item Any leaf $L$ of $\Afn$ is dense in $\Nf$.
\item If all critical points of $f$ are non-recurrent, 
$$
\Afn \ee \Rf \ee \Nf-\skakko{\text{attracting and parabolic cycles}}.
$$
\end{enumerate}
\end{prop}
Note that we may regard the action $\fhat:L \to \fhat(L)$ an affine map by (1), since both $L$ and $\fhat(L)$ are isomorphic to $\C$.

We will apply this proposition to hyperbolic and parabolic quadratic maps in \secref{sec_06}.

\parag{Affine lamination.}
Take a point $\zhat =\bo{z} \in \Af^\n$ and the leaf $L=L(\zhat)$ which contains $\zhat$. Then there is a uniformization $\phi:\C \to L$ with $\phi(0)=\zhat$. (Such a $\phi$ is not unique.) One can easily check that the sequence of maps $\{\psi_{-N}=\pi_f \cc \fhat^{-N} \cc \phi: \C \to \Cbar\}_{N \ge 0}$ determines an element $\psihat=(\psi_{-N})_{N \ge 0}$ in $\Kfhat$ with $\wp(\psihat)=\zhat$. Moreover, the element $\xhat=[\psihat]_\a$ in $\Kfhat^\a$ does not depend on the choice of $\phi$ above and is uniquely determined. Let $\iota=\iota_f: \Afn \to \Kfhat^\a$ denote the inclusion defined by $\iota(\zhat)=\xhat$. Note that $\iota$ is injective and $\wp \cc \iota =\id$ on $\Afn$. 
For later use we set $\pihat_f:=\pi_f \cc \wp:\Kfhat^\a \to \Cbar$.

Now we claim that $\iota$ sends $L=L(\zhat)$ above to a leaf $\Lhat=\Lhat(\xhat)$ of $\Kfhat^\a$ isomorphic to $\C$: In fact, for any $\zetahat \in L$, there exists a unique $a \in \C$ such that $\zetahat=\phi(a)$ thus $\zetahat$ and the element of the form $[\psihat(a+bw)]_\a~(b \in \Cstar)$ in $\Kfhat^\a$ have one-to-one correspondence. This means that the quotient orbit $(\psihat \cc \Aff)/\!\!\sim_\a$ is precisely parameterized by $\C$ and $\iota(\zetahat)$ has coordinate $a$ with respect to $\psihat$. We say $\phihat:=\iota \cc \phi$ is a \textit{uniformization} of $\Lhat$ with respect to $\psihat$.

Let $\Afl$ denote $\iota(\Af^\n)$ in $\Kfhat^\a$. Then $\Afl$ consists of leaves isomorphic to $\C$ and is invariant under $\fhat$. Now we take the closure $\Af$ of $\Afl$ in $\Kfhat^\a$ and call it the \textit{affine lamination} of $f$. There may be some additional orbifold leaves in $\Af-\Afl$ in general. However, in the case where $f$ is hyperbolic or parabolic (quadratic) maps, we will show that $\Afl$ is already closed in $\Kfhat^\a$, and thus $\Af=\Afl$ (\propref{prop_Afn_is_Af_if_para}).

\parag{$\bs{\Hyp^3}$-lamination.}
For $\xhat \in \Af$ and one of its representative $\psihat \in \Kfhat$ (that is, $\xhat=[\psihat]_\a$), the quotient orbit $(\psihat \cc \Aff)/\!\!\sim_\h$ is a leaf of $\Kfhat^\h$ with orbifold chart $(a, |b|) \mapsto [\psihat(a+|b|w)]_\h$ with respect to $\psihat$. Now we gather such leaves and define the \textit{$\Hyp^3$-lamination} associated with $f$ as following:
$$
\Hf~:=~\{ \psihat \in \Kfhat~: ~[\psihat]_\a \in \Af\}~/\! \sim_\h
$$
This is one of our hyperbolic 3-orbifold laminations.
Note that the $\Hyp^3$-lamination is given by taking a suitable $\R^+$-bundle of the affine leaves with vertical projection $\pr: \Hf \to \Af$. 
Indeed, $\Hf$ can be geometrically constructed as an ``scaling bundle" of $\Af$ \cite{KL}. But the virtue of our universal construction is that $\Hf$ is a sublamination of the universal metrized space $\UUhat^\h$ for any rational map $f$.

In this paper, for a subset $X$ of $\Af$, we denote its fiber $\pr^{-1}(X)$ by $X^\h$. (In \cite{KL}, such an $X^\h$ is denoted by $\mathfrak{H}X$ and they say $\mathfrak{H}: X \mapsto X^\h$ is the \textit{hyperbolization functor.}) For example, for a leaf $\Lhat$ of $\Af$, $\Lhat^\h$ is a leaf of $\Hf$ with $\pr(\Lhat^\h)=\Lhat$. For the orbifold chart $\phihat:\C \to \Lhat$ with respect to a fixed $\psihat \in \Kfhat$, let $\phihat^\h$ denote the orbifold chart $\phihat^\h:\Hyp^3 \to \Lhat^\h$ with respect to $\psihat$. Then $\pr:\Lhat^\h \to \Lhat$ is consistent with the vertical projection from $\Hyp^3$ to $\C$; that is, $\pr\cc \phihat^\h(a, |b|)=\phihat(a)$.

\parag{Quotient 3-lamination.}
Now $\Hf$ is invariant under the action $\fhat=\fhat_\h$ over $\Kfhat^\h$. It is known that the action of the cyclic group $\kk{\fhat}$ on $\Hf$ is properly discontinuous (\cite[Proposition 6.2]{LM}). Thus the quotient $\Mf:=\Hf/\fhat$ is Hausdorff, and inherits the laminar structure of $\Hf$. We call $\Mf$ the \textit{quotient 3-lamination} of $f$. This is our second hyperbolic 3-lamination, which corresponds to the hyperbolic 3-orbifold associated with a Kleinian group. 

Since $\fhat$ carries leaves of $\Hf$ isometrically, a leaf $\ell$ of $\Mf$ is isomorphic to an hyperbolic 3-orbifold. If a leaf $\Lhat^\h$ of $\Hf$ is periodic, i.e. $\fhat^n(\Lhat^\h) = \Lhat^\h$ for some $n \in \Z$, then $\ell=\Lhat^\h/\fhat$ is either a solid torus or a rank one cusp; or it has a finite branched cover from $\Hyp^3$. If $\Lhat^\h$ is not periodic, the quotient is isomorphic to $\Hyp^3$. In \cite[Prop.1.22]{KL} one can find a complete classification of such leaves.

\begin{figure}[htbp]
\begin{center}
\includegraphics[width=.7\textwidth]{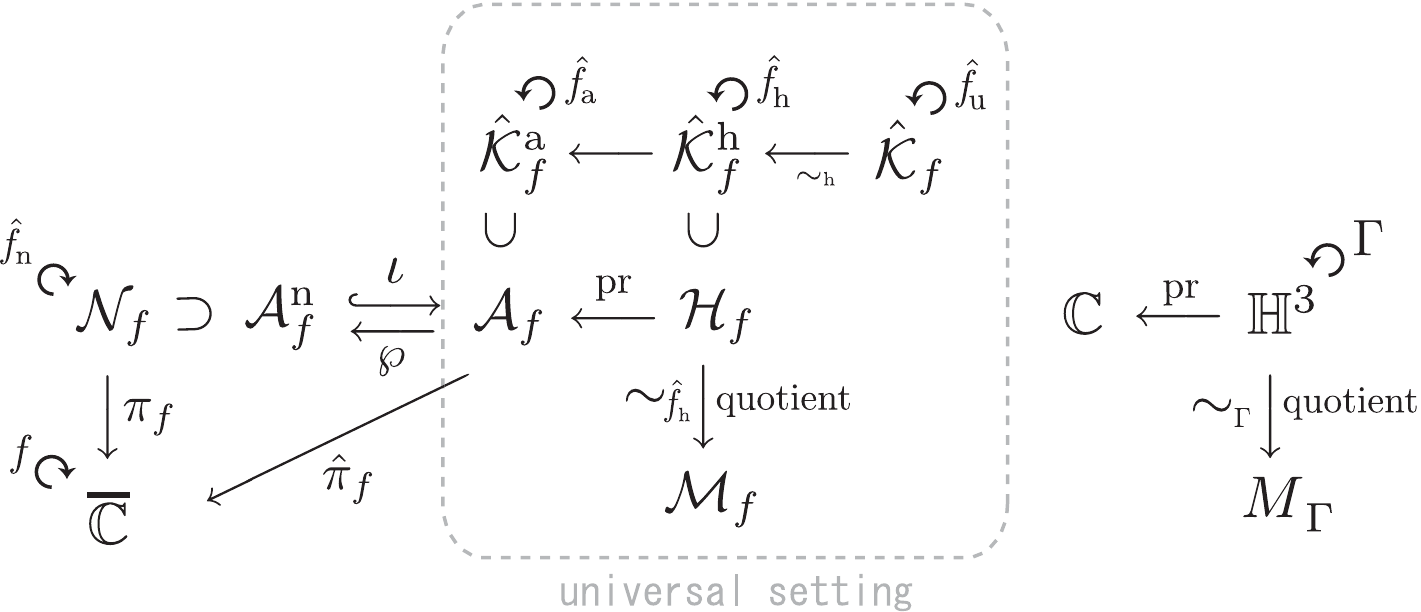}
\end{center}
\caption{The objects in the dashed square are defined by means of the universal laminations. On the right $\Gamma$ is a Kleinian group acting on $\Hyp^3$. }\label{fig_diagram_laminations}
\end{figure}

\figref{fig_diagram_laminations} is a diagram showing the relation between the objects we have defined so far. 

\parag{Conformal boundary and Kleinian lamination.}
According to the Fatou-Julia decomposition $\Cbar=F_f \sqcup J_f$ of the dynamics on $\Cbar$, we have a natural decomposition
$$
\Af \ee \Ff \sqcup \Jf \dee \pihat_f^{-1}(F_f) \sqcup \pihat_f^{-1}(J_f)
$$
of the affine lamination. 
Both $\Ff$ and $\Jf$ are invariant under $\fhat$. 
In particular, $\fhat$ acts on $\Ff$ properly discontinuously \cite[\S\S 3.5]{LM}. The quotient $\partial \Mf:=\Ff/\fhat$ is a Riemann surface lamination (possibly with many components) and called the \textit{conformal boundary} of $\Mf$. 
The union $\overline{\Mf}=\Mf \cup \partial \Mf$ forms a 3-lamination with boundary. We call $\overline{\Mf}$ 
the \textit{quotient 3-lamination with conformal boundary}, 
or simply, 
the \textit{Kleinian lamination} of $f$ after ``Kleinan manifolds" in the theory of hyperbolic 3-manifolds.

\subsection{Affine laminations for expansive rational maps.}
We say a rational map $f$ is \textit{expansive} (or \textit{parabolic}) if any two distinct orbits in the Julia set are eventually a definite distance away from each other. By Denker and Urba\`nski \cite{DU}, it is shown that $f$ is expansive if and only if the Julia set $J_f$ contains no critical points. Equivalently, every critical point of $f$ is attracted to an attracting or parabolic cycle. Note that every hyperbolic ($=$ expanding) map is a special case of expansive map according to this terminology. Here we will establish:

\begin{prop}\label{prop_Afn_is_Af_if_para}
For any expansive rational map $f$, the set $\Afl$ is closed in $\Kfhat^\a$, and thus $\Af=\Afl$.
\end{prop}

Before the proof, let us consider the case of quadratic map $f=f_c: z \mapsto z^2+c$. A direct corollary of this proposition is:
\begin{cor}\label{cor_Afn_is_Af}
If $f=f_c$ has an attracting or parabolic cycle in $\C$, then the affine lamination $\Af$ is equal to $\Afl$. 
\end{cor}
Thus we can identify $\Afn$ in the natural extensions (or in $\Cbar \times \Cbar \times \cdots$) as $\Af$ in the universal affine laminations (or in $\UUhat^\a$). More informally, we may consider that $\zhat$ is just another name of $\xhat=\iota_f(\zhat)$, and vice versa. 

\parag{Proof of \propref{prop_Afn_is_Af_if_para}.}
Take any sequence $\xhat_n \in \Afl$ with $\xhat_n \to \xhat \in \Af$. Set $\zhat_n=(z_{n,-N})_{N \ge 0}:=\wp(\xhat_n)$ and $\zhat=(z_{-N})_{N \ge 0}:=\wp(\xhat)$. By continuity of $\wp$, we have $\zhat_n \to \zhat$ in $\Afn$. Set $\yhat:=\iota(\zhat) \in \Afl$. Now it is enough to show that $\xhat_n \to \yhat$ in $\Af$ and thus $\yhat=\xhat$. If $f$ is hyperbolic, this is straightforward by applying \cite[Propositions 7.5 and 8.2]{LM}. If $f$ has parabolic cycles, we need some modification.

\parag{Backward invariant set $\Omega$.}
Let $O$ be an attracting or parabolic cycle. If $O$ is attracting, we take an open set $E_O$ as a neighborhood of $O$ with $\overline{f(E_O)} \subset E_O$. If parabolic, we take $E_O$ as attracting petals attached to each point of $O$ such that $\overline{f(E_O)}-O \subset E_O$. Let $E$ be the union of such $E_O$ for all attracting and parabolic cycles. Since $f$ is parabolic, there exists an integer $M>0$ such that the compact set $\Omega:=\Cbar-f^{-M}(E)$ contains no critical or postcritical point except parabolic cycles. Note that $f^{-1}(\Omega)-(\text{parabolic cycles}) \subset \Omega^\cc$ and thus every backward orbit except attracting or parabolic cycles is eventually trapped into $\Omega^\cc$ and accumulates on the Julia set.

\parag{Convergence on any large disks.} 
For $\rho>0$, let $\D_\rho$ denote the disk centered at $0$ with radius $\rho$. Let us fix a representative $\psihat=(\psi_{-N})_{N \ge 0} \in \Kfhat$ of $\yhat$. We claim that \textit{for arbitrarily large $N_0>0$ and $\rho>0$, we can choose a representative $\psihat_n=(\psi_{n,-N})_{N \ge 0} \in \Kfhat$ of $\xhat_n$ such that $\psi_{n,-N_0+m}$ converges to $\psi_{-N_0+m}$ uniformly on $\D_\rho$ for $m=1, \ldots, N_0$}.

Take a sufficiently large $\rho'>0$ such that $\e=\rho/\rho' \ll 1$, and set $\V:=\D_{\rho'}$. We first claim that for all $N \gg 0$, $V_{-N}:=\psi_{-N}(\V) \Subset \Omega$. Set $L:=L(\zhat)$ in $\Afn$, the leaf containing $\zhat$. Then there exists a unique uniformization $\phi:L \to \C$ such that $\psi_0=f^N \cc \psi_{-N}=f^N \cc (\pi_f \cc \fhat^{-N} \cc \phi)$. 
Since $\deg (\psi_0|_{\V})$ is constant and finite, so is $\deg(f^N \cc \pi_f \cc \fhat^{-N} \cc \phi) = \deg (f^N \cc \pi_f)$ for all $N \gg 0$. Hence $f:V_{-N} \to V_{-N+1}$ is univalent for all $N \gg 0$. This implies that the backward sequence of open sets $V_0 \leftarrow V_{-1} \leftarrow \cdots$ do not contain the backward orbits corresponding to the invariant lift of the attracting or parabolic cycles. Since $\overline{V_0}$ is compact, eventually $V_{-N}$ is contained in $\Omega^\cc$ for all $N \gg 0$. Note that for such an $N$, the map $\pi_f \cc \fhat^{-N} :\phi(\V) \to V_{-N}$ gives a local chart of $L$ around $\zhat$ (cf. \cite[Lemma 3.1]{LM}) and it implies that $\psi_{-N}:\V \to V_{-N}$ is univalent for all $N \gg 0$.

For each $n$, take a representative $\psihat_n=(\psi_{n,-N})_{N \ge 0} \in \Kfhat$ with $\wp(\psihat_n)=\zhat_n$. Fix a large $N \ge N_0$ so that $V_{-N} \Subset \Omega$. Since $\zhat_n$ converges to $\zhat$ in $\Afn$, we have $z_{n,-N} \in V_{-N}$ for all $n \gg 0$. Thus any pull-back of $V_{-N}$ along $\zhat_n$ is univalent for all $n \gg 0$. Now we can apply the argument of \cite[Proposition 7.5]{LM}: Since each $\psi_{n,-N}$ is also induced from a uniformization of a leaf $L(\zhat_n)$ of $\Afn$, there exists a neighborhood $\V_n \subset \C$ of $0$ such that $\psi_{n,-N}$ maps $\V_n$ univalently onto $V_{-N}$. Consider a univalent map $\chi_n=\chi_{n, -N}: \V \to \V_n$ defined by the univalent branch of $\psi_{n,-N}^{-1} \cc \psi_{-N}$ via $V_{-N}$. Then we have $\psi_{-N+m}=\psi_{n,-N+m} \cc \chi_{n}$ on $\V$ for all $m =0, 1, \ldots, N$. 

Let $a_n$ be the point in $\V$ with $\chi_{n}(a_n)=0$. Set $\tilde{\chi}_{n}(w):=\chi_{n}(w+a_n)$ which is defined on the slided disk $\V- a_n$. By normalizing $\psihat_n$ such that $(\psi_{-N})'(a_n)=(\psihat_{n,-N})'(0)$, we have $\tilde{\chi}_{n}(0)=0$ and $(\tilde{\chi}_{n})'(w)=1$. Since $a_n=\psi_{-N}^{-1}(z_{n,-N})$ tends to $\psi_{-N}^{-1}(z_{-N})=0$, we have 
$$
\D_{\rho}-a_n ~\Subset~ \D_{2\rho} ~\Subset~ \D_{\rho'/2} ~\Subset~ \V-a_n
$$ 
for all $n \gg 0$. By applying the Koebe distortion theorem to $\tilde{\chi}_{n}$ on $\D_{2\rho} \Subset \D_{\rho'/2}$, we have $\tilde{\chi}_{n}(w)=w(1+O(\e))$ on $\D_{2\rho}$ and thus $\chi_n(w)=(w-a_n)(1+O(\e))$ on $\D_{\rho}$. This implies that $\chi_{n} \to \id$ uniformly on $\D_{\rho}$ as $n \to \infty$, and thus $\psi_{n,-N+m} \to \psi_{-N+m}$ for all $m =0, 1, \ldots, N$ on $\D_{\rho}$. We conclude that $\xhat_n$ converges to $\xhat=\yhat$ within $\Afl$. \QED

\parag{Tameness of the laminations.}
After \cite[\S 3.1.7]{KL}, we say a lamination is {\it tame} if it is locally compact. 
Since $\UUhat^\a$ and $\UUhat^\h$ are not locally compact, the tameness of $\Af$ or $\Hf$ of expansive $f$ is not obvious.
However, since $\Afn$ and $\Af$ have the same topology and $\Afn$ is locally compact, we have:

\begin{cor}[Tameness]\label{cor_loc_cpt}
The affine lamination $\Af$ of the expansive rational map $f$ is tame. Hence the hyperbolic 3-laminations $\Hf$ and $\Mf$ are tame too. 
\end{cor}

The local compactness of the quotient 3-lamination would play an important role when one try to extend Lyubich and Minsky's rigidity result of the critically non-recurrent maps with no parabolic points \cite[Theorem 9.1]{LM} in the lamination context. However, Ha\"{\i}ssinsky \cite{Ha} extended their theorem without using the lamination theory to wider class of rational maps called \textit{(uniformly) weakly hyperbolic}, which may have parabolic cycles. For more results on the rigidity, see references in \cite{Ha}.

\section{Stability of hyperbolics}\label{sec_03}

In this section we check the topological/geometric stability of laminations associated with hyperbolic rational maps with superattracting cycles. In our quadratic case, it is justified that the laminations of $f=f_c$ in a hyperbolic component and those of its hyperbolic center have topologically the same structures. 

The key result is that we can perturb superattracting cycles into attracting cycles without changing most part of the dynamics (\thmref{thm_superatt}).  

\parag{Leaf spaces.}
It is convenient to generalize the notion of lamination slightly. 
For a given integer $k \ge 1$, we say a topological space $\LL$ is a \textit{leaf space} if it is a connected topological space that consists of path-connected components (``\textit{leaves}") homeomorphic to $k$-orbifolds. 
For example, the regular leaf space $\Rf$, the affine part $\Afn$, and the affine lamination $\Af$ of a rational map $f$ are leaf spaces with leaves isomorphic to Riemann orbifolds ($=$ complex 1-orbifolds). 
(Technically, we have to remove ``isolated leaves" from $\Af$ when $f$ is Chebyshev or Latt\`es. See \cite[\S 5.4]{LM} and \cite[\S 3.1.1]{KL}. We do not deal with these cases in this paper.)

\parag{Quasi-isometry and quasiconformal maps.}
Let $\LL_1$ and $\LL_2$ be leaf spaces with leaves isomorphic to hyperbolic 3-orbifolds (resp. Riemann orbifolds). Suppose that there is a homeomorphism $H:\LL_1 \to \LL_2$. We say $H$ is (leafwise) \textit{quasi-isometric (resp. quasiconformal)} if there is a uniform $K \ge 1$ such that $H$ sends each leaf in $\LL_1$ to a leaf in $\LL_2$ as a $K$-quasi-isometry (resp. $K$-quasiconformal map). 

If $\LL_1$ and $\LL_2$ support dynamics $F_1:\LL_1 \to \LL_1$ and $F_2:\LL_2 \to \LL_2$ which are conjugated by $H$, we say $\LL_1$ and $\LL_2$ are \textit{ quasi-isometrically (resp. quasiconformally) equivalent}.

The following proposition is given in the proof of \cite[Theorem 9.1, Lemma 9.3]{LM}:

\begin{prop}[Quasi-isometric extension]\label{prop_aff2hyp}
For rational maps $f$ and $g$, suppose that there exists a quasiconformal conjugacy $\hhatA: \Af \to \Ag$. Then it is lifted to their $\Hyp^3$-laminations $\hhatH: \Hf \to \Hg$ as a quasi-isometric conjugacy. Moreover, it descends to a quasi-isometric homeomorphism $\hhatM: \Mf \to \Mg$, which extends to $\hhatM: \Bar{\Mf} \to \Bar{\Mg}$.
\end{prop}

\parag{Perturbation of superattractings.}
Let $\Rat_d$ denote the space of rational maps of degree $d \ge 2$. For a rational map $f \in \Rat_d$, we consider a perturbation $\fe \to f~(\e \to 0)$ within $\Rat_d$ with respect to the topology defined by uniform convergence (equivalently, convergence of coefficients).

Here is the main theorem of this section:
\begin{thm}[Relaxing superattracting cycles]\label{thm_superatt}
Suppose $f$ has a superattracting cycle $O$. Then there exists a perturbation $\fe \to f$ with the following properties:
\begin{itemize}
\item $\fe \to f$ just relaxes $O$ to be an attracting cycle and preserve the topological dynamics except the immediate basin of $O$; and
\item The affine parts $\AAA_{\fe}^\n$ and $\Afn$ are quasiconformally equivalent. 
\end{itemize}
\end{thm} 

Assuming this theorem, we can prove:

\begin{thm}[Stability of hyperbolics]\label{thm_hyp_stability}
Let $f$ be a hyperbolic rational map. For any small perturbation $\fe \to f$, their affine, $\Hyp^3$-laminations are quasiconformally or quasi-isometrically equivalent. In particular, their Kleinian laminations are homeomorphic. 
\end{thm}

\begin{pf}
If $f$ has no superattracting cycle and $\e \ll 1$, it is known that we have a quasiconformal conjugacy between any perturbation $\fe$ and $f$. It naturally extends to a quasiconformal conjugacy on the affine parts $\AAA_{\fe}^\n$ and $\Afn$.

If $f$ has superattracting cycles, by \thmref{thm_superatt}, we can perturb $f$ within any small neighborhood to relax all of the superattractings. We can do the same for $\fe$ in the same neighborhood, and thus there exists a quasiconformal conjugacy on the affine parts $\AAA_{\fe}^\n$ and $\Afn$.

Since $\Afn$ and $\AAA_{\fe}^\n$ are quasiconformally equivalent, the affine laminations $\Af$ and $\AAA_{\fe}$ are also quasiconformally equivalent by \propref{prop_Afn_is_Af_if_para}. (For example, $\iota_f$ gives the conformal conjugacy between $\Afn$ and $\Af$.) By \propref{prop_aff2hyp}, this perturbation is also accompanied by quasi-isometries between their $\Hyp^3$- and quotient 3-laminations. 
Since the quasiconformal conjugacy between the affine laminations preserves the Fatou-Julia decomposition, the isometry between the quotient laminations extends to the conformal boundaries.
\QED
\end{pf}

Thus hyperbolic rational maps are stable in the sense of Lyubich-Minsky laminations. 
Here is an interesting question for people familiar with $J$-stability of holomorphic family of rational maps \cite{MSS}: \textit{Is $J$-stability of a rational map equivalent to quasiconformal stability of its affine lamination?}

In our quadratic case, we have:
\begin{cor}\label{cor_hyp_center}
Let $X \subset \C$ be a hyperbolic component of the Mandelbrot set. For any $c$ and $c'$ in $X$, their affine and hyperbolic 3-laminations are quasiconformally and quasi-isometrically equivalent respectively. Moreover, their quotient laminations are quasi-isometrically homeomorphic including their conformal boundaries. 
\end{cor}

Independently C.Cabrera \cite{Ca} proved the topological equivalence of affine laminations for $c$ and $c'$ as above. A remarkable fact is that he also proved the converse. By the corollary above, we can restate it as a rigidity theorem:

\begin{thm}[Rigidity by Cabrera]\label{thm_carlos}
If there exists an orientation preserving homeomorphism between the affine laminations of hyperbolic $f_c$ and $f_{c'}$, then $c$ and $c'$ are in the same hyperbolic component. Moreover, they are actually quasiconformally equivalent.
\end{thm}

See \cite{Ca} for the proof. 
It seems likely that we can replace ``affine laminations" by ``$\Hyp^3$-laminations". 
However, by the investigation in this paper, it seems unlikely that we can replace it by the ``quotient 3-laminations". 
(It is more likely that we can replace it by ``Kleinian laminations".)

\parag{Proof of \thmref{thm_superatt}.}
We construct such a perturbation by quasiconformal perturbation of $f$. By taking a M\"obius conjugation, we may assume that $0$ is a critical point belonging to $O$. Let $l$ be the period of $O$ and $\al$ be a point in $O$ such that $f(\al)=0$. Note that $\al$ may be also a critical point like $\infty$ of $f(z)=1/z^2$.

Let us fix a small $\lambda>0$ and set $\D_\lambda:=\skakko{w \in \C:|w|<\lambda}$. Then there exists a disk neighborhood $D$ of $0$ and the B\" ottcher coordinate $\Phi:D \to \D_\lambda$ such that $\Phi:z \mapsto w$ conformally conjugates $f^l$ on $D$ to $w \mapsto w^p$ for some $p \ge 2$. Let $D'$ be the component of $f^{-1}(D)$ containing $\al$. If necessary, we replace $\lambda$ by smaller one such that $D'-\skakko{\al}$ contains no critical point. Note that $\Phi \cc f: D' \to \D_\lambda$ semi-conjugates $f^l$ on $D'$ and $w \mapsto w^p$ on $\D_\lambda$. 

\begin{figure}[htbp]
\centering
\vspace{0cm}\hspace{0cm}
\includegraphics[width=0.57\textwidth]{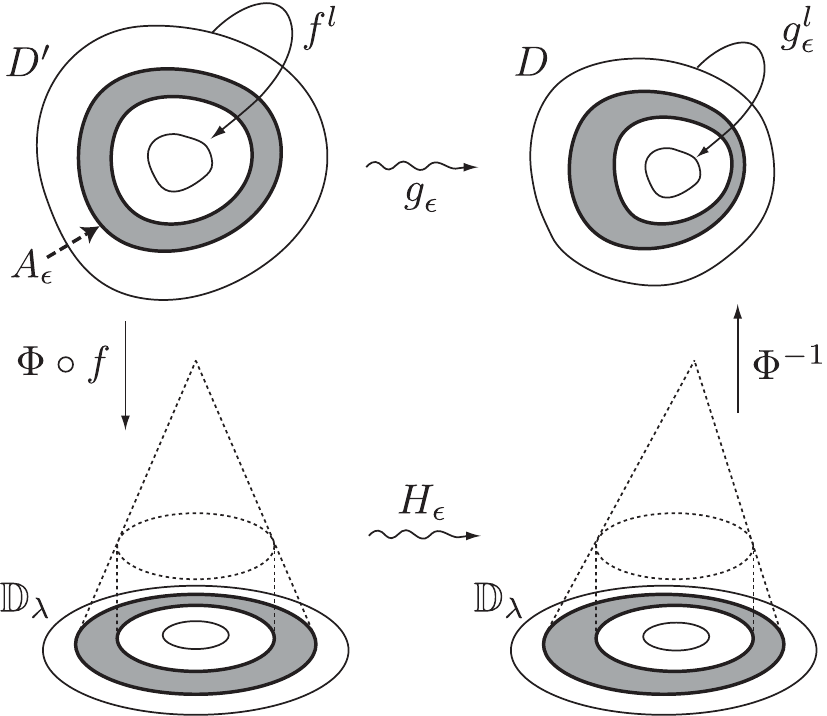}
\caption{This diagram explains the idea of the local perturbation. Note that $f$ and $\Phi \cc f$ in this diagram may have the degree two or more. To get the sliding map $\He$, consider a cone of height $1$ based on $\D_{\rho_2}$. Then slide the summit without changing the base. By truncating this cone at the level of $(\rho_2-\rho_1)/\rho_2$, we obtain a sliding map of a truncated cone. Via vertical projections and this sliding map, we obtain $\He$ above. }\label{fig_sliding}
\end{figure}

\parag{Local perturbation.}
Next we perturb $f:D' \to D$ as following: Fix $\rho_1$ and $\rho_2$ such that $\lambda^p<\rho_1 <\rho_2<\lambda$, and take any $\e>0$ satisfying $\rho_1 + \e < \rho_2$. Now we define a homeomorphism $\He: \D_\lambda \to \D_\lambda$ by 
$$
\He(w):=  \left\{ \begin{array}{rl}
w + \e~~~~~~~~~~~~~ & \text{if $|w| < \rho_1$}  \\
w + \e \cdot \dfrac{\rho_2-|w|}{\rho_2-\rho_1} 
  & \text{if $\rho_1 \le |w| \le \rho_2$}  \\
w  ~~~~~~~~~~~~~~~~~~& \text{if $ \rho_2< |w| < \lambda$}. 
          \end{array}
          \right.
$$
One can check that this map is actually a homeomorphism. See \figref{fig_sliding} for the idea. Note that the $L_\infty$-norm of the Beltrami differential of $\He$ is $O(\e)$.

\parag{Quasiregular map.}
Now we define a branched covering $\gep: \Cbar \to \Cbar$ by:
$$
\gep(z):=  \left\{ \begin{array}{rl}
\Phi^{-1} \cc \He \cc \Phi \cc f(z) &\text{if $z \in D'$}  \\
f(z)~~~~~~~~~~~~~~~~~~~ &\text{if $z \notin D'$}.  
          \end{array}
          \right.
$$
Then the following holds:
\begin{itemize}
\item $\gep \to f$ as $\e \to 0$, since $\He \to \id$. 
\item $\gep$ is quasiregular. Indeed, it is holomorphic except on a compact annulus $A_\e:=\skakko{z \in D': \rho_1 \le \Phi \cc f(z) \le \rho_2}.$ 
\item $\gep^l: D \to D$ is viewed as $w \mapsto w^p+\e$ through $\Phi$. Hence $0 \in D$ is not periodic of period $l$.  
\end{itemize}

\parag{The almost complex structure.}
Let $\sigma_0$ be the standard complex structure of $\Cbar$. We define a $\gep$-invariant almost complex structure $\sigma_\e$ by
$$
\sigma_\e :=  \left\{ \begin{array}{rl}
(\gep^{n})^\ast \sigma_0 & \text{on $\gep^{-n}(D), ~n \ge 0$}  \\
\sigma_0~~~~~~ & \text{otherwise}.  
          \end{array}
          \right.
$$
Since any orbit passes $A_\e$ (where $\gep^n$ gains dilatation) at most once, we can integrate $\sigma_\e$ and obtain a quasiconformal map $\he:\Cbar \to \Cbar$ fixing $0$, $\infty$, and some $b \in \Cstar$ with the following properties:
\begin{itemize}
\item $\he \to \id$ as $\e \to 0$, since the norm of the Beltrami differential tends to $0$.
\item $\he^\ast \sigma_0=\sigma_\e$ a.e.
\item $\he$ is conformal except on the attracting basin of $O$.
\item $\fe:=\he \cc \gep \cc \he^{-1}$ is a rational map which uniformly tends to $f$ as $\e \to 0$. 

\end{itemize}

\parag{Remark.} 
When $f(z)=z^2+c$, we may normalize the integrating map $\he$ such that $\he(z)/z \to 1$ as $|z| \to \infty$ by changing $b \in \Cstar$. Then we have $\fe \to f$ of the form $\fe(z)=z^2+c_\e$.

\parag{Lifting partial conjugacy.} 
Set $D_\e:=\he(D)$. Then we have a commuting diagram
$$
\begin{CD}
 \Cbar -\gep^{-1}(D) @>{\he}>> \Cbar-\fe^{-1}(D_\e) \\
 @V{\gep}VV                      @VV{\fe}V       \\
 \Cbar-D @>{\he}>> \Cbar-D_\e  
\end{CD}
$$
which actually is:
$$
\begin{CD}
 \Cbar -f^{-1}(D) @>{\he}>> \Cbar-\fe^{-1}(D_\e) \\
 @V{f}VV                      @VV{\fe}V       \\
 \Cbar-D @>{\he}>> \Cbar-D_\e.  
\end{CD}
$$
Since any backward orbits other than the (super)attracting cycles eventually get caught in this diagram, we naturally have a conjugacy $\hhat_\e$:
$$
\begin{CD}
 \Afn @>{\hhat_\e}>> \AAA_{\fe}^\n \\
 @V{\fhat}VV                    @VV{\fhat_\e}V       \\
 \Afn @>{\hhat_\e}>> \AAA_{\fe}^\n.  
\end{CD}
$$
More precisely, take any $\zhat=\bo{z} \in \Rf$. Then there exists an $N>0$ such that $z_{-n}$ stays in $\Cbar -D$ for all $n \ge N$. Set
$$
\hhat_\e:\zhat=\bo{z} \mapsto \fe^N(\he(z_N), \he(z_{-N-1}), \ldots).
$$ 
One can easily check that $\hhat_\e$ does not depend on the choice of $N$ and is a conjugacy between $\fhat$ and $\fhat_\e$ on their regular leaf spaces. Since $\he$ is quasiconformal, the lifted map $\hhat_\e$ is leafwise quasiconformal by the definition of the regular leaf spaces. It implies that leaves isomorphic to $\C$ are mapped to leaves isomorphic to $\C$. (Indeed, the map $\hhat_\e$ preserves all types of leaves in $\Rf$.) Thus the affine parts are preserved by this homeomorphism and we have the diagram above. 
\QED

\section{Analogy to quasi-Fuchsian deformations}\label{sec_04}
In this section we give an overview of the quotient (hyperbolic) 3-laminations as an analogue of the hyperbolic 3-manifolds associated with quasi-Fuchsian groups. 
In particular, we roughly explain the relation between the Mandelbrot set and the Bers slice (or the Bers simultaneous uniformization) in lamination context.
This gives an motivation to define upper and lower ends for our quotient 3-laminations.

\subsection{Quasi-Fuchsian deformation of hyperbolic 3-manifolds}
The aim of this paper is to observe topological changes of the quotient 3-lamination as the parameter $c$ of $f_c$ moves from one hyperbolic component to another via parabolic $g=f_\sigma$. 
Before proceeding to the investigations, we remark an analogy between quadratic maps and quasi-Fuchsian groups, which is already pointed out in \cite[\S 6.4, \S 10]{LM}. For the basics of quasi-Fuchsian groups as $\pi_1$ of hyperbolic $3$-manifolds, see \cite[\S 4.3.3]{MT} and \cite[\S 3.2]{McBook2}.

\parag{3-Manifolds for once-punctured torus groups.}
Let $\Gamma \subset \text{\textit{PSL}}(2,\R)$ be a Fuchsian group acting on the upper plane $\Hyp$. 
Then $\Gamma$ also stabilizes the lower half plane $\Hyp^-$ and the limit set $\Lam_\Gamma$.
For example, suppose that $S=\Hyp/\Gamma$ and $\overline{S}=\Hyp^-/\Gamma$ be once-punctured tori, where $\overline{S}$ is $S$ with its orientation reversed. 
In this case the limit set is $\R \cup \skakko{\infty}$, and its complement $\Omega_\Gamma$ (the region of discontinuity) is given by $\Hyp \sqcup \Hyp^-$. 
The Fuchsian group $\Gamma$ acts on $\Hyp^3$ properly discontinuously and the quotient hyperbolic 3-manifold $M=\Hyp^3/\Gamma$ is homeomorphic to $S \times (0,1)$. 
(Topologically it is a bit more natural to consider it a thickened surface.)
Set $\partial M:= \Omega_\Gamma/\Gamma =S \cup \Bar{S}$ and we call it the \textit{conformal boundary} of $M$. 
We say $\Bar{M}:=M \cup \partial M$ is the \textit{Kleinian manifold} of $M$. 
In this case, we have a natural product structure $\Bar{M}=S \times [0,1]$.
Let us call $S$ and $\overline{S}$ the \textit{upper} and \textit{lower ends}. 
It is well known that the deformation space of $\Gamma$ is identified as $T(S) \times T(\overline{S})$ (where $T(S)$ is the Teichm\"uller space of $S$ and isomorphic to $\Hyp$) which corresponds to the deformation of the upper end $S$ and the lower end $\overline{S}$. 

\begin{figure}[htbp]
\begin{center}
\includegraphics[width=.75\textwidth]{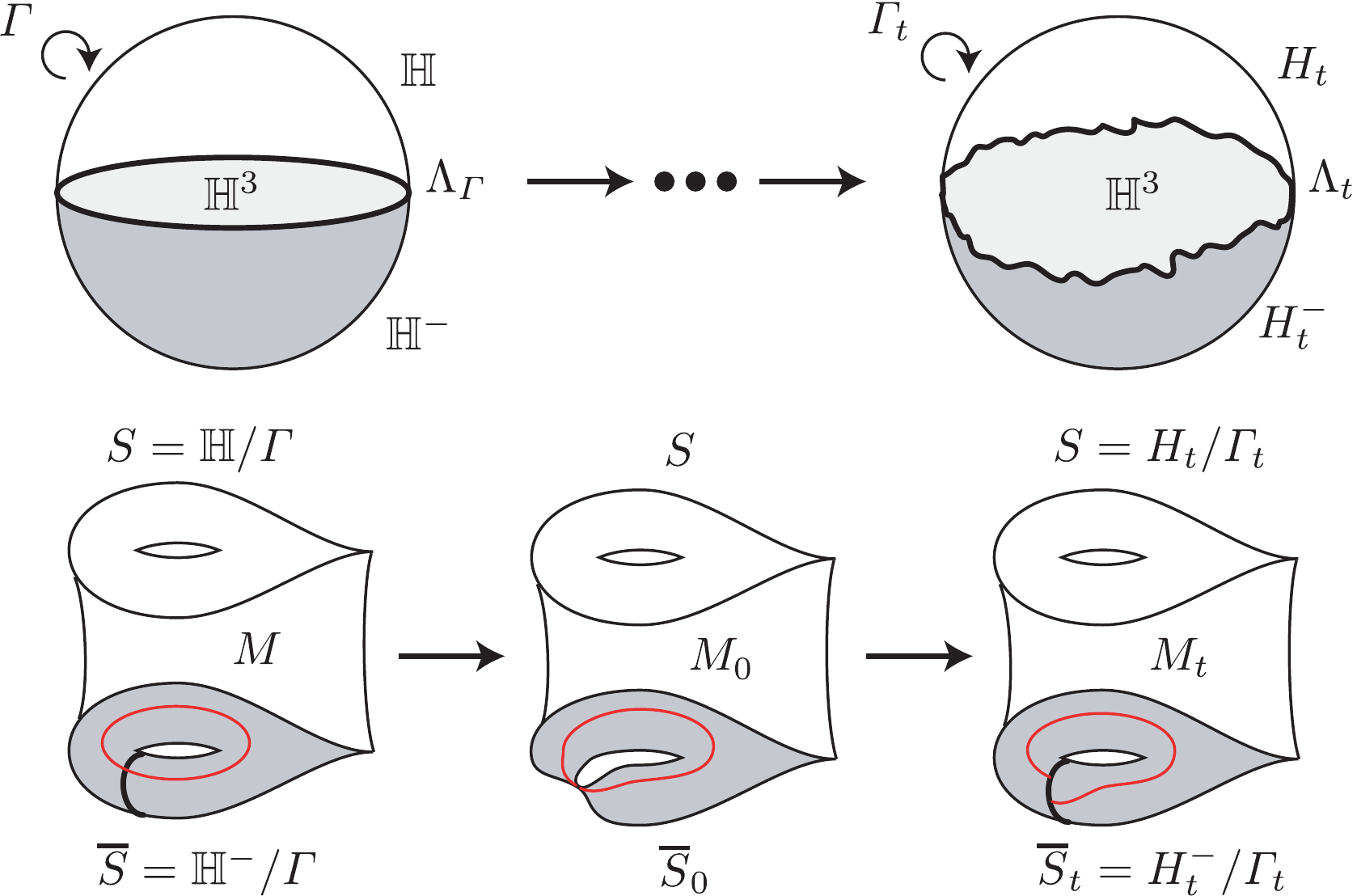}
\end{center}
\caption{Deformation of the hyperbolic 3-manifold associated with a once-punctured torus group. The heavy curves below show $\eta$ and $\eta_t$.}
\label{fig_fuchsian}
\end{figure}

\parag{Quasiconformal deformation.}
Now we consider a continuous family of quasiconformal maps $\skakko{\phi_t:\Cbar \to \Cbar}_{0 < |t| \le 1}$ such that: $\phi_{-1} = \id$; each $\phi_t$ is conformal on $\Hyp$; and the group $\Gamma_t=\phi_t \Gamma \phi_t^{-1}~(\Gamma_{-1}=\Gamma)$ is a Kleinian group which stabilizes $H_t=\phi_t(\Hyp)$ and $H_t^-=\phi_t(\Hyp^-)$ ($ =: $ \textit{a quasi-Fuchsian group}). 
Then the quotient $H_t/\Gamma_t$ is the same Riemann surface as $S$, but $\overline{S}_t:=H_t^-/\Gamma_t$ is quasiconformally deformed. 
(Such a deformation is realized by the Bers simultaneous uniformization \cite{Be}.)
In addition, we may assume that this deformation satisfies the following (see \figref{fig_fuchsian}) :
\begin{enumerate}
\item There exists a simple closed geodesic $\eta=\eta_{-1}$ on $\overline{S}=\overline{S}_{-1}$ such that $|t|$ is the hyperbolic length of the geodesic $\eta_t$ on $\overline{S}_t$ homotopic to $(\phi_t)_\ast(\eta)$.
\item By letting $t$ increase from $-1$ to $0$, we have an extra cusp at the limit. (Now $\overline{S}_0$ is a thrice-punctured sphere in this case.) 
\item Next as $t$ increases from $0$ to $1$, we plump the geodesic again but we twist it by a certain amount of length. More precisely, to have $\Bar{S}_t$, we fix a $\theta \in \R$ and we cut $\overline{S}_{-t}$ along $\eta_{-t}$ and twist by distance $t\theta$ to the right, and glue it again. 
\end{enumerate}

These pinching-and-plumping deformation eventually changes only the complex structure and the marking of the lower end $\overline{S}$. 
However, it is known that $M_t \approx S \times (0,1)$ for all $t$ (including $t = 0$ as a limit) and $\Bar{M}_t \approx S \times [0,1]$ for all $t \neq 0$. 
While this deformation process the main topological and combinatorial change is observed only in the lower ends.

\parag{The Bers slice vs. the Mandelbrot set.}
In the deformation space $T(S) \times T(\overline{S})$, the pinching-and-plumping deformation process above is realized in the \textit{Bers slice}, which is an embedded image of $\skakko{[S]} \times T(\overline{S}) \approx T(S)$ in a complex Banach space. 
One may find two distinct paths landing at the same boundary point in the Bers slice that correspond to pinching and plumping deformations.

In the following sections we will observe a similar phenomenon as $f_c$ moves from one hyperbolic component to another in the Mandelbrot set. 
(Not in the same component: This is a different point from the Bers slice.) 
Let us formalize our quadratic setting according to the quasi-Fuchsian setting. 
\begin{figure}[htbp]
\begin{center}
\includegraphics[width=.60\textwidth]{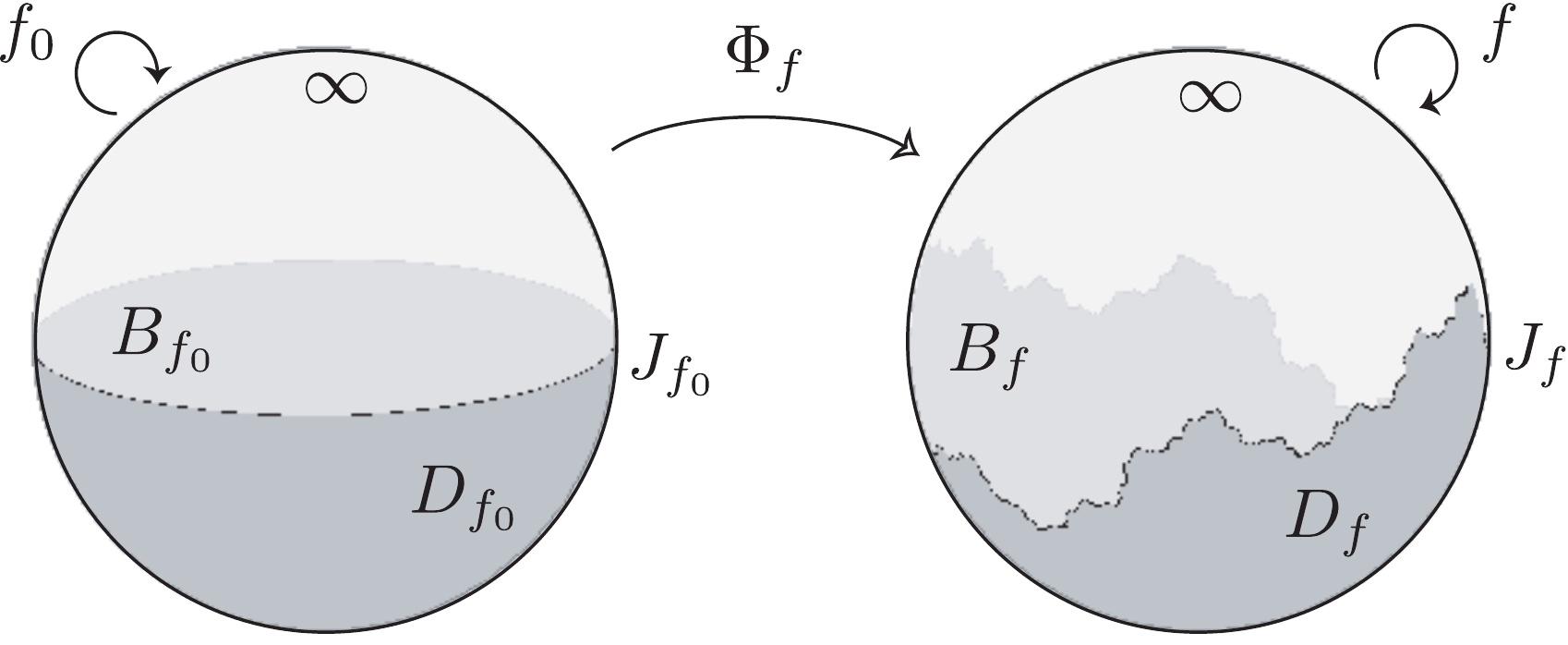}
\end{center}
\caption{The B\" ottcher coordinate $\Phi_f$ gives a conformal conjugacy between $f_0|_{B_{f_0}}$ and $f|_{B_f}$.}
\label{fig_basin_at_infty}
\end{figure}

We start with $f_0(z)=z^2$. This map stabilizes the sets $B_{f_0}:=\Cbar-\Dbar$, $J_{f_0}:=\skakko{|z|=1}$ and $D_{f_0}:=\D$. 
Let us try to regard the Mandelbrot set as the deformation space of this map. 
Set $f(z):=f_c(z)=z^2+c$ with $c$ in the Mandelbrot set. 
Then the dynamics of $f_c$ gives a natural division of the sphere
$$
\Cbar \ee B_f \sqcup D_f \sqcup J_f
$$
where $B_f:=\Cbar-K_f$, \textit{the basin at infinity}, and $D_f:=\Kfc$. 
Since the filled Julia set $K_f$ is connected, $B_f$ is a simply connected domain. 
Moreover, it is well-known that there exists a conformal conjugacy (the \textit{B\" ottcher coordinate}) $\Phi_f:B_{f_0} \to  B_f$ that conjugates $f_0|_{B_{f_0}}$ and $f|_{B_f}$ (\figref{fig_basin_at_infty}). 
On the other hand, the structure of $D_f$ may change drastically. (It can even be empty.)

For simplicity, let $c$ be a parameter in the main cardioid. Then the Julia set $J_f$ is a quasi-circle, so $D_f$ is a quasidisk. 
However, $c=0$ is the unique parameter with a superattracting fixed point (except $\infty$) so there is no quasifoncormal map that conjugates $f_0|_{D_{f_0}}$ and $f|_{D_f}$. 
This implies that we cannot regard the main cardioid as the direct analogue of the Bers slices contained in a quasiconformal deformation space of the holomorphic dynamics on the sphere. 

One solution to improve this situation is to consider the dynamics only near the Julia sets. 
(This is simply a problem of $J$-stability.)
Another solution we promote is to introduce the laminations.

\subsection{Deformation of the quotient 3-lamination of $\bs{f_0}$.}

\paragraph{Conformal boundary, upper and lower ends.}
A little more generally, suppose that $f=f_c$ has an attracting or parabolic cycle. 
Then $D_f$ is non-empty open set containing the critical point $z=0$.
Since $B_f$, $J_f$, and $D_f$ are completely invariant, it is natural to consider their lift to the affine lamination 
$$
\Af \ee \Bf \sqcup \Df \sqcup \Jf
$$
where $\Bf:=\pihat_f^{-1}(B_f)$, $\Df:=\pihat_f^{-1}(D_f)$, and $\Jf:=\pihat_f^{-1}(J_f)$. (Recall that $\pihat_f=\pi_f \cc \wp$. See the diagram in \secref{sec_02}.) 

By the B\" ottcher coordinate, the dynamics of $f: B_f \to B_f$ is conformally conjugate to $f_0(z)=z^2$ on $B_{f_0}=\C-\Dbar$. This conjugacy is naturally lifted to the natural extension. By restricting it on the affine part, or equivalently, on the affine lamination, we have a conformal conjugacy between the invariant sublaminations $\BB_{f_0}$ of $\AAA_{f_0}$ and $\Bf$ of $\Af$. 

The structure of $\BB_{f_0}$ is described as follows: (See also \cite[\S 11]{LM} and \cite{S}.)  
Under the identification of $\T=\R/\Z$ as $\partial \D=J_{f_0}$, the angle doubling $\delta: \T \to \T,~ \theta \mapsto 2 \theta$ is considered as a covering of degree two on the unit circle. 
Now the Julia set upstairs $\JJ_{f_0}$ is identified as the inverse limit $\TThat$ of the dynamics $\delta: \T \to \T$, which is a one-dimensional lamination. (See \S 6.1.)
Each leaf of $\TThat=\JJ_{f_0}$ is isomorphic to $\R$ with the natural affine structure coming from $\T$. We can attach $\Hyp$ to each leaf of $\TThat$, and this is exactly what $\BB_{f_0} \sqcup \JJ_{f_0}$ is. The action $\fhat_0:\BB_{f_0} \to \BB_{f_0}$ is leafwise isomorphic and properly discontinuous. The quotient $\BB_{f_0}/\fhat_0$ is called \textit{Sullivan's solenoidal Riemann surface lamination} $\SSS_0$, which is compact and consists of leaves isomorphic to $\Hyp$ or annuli. 

By conformal equivalence, the Riemann surface lamination $\Bf/\fhat$ is equivalent to $\SSS_0$ for any $f$. In addition, by the argument in \cite[\S\S 3.5]{LM} or \cite{S}, we have: 
\begin{prop}[Conformal boundary]\label{prop_Bf_Df}
Suppose $f=f_c$ has an attracting or parabolic cycle. Then both $\Bf$ and $\Df$ are Riemann surface sublaminations of $\Af$ and that the actions of $\fhat$ on $\Bf$ and $\Df$ are properly discontinuous. 
Thus the quotients $\Bf/\fhat$ and $\Df/\fhat$ are again Riemann surface laminations. 
In particular, $\Bf/\fhat$ is conformally the same as $\SSS_0:=\BB_{f_0}/\fhat_0$. 
\end{prop}
Hence $\Bf/\fhat$ and $\Df/\fhat$ form the two ends of $\Mf$ corresponding to the ends $S$ and $\overline{S_t}$ of $M_t$ above. 
We say $\SSS_0=\Bf/\fhat$ the \textit{upper end}, and $\Sf:=\Df/\fhat$ the \textit{lower end}.  
Now the union $\SSS_0 \cup \Sf$ is $\partial \Mf=\Ff/\fhat$, the \textit{conformal boundary} of $\Mf$ (\secref{sec_02}).
Recall that we say $\overline{\Mf} =\partial \Mf \cup \Mf$ is the Kleinian lamination associated with $f$.

\begin{figure}[htbp]
\begin{center}
\includegraphics[width=.70\textwidth]{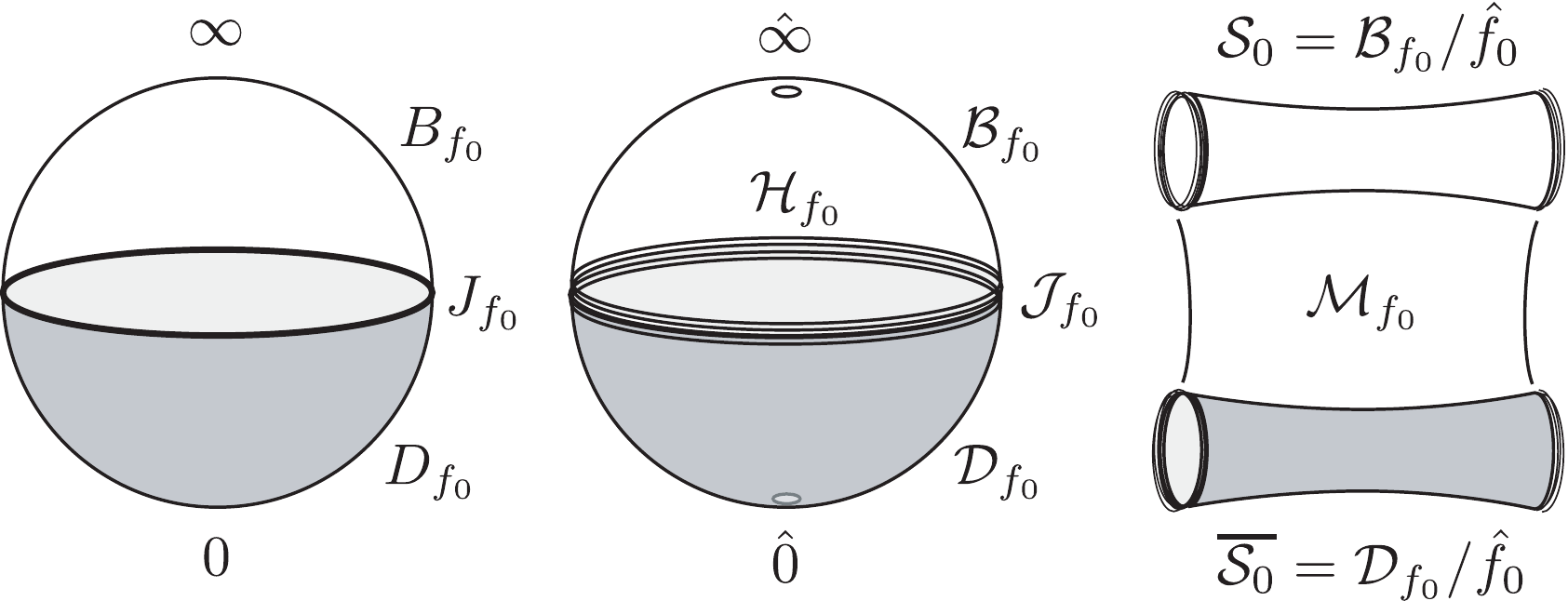}
\end{center}
\caption{The quotient 3-lamination of $f_0$.}
\label{fig_zz_lamins}
\end{figure}

\parag{Example: $\bs{c=0}$.}
Since $f_0$ has a symmetric dynamics with respect to $J_{f_0}$ (the unit circle), $D_{f_0}$ and $\DD_{f_0}$ are the mirror images of $B_{f_0}$ and $\BB_{f_0}$. 
Thus the lower end $\SSS_{f_0}=\DD_{f_0}/\fhat_0$ is the mirror image $\overline{\SSS_0}$ of $\SSS_0$, so we have a natural product structure $\overline{\MM_{f_0}} \approx \SSS_0 \times [0,1]$ (\figref{fig_zz_lamins}). 

\parag{Deformation of the Kleinian lamination $\bs{\Bar{\MM_{f_0}}}$.}
Now let us consider a deformation $\Mfbar$ of $\Bar{\MM_{f_0}}$ within the main cardioid. 
As an analogue to quasi-Fuchsian groups, we have:

\begin{thm}[Deformation in the main cardioid]
\label{thm_prod_str_of_M}
For any $f=f_c$ with $c$ in the main cardioid, the laminations $\AAA_{f}$, $\BB_{f}$, $\DD_{f}$, $\HH_{f}$, and $\MM_{f}$  are quasiconformally or quasi-isometrically equivalent to $\AAA_{f_0}$, $\BB_{f_0}$, $\DD_{f_0}$, $\HH_{f_0}$, and $\MM_{f_0}$ respectively. 

In particular, the Kleinian lamination $\overline{\Mf}=\Mf \cup \partial \Mf$ is homeomorphic to the product $\SSS_0 \times [0,1] \approx \Sf \times [0,1]$.
\end{thm}
This theorem emphasize the similarity between the Bers slice and the main cardioid of the Mandelbrot set. 
Note that the existence of the product structure in $\overline{\Mf}$ is already shown in \cite[\S 6.4]{LM}. Here we give an alternative proof:
\begin{pf}
The first paragraph directly comes from \corref{cor_hyp_center}. 
Then last sentence follows from the product structure $\overline{\MM_{f_0}} \approx \SSS_0 \times [0,1]$. 
\QED
\end{pf}

\parag{Remarks.}
\begin{itemize}\item
In \secref{sec_10} we will show that the quadratic Kleinian lamination $\Bar{\Mf}$ has a product structure if and only if $c$ of $f = f_c$ is in the main cardioid (\thmref{thm_prod_str_vs_main_cardioid}).

\item 
The theorem above is true even if we replace ``$f = f_c$ in the main cardioid" by ``$f$ in the hyperbolic component of $f_0(z) = z^d$ in the space of rational function of degree $d \ge 2$". 
For this case, the Riemann surface lamination $\Df$ is also dealt with in \cite[\S 10]{Mc2} as an analogue of the foliated unit tangent bundle over a Riemann surface. 
\end{itemize}

\parag{}
The following sections are devoted for the thorough study of the degeneration process of the affine, $\Hyp^3$-, and quotient 3-laminations and their conformal boundaries. 
According to the deformation process in the quasi-Fuchsian example above, the interior of the quotient 3-space would be topologically preserved, but its conformal boundary (the lower end) would change. 
We will also study the bifurcation of the degenerated laminations, which would not happen in the quasi-Fuchsian deformation space. 


\section{Degeneration pairs and Tessellation}\label{sec_05}
From this section we return to the setting of Part I \cite{Ka3} and investigate degenerations and bifurcations of quadratic maps in terms of the Lyubich-Minsky laminations.

\parag{Abstract of this section.}
We summerize the notation and results in Part I as follows:
\begin{itemize}
\item We formalize hyperbolic-to-parabolic degeneration of a quadratic map in the Mandelbrot set in terms of a \textit{degenetion pair} $(f \to g)$ with hyperbolic $f$ and parabolic $g$. 
We divide the degeneration pairs in two types: Case (a) and Case (b).
\item For hyperbolic $f$ above, we define \textit{degenerating arc system} $I_f$ which corresponds to the grand orbit $I_g$ of the parabolic cycle of $g$. 
\item We define \textit{tessellations} $\Tess(f)$ and $\Tess(g)$ for $(f \to g)$. Tessellations are tilings of the interiors of the filled Julia sets. Each tile is of the form $T=T(\theta, m, \ast)$ parameterized by \textit{angle} $\theta$, \textit{level} $m$, and \textit{signature} $\ast=\pm$ (\thmref{thm_tessellation}). 
\item The boundary of each tile is the union of three or four analytic arcs (``edges"). We classify them into \textit{equipotential, critical}, and \textit{degenerating edges}. 
\item By using tessellations, we have a pinching semiconjugacy $h:\Cbar \to \Cbar$ from $f$ to $g$ that just pinches $I_f$ to $I_g$ (\thmref{thm_semiconj}). 
\end{itemize}
The objects above will play key roles in the following sections. 
This section is a little more expository than the original \cite{Ka3}, where detailed definitions and proofs are given.

\subsection{Degeneration pairs: Case (a) vs. Case (b) (\cite[\S 2]{Ka3})}

\parag{Degeneration pair.}
Let $X$ be a hyperbolic component of the Mandelbrot set. By Douady and Hubbard's uniformization theorem \cite[Theorem 6.5]{Mi}, there exists a conformal map $\lam_X$ from $\D$ onto $X$ that parameterize the multiplier of the attracting cycle $O_f$ of $f=f_c$ for $c \in X$. Moreover, the map $\lam_X$ has a unique homeomorphic extension $\lam_X: \Dbar \to \bar{X}$ such that $\lam_X(e^{2 \pi i p/q})$ is a parabolic parameter for all $p, ~q \in \N$. A \textit{degeneration pair} $(f \to g)$ is a pair of hyperbolic $f=f_c$ and parabolic $g=f_\sigma$ where $(c, \sigma)=(\lam_X(re^{2 \pi i p/q}), \lam_X(e^{2 \pi i p/q}))$ for some $0<r<1$ and coprime $p, q \in \N$. By letting $r \to 1$, we have a convergence $f \to g$ uniformly on $\Cbar$.

\parag{Case (a) and Case (b).}
Let us introduce more notation that will be used throughout this paper. For a degeneration pair $(f \to g)$, let $O_f:=\skakko{\llist{\al}}$ (taking subscripts modulo $l$) be the attracting cycle of $f$ and $O_g:=\skakko{\beta_1, \ldots, \beta_{l'}}$ (taking subscripts modulo $l'$) be the parabolic cycle of $g$. Let $r e^{2 \pi i p/q}$ and $e^{2 \pi i p'/q'}$ denote the multipliers of $O_f$ and $O_g$ with coprime $p'$ and $q'$ in $\N$. (Then $O_g$ is a parabolic cycle with $q'$ repelling petals.) Then we have 

\begin{prop}[{{\cite[Proposition 2.1]{Ka3}}}]
\label{prop_case_a_case_b}
Any degeneration pair $(f \to g)$ satisfies either
\begin{itemize}
\item[] {\bf Case (a):} ~~$q=q'$ and $l=l'$; or
\item[] {\bf Case (b):} ~~$q=1<q'$ and $l=l'q'$.  
\end{itemize}
For both cases, we have $lq=l'q'=:\lbar$. 
\end{prop}

\begin{figure}[htbp]
\begin{center}
\includegraphics[width=.40\textwidth]{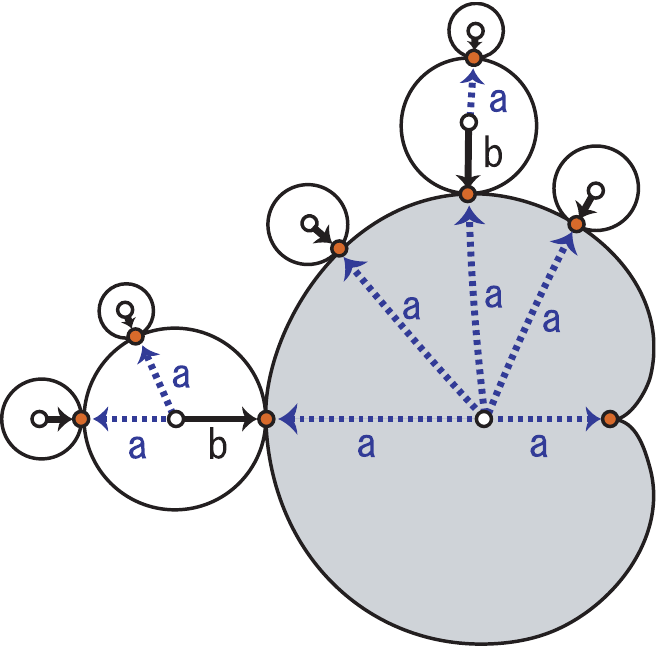}
\end{center}
\caption{Suppose that this is a part of a primitive copy of the Mandelbrot set. (For example, the cuspidal point on the boundary of the shadowed cardioid can be $c=1/4$, $c=-7/4$, etc.) The blue dotted arrows and black arrows show the Case (a) and Case (b) degeneration pairs respectively.}\label{fig_case_a_case_b}
\end{figure}

For example, if $f$ is in the main cardioid of the Mandelbrot set and $g$ is the root of a satellite limb, then $(f \to g)$ is Case (a). Indeed, any Case (a) degeneration pair is either a tuned copy of such $(f \to g)$; or $g$ corresponds to a cusped point of a small copy of the Mandelbrot set. (See \figref{fig_case_a_case_b}.) It would be helpful to remember that \textit{Case (a) happens iff the attracting points of $f$ and the parabolic points of $g$ have one-to-one correspondence.} 

\parag{Examples: The Cauliflowers.}
(In this paper we often identify the quadratic map $f_c(z)=z^2+c$ and the parameter $c$.)
The simplest example is given when $g=f_{1/4}$ and $f=f_c$ with $0 < c < 1/4$. 
In this case we have $p = p' = 1$, $q = q' = 1$, and $l = l' = 1$ thus Case (a). 
We call such a degeneration pair \textit{the Cauliflowers}.

\parag{Examples: The Rabbits.}
By {\it Douady's rabbit} we mean the quadratic map $f_{c_\mathrm{rab}}$ which is the center of the $1/3$-limb of the Mandelbrot set.
By {\it the Rabbits} we mean the simple deformation process from $f_0(z) = z^2$ to $f_{c_\mathrm{rab}}(z) = z^2 + c_\mathrm{rab}$ that can be represented by two degeneration pairs $(f_1 \to g)$ and $(f_2 \to g)$ as follows: Let $g=f_\s$ be the root of the $1/3$-limb. 
Then we can take a Case (a) degeneration pair $(f_1 \to g)$ with $f_1=f_{c_1}$ in the main cardioid and a Case (b) degeneration pair $(f_2 \to g)$ with $f_2=f_{c_2}$ in the same hyperbolic component as Douady's rabbit. 
(See also Figure 1 of Part I \cite{Ka3}.)

For the Case (a) Rabbits $(f_1 \to g)$, we have $p=p'=1$, $q= q'=3$, and $l= l'=1$. 
For the Case (b) Rabbits $(f_2 \to g)$, we also have $p=p'=1$ but $q= 1<q'=3$ and $l=3> l'=1$. 
(See Figures \ref{fig_deg_bif} and \ref{fig_deg_arc_sys} below.)
This represents a standard example of degeneration and bifurcation processes.
Note that $lq = l'q' = 3$ is invariant under this example.

\parag{Examples: The Airplanes.}
The only real $f_{c_\mathrm{air}}(z) = z^2 +  c_\mathrm{air}$ with superattracting cycle with period three is called \textit{the Airplane}. This parameter is the center of a primitive copy of the Mandelbrot set.
We say a degeneration pair $(f \to g) = (f_c \to f_\s)$ is \textit{the Airplanes} if $c_\mathrm{air} < c < -7/4$ and $\s = -7/4$. 
In this case we have $p = p' = 1$, $q = q' = 1$, and $l = l' = 3$ thus Case (a). 

\subsection{Degenerating arc system (\cite[\S 2]{Ka3})}

\parag{Degenerating arc system.}
For a degeneration pair $(f \to g)$ with $r \approx 1$, the parabolic cycle $O_g$ is approximated by an attracting or repelling cycle $O_f'$ with the same period $l'$ and multiplier $\lam' \approx e^{2 \pi i p'/q'}$. Let $\al_1'$ be the point in $O_f'$ with $\al_1' \to \beta_1$ as $r \to 1$.

In Case (a), the cycle $O_f'$ is attracting (thus $O_f'=O_f$) and there are $q'$ symmetrically arrayed repelling periodic points around $\al_1=\al_1'$. Then there exits an $f^{l'}$-invariant star-like graph $I(\al_1')$ that joins $\al_1'$ and the repelling periodic points by $q'$ arcs. 

In Case (b), the cycle $O_f'$ is repelling and there are $q'=l/l'$ symmetrically arrayed attracting periodic points around $\al_1'$. Then there exits an $f^{l'}$-invariant star-like graph $I(\al_1')$ that joins $\al_1'$ and the attracting periodic points by $q'$ arcs. 

\begin{figure}[htbp]
\begin{center}
\includegraphics[width=.60\textwidth]{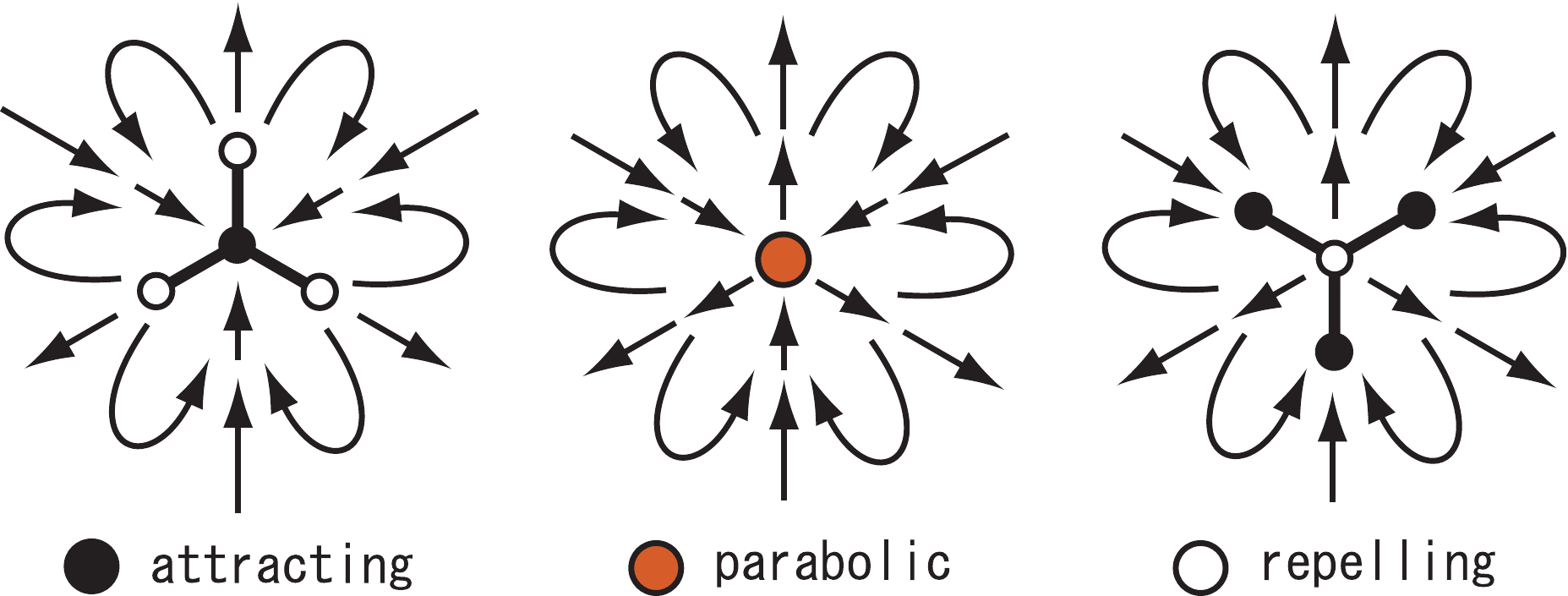}
\end{center}
\caption{Invariant graphs for the Rabbits. The arrows show the actions of $f_1^3, g^3$, and $f_2^3$.}
\label{fig_deg_bif}
\end{figure}

In both cases, the choice of $I(\al_1')$ has an ambiguity. Here we make a specific choice according to \cite[Lemma 2.3]{Ka3}: The arcs in $I(\al_1')$ and the critical orbit are lying on symmetrically arrayed radial rays in the K\"onigs coordinate. 
We define the \textit{degenerating arc system} $I_f$ by 
$$
I_f \dee \bigcup_{n \ge 0} f^{-n}(I(\al_1')).
$$
As its (expected) limit as $f \to g$, we define $I_g$ by 
$$
I_g \dee \bigcup_{n \ge 0} g^{-n}(\skakko{\beta_1}).
$$
Note that $I_f$ and $I_g$ are forward and backward invariant sets. For later use, let $\al_f$ be the set of all points which eventually land on $O_f$. For $\zeta \in \al_f$, let $I(\zeta)$ denote the component of $I_f$ (simply, $I_f$-component) that contains $\zeta$. In Case (b), each $I_f$-component contains $q'$ points in $\al_f$. 

\parag{Example: The Cauliflowers.}
Let $\al_1$ and $\gam_1$ be the attracting and repelling fixed points of $f$ respectively. Let $\beta_1 = 1/2$ be the parabolic fixed point of $g$. 
In this case we may regard $\al_1$ as $\al_1'$ and the invariant arc $I(\al_1')$ is just the interval $[\al_1, \gam_1]$ on the real axis. 
The degenerating arc system $I_f$ is the union of Jordan arcs that eventually land on the interval $[\al_1, \gam_1]$, and $I_g$ is the union of points that eventually land on $\beta_1$. 
(See \figref{fig_cauliflower}, the thick segments.)
The case of the Airplanes is similar.

\begin{figure}[htbp]
\centering
\includegraphics[width=0.8\textwidth]{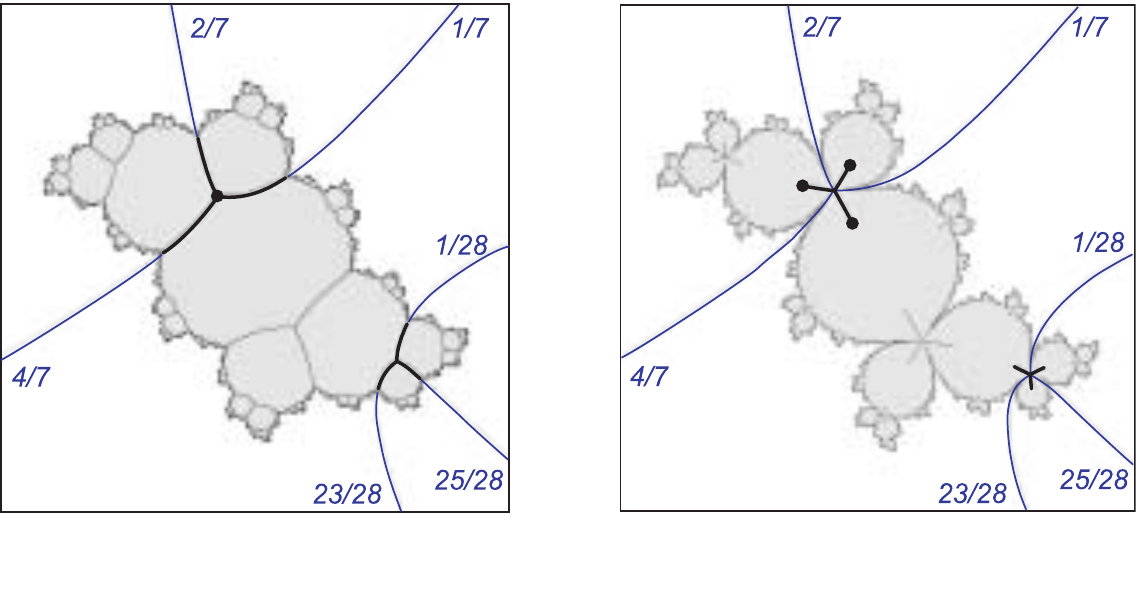}
\caption{Left, the filled Julia set of $f_1$ with their degenerating arc system roughly drawn in. Right, that for $f_2$. Attracting cycles are shown in heavy dots. Degenerating arcs with types $\skakko{1/7, 2/7/ 4/7}$ and $\skakko{1/28, 23/28, 25/28}$ are emphasized. See \S \ref{sec_05}.3.}
\label{fig_deg_arc_sys}
\end{figure}

\parag{Example: The Rabbits.}
The invariant arc $I(\al_1')$ with nearby dynamics for the Rabbits $f_1$ and $f_2$ are drawn in the left and right hand sides of \figref{fig_deg_bif} respectively. By taking its preimages, we have $I_{f_1}$ and $I_{f_2}$ as in \figref{fig_deg_arc_sys}.

\subsection{Tessellations associated with degeneration pairs (\cite[\S 3]{Ka3})}

\parag{Types.} 
Let $J_f$ and $K_f$ denote the Julia set and the filled Julia set of $f$. For any degeneration pair $(f \to g)$, we have $I_f \subset K_f$ and $I_g \subset J_g$. It is known that both $J_f$ and $J_g$ are locally connected so all the external rays land at these Julia sets.
Let $R_f(\theta)$ denote the external ray of angle $\theta$ and $\gam_f(\theta)$ denote its landing point in the Julia set.

For a point $z \in J_f$ or $J_g$, its \textit{type} $\Theta(z)$ is the set of angles of external rays that land on $z$. We abuse this term for more general subsets in the filled Julia sets $K_f$ and $K_g$. 
For any subset $E$ of the filled Julia set, its \textit{type} $\Theta(E)$ is the set of angles of external rays that land on points in $E$. Let $\Theta_f$ and $\Theta_g$ denote $\Theta(I_f)$ and $\Theta(I_g)$  respectively. 

\parag{Example: The Cauliflowers.} 
We again consider the case of the Cauliflowers. 
Since the repelling fixed point $\gam_1$ is $\gam_f(0)$ and the parabolic fixed point $\beta_1$ is $\gam_g(0)$, both $\Theta_f$ and $\Theta_g$ consist of the angles that eventually land on $0$ by the angle doubling map $\delta: \theta \mapsto 2\theta$. 
Hence $\Theta_{f}=\Theta_g$ and they are dense subsets of $\T=\R/\Z$. 

\parag{Example: The Rabbits.} 
For $(f_1 \to g)$, the set $\Theta_{f_1}$ consists of the angles that eventually land on $\skakko{1/7, 2/7,4/7}$ by angle doubling. 
One can also show that $\Theta_{f_1}=\Theta_g$ and they are dense subsets of $\T=\R/\Z$. 
We also have $\Theta_{f_2}=\Theta_{g}$ for $(f_2 \to g)$. 

\parag{Example: The Airplanes.}
We also have $\Theta_f = \Theta_g$ and this set consists of angles that eventually land on $\skakko{1/7, 2/7, 3/7, 4/7, 5/7, 6/7}$. See \cite[Figure 3]{Ka3}.

In general, we have:
\begin{prop}[{{\cite[Proposition 2.5]{Ka3}}}]
\label{prop_Prop2.5_of_Part1}
For any degeneration pair $(f \to g)$, the set $\Theta_f = \Theta(I_f)$ coincides with $\Theta_g = \Theta(I_g)$.
\end{prop}
We denote this set $\Theta_f=\Theta_g$ by $\Theta = \Theta_{(f \to g)}$.

\parag{Tessellation.}
The dynamics on/outside the Julia set is organized by the external rays and their landing points, satisfying
$$
f(R_f(\theta)) \ee R_{f}(2 \theta)
~~~\text{and}~~~
f(\gam_f(\theta)) \ee \gam_{f}(2 \theta).
$$
Tessellation is a method that organizes the dynamics inside the Julia set. Each tile is an element like a external ray outside:

\begin{thm}[Tessellation {{\cite[Theorem 1.1]{Ka3}}}]\label{thm_tessellation}
Let $(f \to g)$ be a degeneration pair. There exist families $\Tess(f)$ and $\Tess(g)$ of simply connected sets with the following properties:
\begin{itemize}
\item[\rm (1)] Each element of $\Tess(f)$ is called a {\rm tile} and identified by an {\rm angle} $\theta$ in $\Theta = \Theta_{(f \to g)}$, a {\rm level} $m$ in $\Z$, and a {\rm signature} $\ast=+$ or $-$.
\item[\rm (2)] Let $T_f(\theta, m, \ast)$ be such a tile in $\Tess(f)$. Then $f(\Tf(\theta, m, \ast))=\Tf(2\theta, m+1, \ast)$.
\item[\rm (3)] The interiors of tiles in $\Tess(f)$ are disjoint topological disks. Tiles with the same signature are univalently mapped each other by a branch of $f^{-i} \cc f^j$ for some $i,~j > 0$;
\item[\rm (4)] Let $\Pi_f(\theta, \ast)$ denote the union of tiles with angle $\theta$ and signature $\ast$. Then its interior $\Pi_f(\theta, \ast)^\cc$ is also a topological disk and its boundary contains the landing point $\gam_f(\theta)$ of $R_f(\theta)$. In particular, $f(\Pi_f(\theta, \ast))=\Pi_f(2\theta, \ast)$.
\end{itemize}
The properties above hold if we replace $f$ by $g$. Moreover:
\begin{itemize}
\item[\rm (5)] The union of tiles in $\Tess(f)$ is $\Kfc-I_f$. On the other hand, the union of tiles in $\Tess(g)$ is $\Kgc$.
\item[\rm (6)] The boundaries of $T_f(\theta, m, \ast)$ and $T_f(\theta', m', \ast')$ in $\Kfc -I_f$ intersect iff so do the boundaries of $T_g(\theta, m, \ast)$ and $T_g(\theta', m', \ast')$ in $\Kgc$.  
\end{itemize}
\end{thm}

One can easily understand the idea of tessellation by looking at the following figures (Figures \ref{fig_cauliflower} and \ref{fig_rabbit_tess}, \cite[Figure 2]{Ka3}, etc.), but the construction is not that easy.
Here we roughly explain the construction and the combinatorics of the addresses by showing some standard examples.

\parag{Example: The Cauliflowers.}
Their tessellations are shown in \figref{fig_cauliflower}. 
\begin{figure}[htbp]
\begin{center}
\includegraphics[width=.92\textwidth]{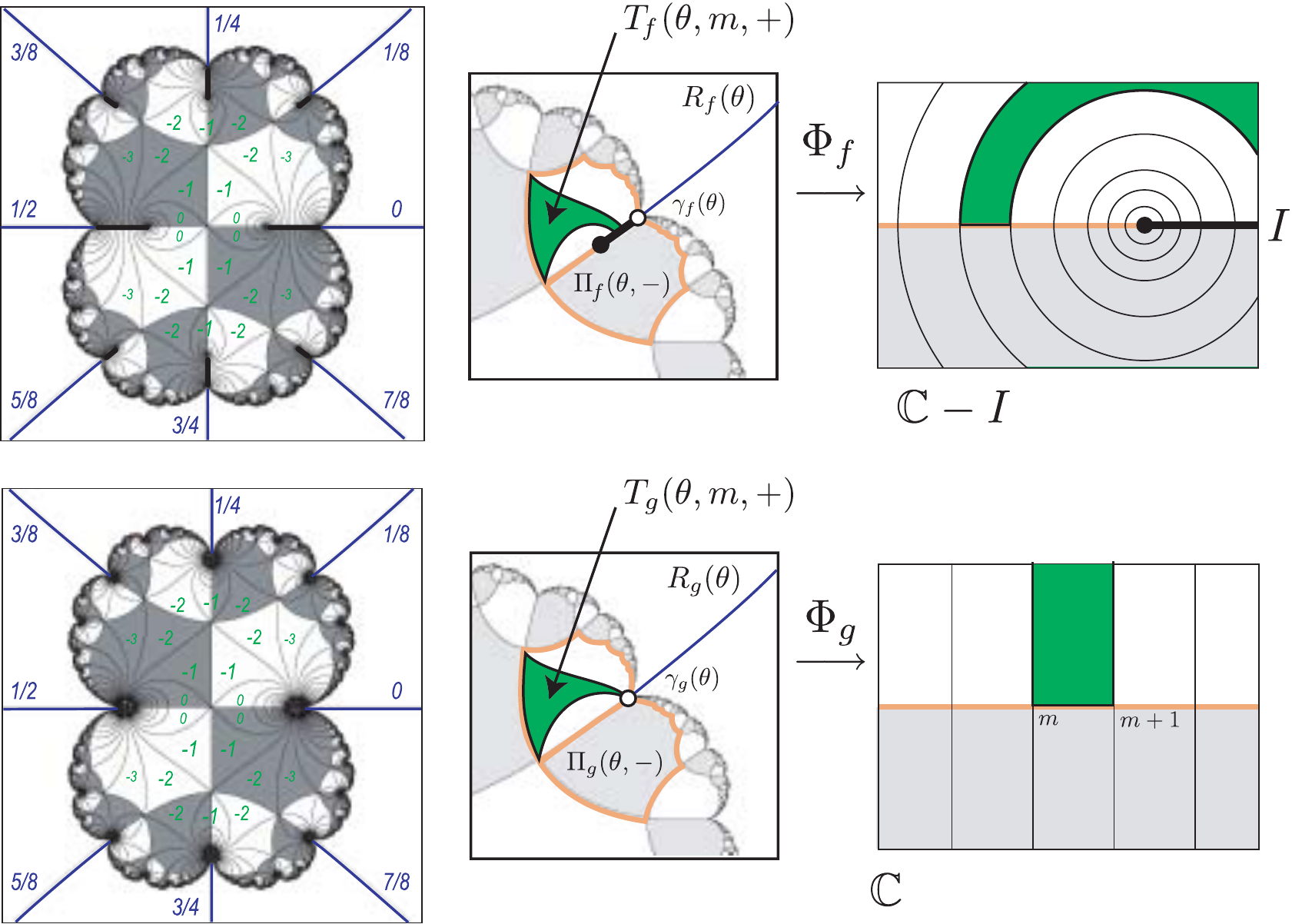}
\end{center}
\caption{Tessellations for the Cauliflowers and their construction. The blue curves with numbers are external rays with its angles. The numbers in green indicate the levels of tiles. The black thick curves show the degenerating arcs. The orange curves correspond to the critical edges of the tiles. See \figref{fig_edges}.}
\label{fig_cauliflower}
\end{figure}
The tiles are based on fundamental regions of the model dynamics given by $F:W \mapsto rW + 1~(0<r<1)$ and $G:W \mapsto W + 1$ on $\C$. 
The map $F$ has an attracting fixed point $a = 1/(1-r)$, and an invariant interval $I=[a, \infty) \subset \R$. 
An important fact is that $F: \C-I \to \C-I$ is topologically conjugate to $G: \C \to \C$. 
For example, we can find a conjugacy that sends the ``semi-fundamental region"
$$ 
A(m, \pm) \dee
\skakko{W \in \C -I ~:~ r^{m+1}a \le |W-a| \le r^{m}a, ~\pm \ip W \ge 0}
$$
to
$$
C(m, \pm) \dee
\skakko{W \in \C ~:~ m \le \rp W \le m+1, ~\pm \ip W \ge 0}
$$
for any $m \in \Z$ and any choice of $+$ or $-$.

We can relate these models and the original dynamics as follows. 
Now we may regard $r$ as the multiplier of the attracting cycle $\al_1$. The dynamics near $\al_1$ is conformally conjugate to $w \mapsto rw$ by a K\"onigs coordinate, but we take an extra conjugacy by an affine map so that it is conformally conjugate to $F$ near $a$. In fact, there exists a modified K\"onigs coordinate $\Phi_f: \Kfc \to \C$ such that: $\Phi_f(0)=0$; $\Phi_f$ is a branched covering semiconjugating $f|_{\Kfc}$ and $F$ on $\C$; the degenerating arc system $I_f$ is mapped onto the interval $I$; and we can take countably many univalent branches of $\Phi_f^{-1}$ over the upper and lower half-planes $\Hyp = \Hyp^+$ and $\Hyp^-$.

Let us choose specific branches $\Psi_\theta^\pm$ over $\Hyp^\pm$ so that the images of $\Hyp^\pm$ contain $\gam_f(\theta) \in I_f \subset J_f$ in their boundaries.
These images actually give the interiors of the panels $\Pi_f(\theta, +)$ or $\Pi_f(\theta, -)$. (Thus $\theta$ must be in $\Theta$.)

One can easily see that we have a homeomorphic extension $\Psi_\theta^\pm:\Bar{\Hyp^\pm} \to \C$.
Now the tile of address $(\theta, m , \pm)$ in $\Tess(f)$ is defined by
$$
T_f(\theta, m , \pm) \dee 
\Psi_\theta^\pm\kakko{A(m, \pm)}.
$$ 

Correspondingly, we have the normalized Fatou coordinate $\Phi_g: \Kgc \to \C$ so that $\Phi_g(0)=0$ and $\Phi_g$ is a branched covering semiconjugating $g|_{\Kgc}$ and $G(W)=W+1$ on $\C$. 
In particular, we can also take countably many univalent branches of $\Phi_g^{-1}$ over $\Hyp^\pm$. 
To define the tile $T_g(\theta, m , \pm)$, we take the branches $\Psi_\theta'^\pm$ over $\Hyp^\pm$ so that the images of $\Hyp^\pm$ contain $\gam_g(\theta) \in I_g \subset J_g$ in their closures.
Since we also have a homeomorphic extension $\Psi_\theta'^\pm:\Bar{\Hyp^\pm} \to \C$, we define the tile of address $(\theta, m , \pm)$ in $\Tess(g)$ by
$$
T_g(\theta, m , \pm) \dee 
\Psi_\theta'^\pm\kakko{C(m, \pm)}.
$$

\parag{Example: The Rabbits.}
\figref{fig_rabbit_tess} shows tessellations $\Tess(f_1)$ and $\Tess(f_2)$. 
It is difficult to draw a nice picture of the tessellation for $g$ because a terrible moire pattern appears. 
But one can imagine it by pinching $I_{f_1}$ or $I_{f_2}$. 
The construction of these tessellations are again based on the modified K\"onigs coordinates and the normalized Fatou coordinates. 
(Here we skip the details.)
\begin{figure}[htbp]
\begin{center}
\includegraphics[width=.9\textwidth]{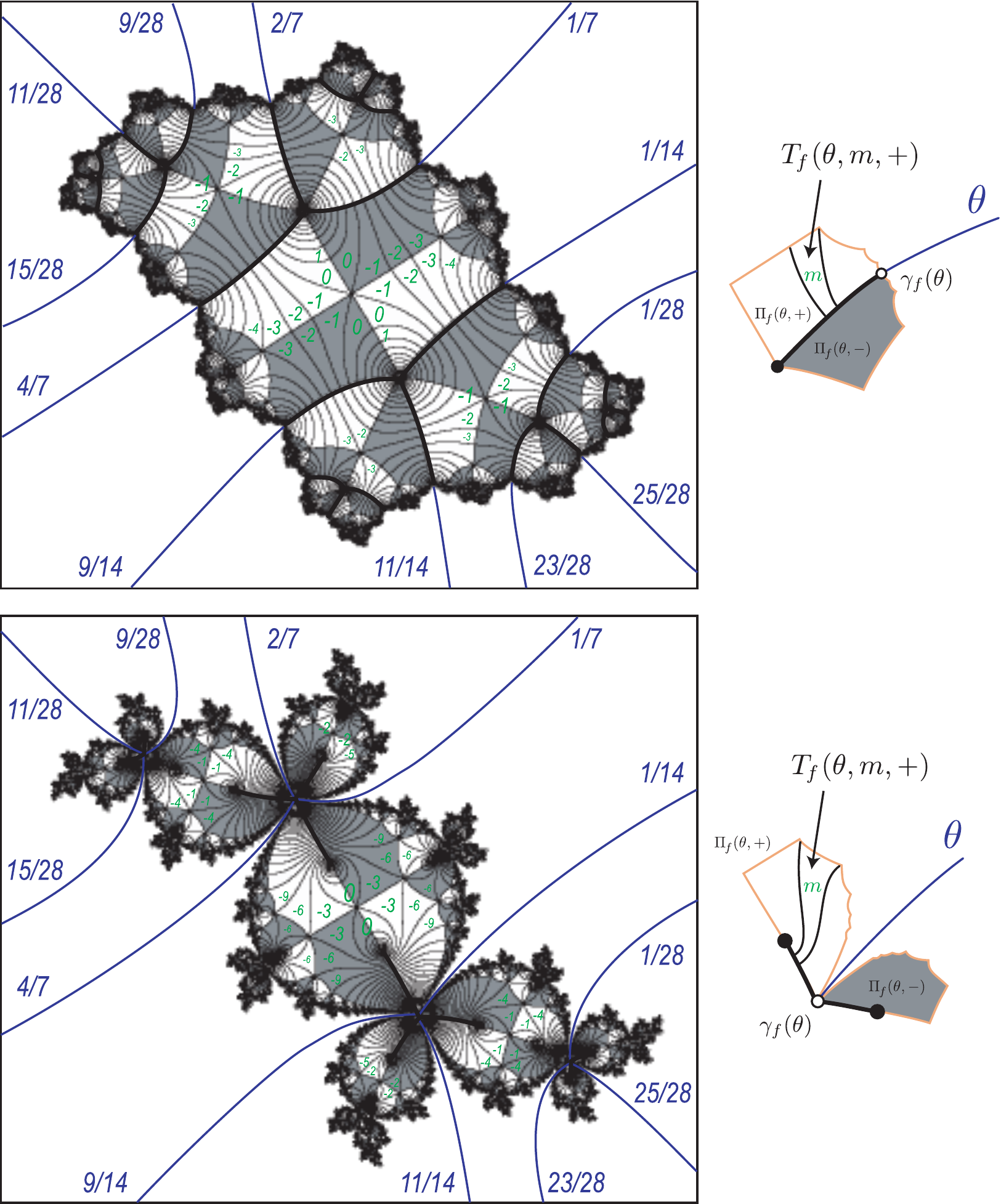}
\end{center}
\caption{Tessellations for the Rabbits $f_1$ and $f_2$. (cf. The Basilicas \cite[Figure 7]{Ka3}) 
By \thmref{thm_semiconj}, the tessellations for $g$ associated with $(f_1 \to g)$ and $(f_2 \to g)$ are given by just pinching $I_{f_1}$ and $I_{f_2}$. }
\label{fig_rabbit_tess}
\end{figure}

By looking at \figref{fig_rabbit_tess} carefully, one would find that the tessellations for $g$ associated with $(f_1 \to g)$ and $(f_2 \to g)$ would be different. 
In fact, the former is a subdivided tessellation of the latter. 
The difference comes from the difference of $p,~q$ and $l$ for $f_1$ and $f_2$. 
(We will deal with this difference in \secref{sec_10}.)

\begin{figure}[htbp]
\begin{center}
\includegraphics[width=.50\textwidth]{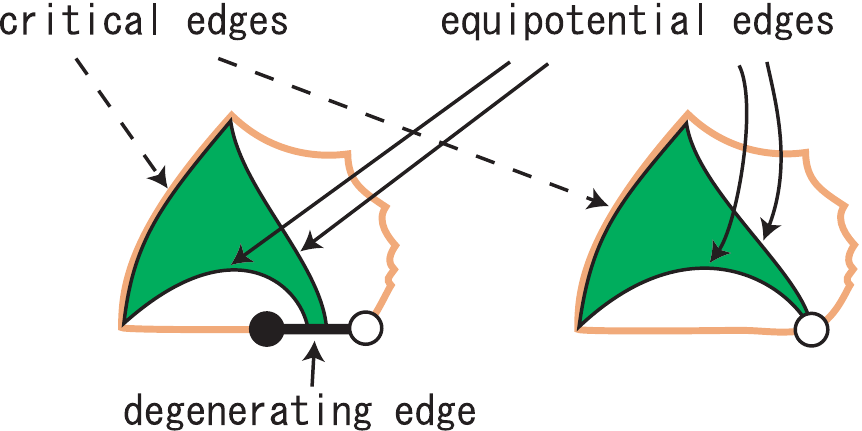}
\end{center}
\caption{Naming the edges of tiles.}
\label{fig_edges}
\end{figure}

\parag{Edges of tiles  (\cite[\S\S 3.2 - 3.3]{Ka3}).}
For the degeneration pair $(f \to g)$, every tile $T$ in $\Tess(f)$ (resp. $\Tess(g)$) is quadrilateral (reps. trilateral) with piecewise analytic boundary arcs, which we call \textit{edges}. 
It is useful to name these edges as follows (see \figref{fig_edges}): 
\begin{itemize}
\item {\bf Degeneration edges}: For each tile $T$ in $\Tess(f)$, its boundary $\partial T$ contains an arc shared with $I_f$. 
We call such an arc the \textit{degenerating edge} of $T$. 
\textit{Note that the degenerating edge is not contained in $T$ itself.} 
On the other hand, for each tile $T'$ in $\Tess(g)$, the part corresponding to the degenerating edge in its boundary $\partial T'$ is the vertex contained in $I_g$. \textit{Note also that this point is not contained in $T'$ itself.}
\item {\bf Equipotential edges}: For each $\partial T$ and $\partial T'$ there are two edges lying on the equipotential curves with respect to the linearizing (K\"onigs or Fatou) coordinate. When we make panels these edges are glued. We call them the \textit{equipotential edges} or the \textit{circular edges.}
\item {\bf Critical edges}: For each $\partial T$ and $\partial T'$ there is one edge left. Such an edge is contained in an arc joining two points in the grand orbit of the critical point $z=0$. 
We call such arcs the \textit{critical edges}. 
\end{itemize}

The combinatorics of edges would play an important role in Sections \ref{sec_09} and \ref{sec_10}. Here we give a refinement of \thmref{thm_tessellation}(6):

\begin{prop}[{{\cite[Proposition 3.2]{Ka3}}}]
\label{prop_edges}
The edges of tiles have the following properties:
\begin{itemize}
\item The equipotential (resp. critical) edges of $T_f(\theta, m, \ast)$ and $T_f(\theta', m', \ast')$ are shared iff so do the equipotential (resp. critical) edges of $T_g(\theta, m, \ast)$ and $T_g(\theta', m', \ast')$.
\item If the degeneration edges of $T_f(\theta, m, \ast)$ and $T_f(\theta', m', \ast')$ are shared, then $T_g(\theta, m, \ast)$ and $T_g(\theta', m', \ast')$ share their boundary points in $I_g$. That is, $\gam_g(\theta)=\gam_g(\theta')$.
\end{itemize}
\end{prop}

One can easily observe these properties in Figures \ref{fig_cauliflower} and \ref{fig_rabbit_tess}.

\subsection{Pinching semiconjugacy (\cite[\S 4 - 5]{Ka3})}

As an application of tessellation, we can show that there exists a pinching semiconjugacy from $f$ to $g$ for the degeneration pair $(f \to g)$:
\begin{thm}[Theorem 1.2 of \cite{Ka3}]\label{thm_semiconj}
Let $(f \to g)$ be a degeneration pair. There exists a semiconjugacy $h: \Cbar \to \Cbar$ from $f$ to $g$ such that:
\begin{enumerate}[{\rm (1)}]
\item $h$ only pinches each $I_f$-component to the point in $I_g$ with the same type. 
\item $h$ sends all possible $T_f(\theta, m, \ast)$ to $T_g(\theta, m, \ast)$, $R_f(\theta)$ to $R_g(\theta)$, and $\gam_f(\theta)$ to $\gam_g(\theta)$. 
\item $h$ tends to the identity as $f$ tends to $g$. 
\end{enumerate}
\end{thm}
So we have a precise picture for the degeneration of the dynamics. 
An interesting fact is that this theorem is proved \textit{without} using quasiconformal deformation. 

In the following sections we will investigate degeneration of the laminations by using the tessellations and this semiconjugacy.

\section{Degeneration of the affine laminations}\label{sec_06}

In this section we consider degeneration of the natural extentions, affine parts, and affine laminations for $(f \to g)$.

We first lift all the objects defined in the previous section to the natural extensions. 
For example, we will consider the external rays, degenerating arcs, and tessellations upstairs. 
(Such objects are nicely parameterized by the natural extention $\TThat$ of the angle doubling map $\delta: \T \to \T$.)

Next in \thmref{thm_natural_extension}, we construct a natural lift $\hhat:\Nf \to \Ng$ of the pinching semiconjugacy $h:\Chat \to \Chat$ given by \thmref{thm_semiconj}. 
This lifted pinching semiconjugacy $\hhat$ will give a precise picture of the degeneration of $\Nf$ to $\Ng$. 
By checking the action of $\hhat$ on the affine part $\Afn$, we also describe the degeneration of $\Afn$ to $\Agn$. 
Then we will conclude that as $f \to g$ the leafwise topology of $\Afn$ are preserved except on some specific leaves, called the \textit{principal leaves} (\thmref{thm_affine_part}). 
The descriptions on the affine parts are directly translated to the affine laminations  $\Af$ and $\Ag$ by \corref{cor_Afn_is_Af}. 

The results in this section are closely related to the backward Julia sets of the H\'enon mappings. See Appendix \ref{sec_11} for more details.

\subsection{Semiconjugacies in the natural extension}
First we describe the degeneration of the natural extensions associated with $(f \to g)$. Let us start with lifting the objects in $\Cbar$ (``downstairs") to the natural extensions $\Nf$ and $\Ng$ (``upstairs").

\parag{Natural extension of the angle doubling.}
To parameterize the objects we will deal with, let us consider the natural extension $\TThat$ of the angle doubling map $\delta: \T \to \T$, $\theta \mapsto 2\theta$. That is,
$$
\TThat \dee \{ \thetahat=(\theta_0,\theta_{-1}, \ldots)
:~\theta_0 \in \T,~2\theta_{-n-1}=\theta_{-n} \}
~\subset~\T \times \T \times \cdots 
$$
with a natural homeomorphic action $\delhat:\TThat \to \TThat$ given by
$$
\delhat(\thetahat)
\dee(2 \theta_0, \theta_0, \theta_{-1}, \ldots).
$$
We also call the elements of $\TThat$ \textit{angles}. 
Let $\Thetahat = \Thetahat_{(f \to g)}$ denote the invariant lift of $\Theta=\Theta_{(f \to g)}$ in $\TThat$. 
This subset of angles upstairs will play exactly the same role as $\Theta$ downstairs. 

\parag{External rays upstairs.}
External rays outside the Julia sets are naturally lifted to the natural extensions. For $\thetahat=(\theta_{-n})_{n \ge 0} \in \TThat$, we say $\zhat=(z_{-n})_{n \ge 0} \in  R_f(\thetahat)$ if $z_{-n} \in R_f(\theta_{-n})$ for all $n \ge 0$. That is, the external ray $R_f(\thetahat)$ is the backward orbits (or history) along the sequence 
$$
R_f(\theta_0) \leftarrow  R_f(\theta_{-1}) \leftarrow R_f(\theta_{-2}) \leftarrow \cdots.
$$
Note that any $\zhat \in R_f(\thetahat)$ is regular, thus $R_f(\thetahat) \subset \Rf$ ($= \Afn$ by \propref{prop_LM_4.5}). 
Since every external ray $R_f(\theta_{-n})$ lands on $J_f$, it is reasonable to denote the backward orbit $(\gam_f(\theta_0), \gam_f(\theta_{-1}), \ldots) \in \pi_f^{-1}(J_f)$ by $\gam_f(\thetahat)$. Hence the Julia set upstairs $\JJ_f^\n:=\pi_f^{-1}(J_f)$ is parameterized by $\TThat$. Indeed, one can check that $\gam_f: \TThat \to \JJ_f^\n$ is a continuous map. We define $R_g(\thetahat)$ and $\gam_g(\thetahat)$ in the same way.

\parag{Types upstairs.}
Let $\zhat=(z_{0}, z_{-1}, \ldots) \in \NN_f$ be a point in $\JJ_f^\n$. We define the \textit{type} $\Theta(\zhat)$ of $\zhat$ by the set of angles of the form $\thetahat = (\theta_{-n})_{n \ge 0} \in \TThat$ satisfying $\theta_{-n} \in \Theta(z_{-n})$ for all $n \ge 0$. One can easily check that $\delhat(\Theta(\zhat))=\Theta(\fhat(\zhat))$. (Note that the critical orbit of $f$ is bounded distance away from $J_f$. Thus $\card(\Theta(\zhat))$ is finite for each backward history $\zhat$ in $J_f$. See \cite[\S 2]{Ka3}.) Types for backward orbits in $J_g$ are also defined in the same way.

\parag{Degenerating arcs upstairs.}
Let us set $\alhat_f=\pi_f^{-1}(\al_f),~ \Ihat_f:= \pi_f^{-1}(I_f)$ and $\Ihat_g:= \pi_g^{-1}(I_g)$. For $\zetahat=(\zeta_{0}, \zeta_{-1}, \ldots) \in \alhat_f$, we define $I(\zetahat)$ by the set of backward orbits $\zhat=(z_{-n})_{n \ge 0}$ with $z_{-n} \in I(\zeta_{-n})$ for all $n \ge 0$. Equivalently, the set $I(\zetahat)$ is the path-connected component of $\Ihat_f$ (simply, \textit{$\Ihat_f$-component}) containing $\zetahat$ in $\NN_f$. The \textit{type} $\Theta(I(\zetahat))$ of $I(\zetahat)$ is the finite set of angles of the form $\thetahat = (\theta_{-n})_{n \ge 0} \in \Thetahat$ satisfying $\theta_{-n} \in \Theta(I(\zeta_{-n}))$ for all $n \ge 0$.

\parag{Tiles upstairs.}
The tessellations are naturally lifted to the natural extensions. For $\thetahat=(\theta_{-n})_{n \ge 0} \in \Thetahat$ and $\zhat=(z_{-n})_{n \ge 0} \in \Nf$, we say $\zhat \in T_f(\thetahat, m, \ast)$ if $T_f(\theta_0, m, \ast) \in \Tess(f)$ and $z_{-n} \in T_f(\theta_{-n}, m-n, \ast)$ for all $n \ge 0$. 
Since tiles in $\Tess(f)$ are homeomorphically mapped each other, the \textit{tile} $T_f(\thetahat, m, \ast)$ is contained in $\Rf= \Afn$ and the projection $\pi_f$ on $T_f(\thetahat, m, \ast)$ is a homeomorphism onto $T_f(\theta_0, m, \ast)$.
We gather all such tiles and denote it by $\Tess^\n(f)$. One can easily check that $\Tess^\n(f)$ actually tessellates $\pi^{-1}_f(\Kfc-I_f)$. 
For $g$, we define tiles and tessellation upstairs in the same way.

\parag{Pinching semiconjugacy upstairs.}
Now we have a natural lift of \thmref{thm_semiconj}:

\begin{thm}[Natural extensions] \label{thm_natural_extension} 
For degeneration pair $(f \to g)$, there exists a semiconjugacy $\hhat=\hhat_\n:\Nf \to \Ng$ from $\fhat$ to $\ghat$ with the following properties:
\begin{enumerate}[\rm (1)] 
\item $\hhat$ only pinches each $\Ihat_f$-component to the point in $\Ihat_g$ with the same type. 
\item $\hhat$ sends all possible $T_f(\thetahat, m, \ast)$ to $T_g(\thetahat, m, \ast)$, $R_f(\thetahat)$ to $R_g(\thetahat)$, and $\gam_f(\thetahat)$ to $\gam_g(\thetahat)$. 
\item $\hhat$ tends to the identity as $f$ tends to $g$. More precisely, the distance $\hat{d}(\hhat(\zhat), \zhat)$ tends to $0$ uniformly for any $\zhat \in \Nf$ where $\hat{d}(\cdot, \cdot)$ is a distance function on $\Cbar \times \Cbar \times \cdots$.
\end{enumerate}
\end{thm}
We will use this lifted semiconjugacy to investigate the degenerations of the laminations.

\begin{pf}
For $\zhat=(z_0, z_{-1}, \ldots) \in \Nf$, set
$$
\hhat(\zhat) \dee (h(z_0), h(z_{-1}), \ldots) ~\in~ \Ng.
$$ 
Since $h$ is a semiconjugacy from $f$ to $g$, one can easily check that $\hhat$ is surjective, continuous, and satisfies $\hhat \circ \fhat=\ghat \circ \hhat$. Thus $\hhat$ is a semiconjugacy from $\fhat$ to $\ghat$ on their respective natural extensions. 

Now properties (1)--(3) are almost straightforward. To show (3), one can choose the distance as follows: For $\zhat=(z_{-n})_{n \ge 0}$ and $\zetahat=(\zeta_{-n})_{n \ge 0}$ in $\Cbar \times \Cbar \times \cdots$, set $\hat{d}(\zhat, \zetahat):=\sum_{n \ge 0}d_{\Cbar}(z_{-n}, \zeta_{-n})/2^n$, where $d_{\Cbar}(\cdot, \cdot)$ is the spherical distance on $\Cbar$. Then $\hat{d}(\cdot, \cdot)$ determines the same topology as the natural topology of $\Cbar \times \Cbar \times \cdots$.
\QED
\end{pf}

\parag{Remark.}
Note that property (3) of the theorem implies $\Nf \to \Ng$ as $f \to g$ in the Hausdorff topology of non-empty compact subsets in $\Cbar \times \Cbar \times \cdots$. Indeed, for any uniform convergence of rational maps $f_n \to f$, it is always true that $\NN_{f_n} \to \Nf$ in the Hausdorff topology. 

\subsection{Principal leaves}
Next we consider degeneration of the affine part $\Afn$ to $\Agn$. 
To state the main theorem (\thmref{thm_affine_part}), we need to define some key objects: The \textit{irregular cycles} $\Ohat_f$ and $\Ohat_g$, \textit{principal leaves} $\Lam_f$ and $\Lam_g$, and the {\it cyclic degenerating arcs} $I(\Ohat_f)$.


\parag{Irregular cycles.}
For the attracting cycle $O_f= \skakko{\llist{\al}}$ with $f(\al_j)=\al_{j+1}$ (taking subscripts modulo $l$), we define its cyclic lift $\Ohat_f=\skakko{\llist{\alhat}}$ by
\begin{align*}
\alhat_j&\dee(\al_j, \al_{j-1}, \al_{j-2}, \ldots) ~\in~\Nf
\end{align*}
for each $j$ modulo $l$. For the parabolic cycle $O_g=\skakko{\beta_1, \ldots, \beta_{l'}}$ with $g(\beta_{j'})=\beta_{j'+1}$, we define $\Ohat_g=\{\betahat_1, \ldots, \betahat_{l'}\}$ in the same way by taking subscripts modulo $l'$. 
In the case of $q=1<q'$ (``Case (b)" in \secref{sec_05}), the pinching semiconjugacy $h$ gives a natural one-to-one correspondence between $l$ attracting points in $O_f$ and $l=q'l'$ attracting directions associated with $O_g$. So we assume that $\al_{j'}, \al_{j'+l'}, \ldots, \al_{j'+(q'-1)l'}$ are pinched to $\beta_{j'}$ by $h$; or equivalently,
\begin{equation}
I(\al_{j'}) \ee I(\al_{j'+l'}) \ee \cdots \ee I(\al_{j'+(q'-1)l'})
\tag{\ref{sec_06}.1}
\end{equation}
and this $I_f$-component has the same type as $\beta_{j'}$. Note that we always have the relation $f(I(\al_{j'}))=I(\al_{j'+1})$, which corresponds to $g(\beta_{j'})=\beta_{j'+1}$ by the semiconjugacy $h$.

\parag{Remark.}
In both Case (a) and Case (b), it is economical to use the subscript $j'$ modulo $l'$. This implies that the properties of $g$ are dominant when we consider degeneration of the affine parts. 

\parag{Examples.}
In the Cauliflowers $(f \to g)$ have the irregular cycles $\Ohat_f$ and $\Ohat_g$ which are just the lifted attracting and parabolic cycles $\{ \alhat_1 \}$ and $\{ \betahat_1 \}$. 
The same is true for the Case (a) Rabbits $(f_1 \to g)$, but
the Case (b) Rabbits $(f_2 \to g)$ have more irregular points in $\Ohat_{f_2}$: We have 
$\Ohat_{f_2} = \skakko{\alhat_1, \alhat_2, \alhat_3}$, 
the cyclic lift of the attracting cycle $\skakko{\al_1, \al_2, \al_3}$ with $\pi_f(\alhat_j)=\al_j$. 
Note that the equality (\ref{sec_06}.1) corresponds to $I(\al_1) = I(\al_2) =I(\al_3)$.  

\parag{The affine parts as the affine laminations.}
For our degeneration pair $(f \to g)$, it is fairly easy to specify the affine parts $\Afn$ and $\Agn$:

\begin{prop}\label{prop_Afn_Agn}
The affine parts of $(f \to g)$ are given by
$$
\Afn \ee \Nf-\hat{O}_f \cup \{\hat{\infty}\}  
\text{~~and~~}
\Agn \ee \Ng-\hat{O}_g \cup \{\hat{\infty}\}.
$$
Moreover, $\Afn$ and $\Agn$ are Riemann surface laminations with leaves isomorphic to $\C$, which are conformally equivalent to the affine laminations $\Af$ and $\Ag$ in the universal setting. 
\end{prop}
\begin{pf}
For $f$ and $g$, there is a common irregular point $\hat{\infty}=(\infty,\infty, \ldots)$ in their natural extensions. By \propref{prop_LM_4.5}, we have the equalities as above. 

Recall that the affine parts $\Afn$ and $\Agn$ are the unions of all leaves in the regular leaf spaces isomorphic to $\C$. 
Moreover, we may identify $\Afn$ and $\Agn$ as the affine laminations $\Af$ and $\Ag$ by \corref{cor_Afn_is_Af}. 
More precisely, the map $\wp|_{\Af}=\iota_f^{-1}:\Af \to \Afn$ gives the conformal conjugacy between $\fhat_\a$ and $\fhat_\n|_{\Afn}$ thus $\Af$ is conformally equivalent to $\Afn$. The same holds for $g$. 
\QED
\end{pf}
By using this identification, we will translate the following investigation on the affine parts to the affine laminations.

\parag{Principal leaves.}
Let us introduce the important notion of \textit{principal leaves} where we will observe the main difference caused by the hyperbolic-to-parabolic degeneration.

For each $j'$ modulo $l'$, we define the following subsets in the affine parts:
\begin{align*}
\Lam(\alhat_{j'})&\dee \skakko{\zhat=(z_0, z_{-1}, \ldots) \in \Afn
                   : \text{$\skakko{z_{-nl'}}_{n \ge 0}$ accumulates on $I(\al_{j'})$}};~\text{and}\\
\Lam(\betahat_{j'})&\dee \skakko{\zhat=(z_0, z_{-1}, \ldots) \in \Agn
                   : \text{$\skakko{z_{-nl'}}_{n \ge 0}$ accumulates on $\beta_{j'}$}}.
\end{align*}
Note that they do not contain the irregular points $\alhat_{j'}$ and $\betahat_{j'}$ themselves. 
In addition, we can replace ``$I(\al_{j'})$" by ``the repelling points in $I(\al_{j'})$" in the the definition of $\Lam(\alhat_{j'})$. 

\begin{figure}[htbp]
\begin{center}
\includegraphics[width=.7\textwidth]{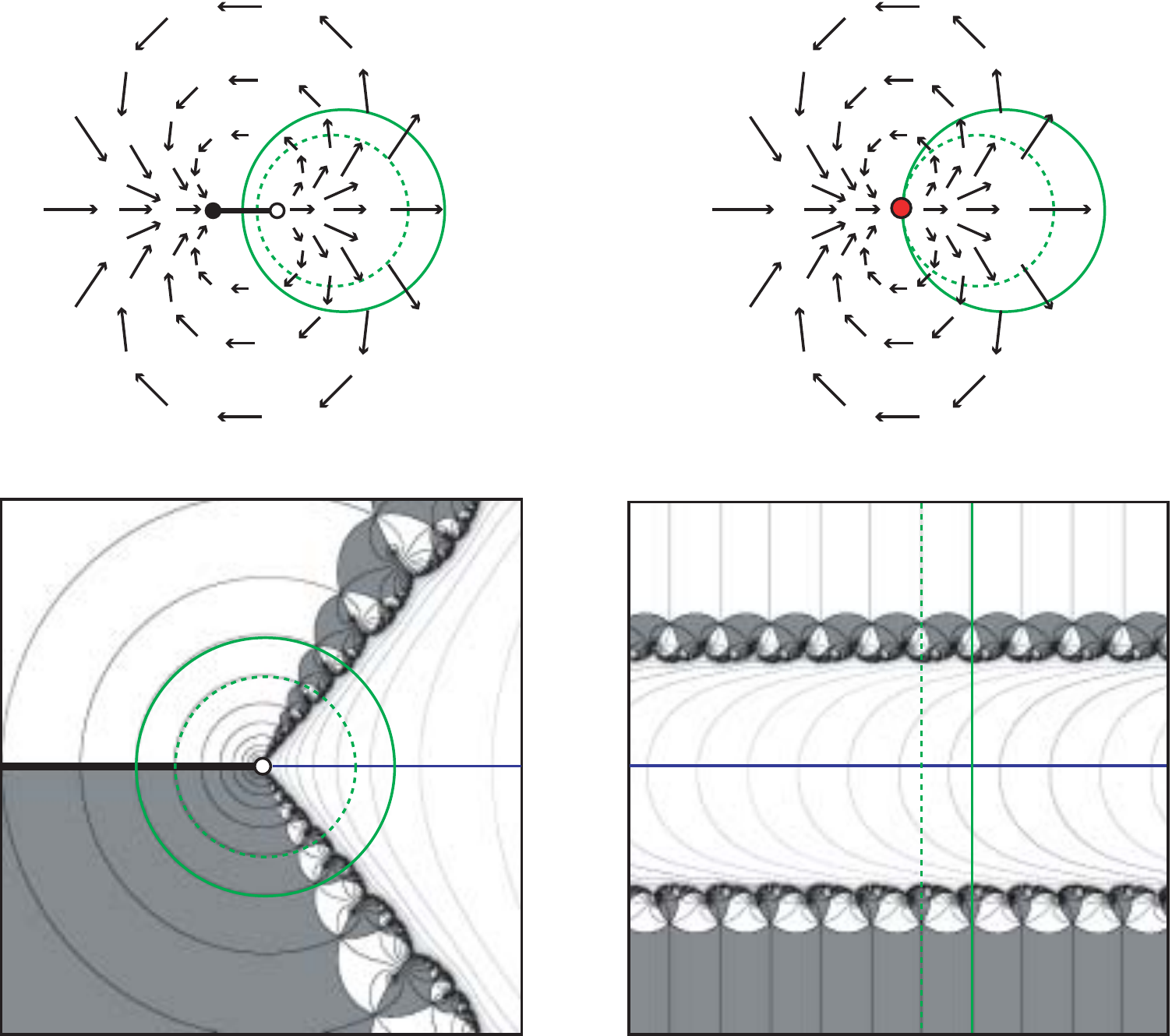}
\end{center}
\caption{The principal leaves for the Cauliflowers $(f \to g)$ below. The arrows above show the local dynamics near $\al_1$ and $\beta_1$.  The green curves show the fundamental regions. The black heavy line is $I(O_f)$, and the blue line is the external ray of angle $\thetahat=(0,0,\ldots)$.}
\label{fig_cauli_leaves}
\end{figure}

Then the dynamics downstairs implies:
\begin{prop}[Principal leaves]\label{prop_L}
The objects defined above have the following properties:
\begin{enumerate}[{\rm (1)}]
\item For each $j'$ modulo $l'$, $\fhat(\Lam(\alhat_{j'}))=\Lam(\alhat_{j'+1})$ and $\ghat(\Lam(\betahat_{j'}))=\Lam(\betahat_{j'+1})$. 
\item {\bf Case (a):} If $q=q'$ and $l=l'$, then for each $j'$ modulo $l'$, the set $\Lam(\alhat_{j'})$ consists of $q'$ leaves cyclic under $\fhat^{l'}$, and each of such leaves consists of the backward orbits converging to one of the $q'$ repelling periodic points in $I(\al_{j'})$. Similarly, the set $\Lam(\betahat_{j'})$ also consists of $q'$ leaves cyclic under $\ghat^{l'}$, and each of such leaves consists of the backward orbits converging to one of the $q'$ repelling petals of $\beta_{j'}$. 
\item {\bf Case (b):} If $q=1<q'$ and $l=l'q'$, then for each $j'$ modulo $l'$, the set $\Lam(\alhat_{j'})$ is just one leaf invariant under $\fhat^{l'}$ which consists of the backward orbits converging to the central repelling periodic point of $I(\al_{j'})$. Moreover, we have
$$ 
\Lam(\alhat_{j'}) \ee 
\Lam(\alhat_{j'+l'}) \ee \cdots \ee 
\Lam(\alhat_{j'+(q'-1)l'}).
$$
On the other hand, each $\Lam(\betahat_{j'})$ consists of $q'$ leaves described as above. 
\item In both cases above, we have $lq=l'q'=:\lbar$. Every leaf in $\Lam_f(\alhat_{j'})$ or $\Lam_g(\betahat_{j'})$ is invariant under $\fhat^\lbar$ or $\ghat^\lbar$ respectively. 
\item There exists a unique $\lam'$ with $|\lam'|>1$ such that for any $j'$ modulo $l'$ and any leaf $L$ of $\Lam_f(\alhat_{j'})$, the map $\fhat^{\lbar}:L \to L$ is conformally conjugate to the hyperbolic action $w \mapsto \lam' w$ on $\C$. Similarly, for any $j'$ modulo $l'$ and any leaf $L'$ of $\Lam_g(\betahat_{j'})$, the map $\ghat^{\lbar}:L' \to L'$ is conformally conjugate to the parabolic action $w \mapsto w+1$. Moreover, the multiplier $\lam'$ tends to $1$ as $f \to g$.
\end{enumerate}
\end{prop}
The proof is almost straightforward and left to the reader. (The last property is given by the existence of the linearizing coordinates near repelling cycles and petals.) To simplify the notation, we set 
$$
\Lam_f \dee 
\bigsqcup_{j'=1}^{l'} \Lam(\alhat_{j'})~~\text{and}~~ 
\Lam_g \dee \bigsqcup_{j'=1}^{l'} \Lam(\betahat_{j'}).
$$
We call $\Lam_f$ and $\Lam_g$ (and the leaves in these sets) the \textit{principal leaves}, where we will observe a significant degeneration associated with $(f \to g)$. 
The following example is helpful:

\parag{Example: The Cauliflowers.}
For the Cauliflowers $(f \to g)$, the principal leaves $\Lam_f$ and $\Lam_g$ are backward orbits converging to the repelling and parabolic fixed points respectively. They are invariant leaves of $\fhat$ and $\ghat$ with hyperbolic and parabolic dynamics respectively (\figref{fig_cauli_leaves}). Essentially the uniformizations of these principal leaves are given by uniformizing the quotient torus and the cylinder of the dynamics downstairs. See \cite{Ka1} and \cite[\S 4.2]{LM}. 

\begin{figure}[htbp]
\begin{center}
\includegraphics[width=.8\textwidth]{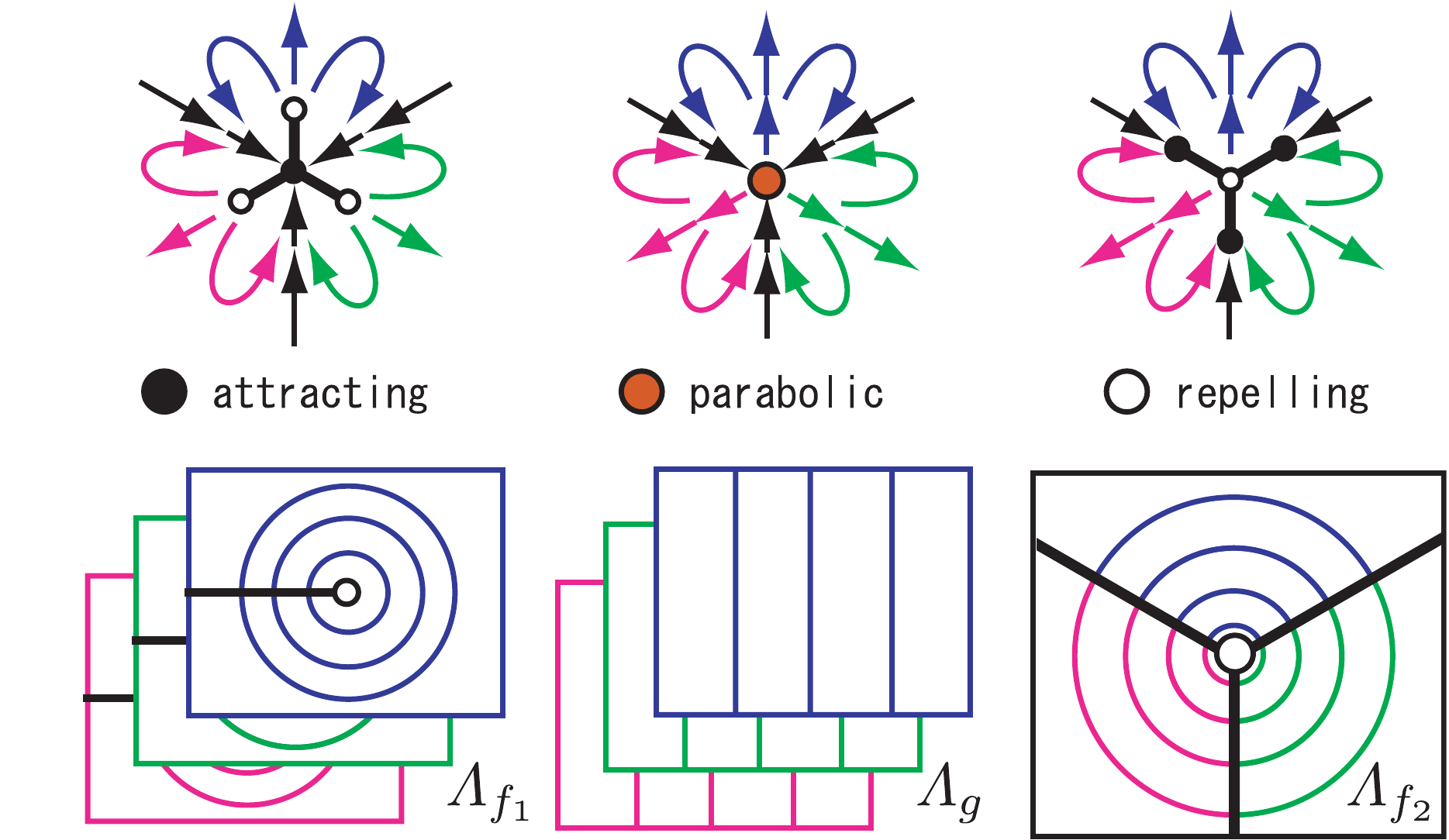}
\end{center}
\caption{The principal leaves of the Rabbits $f_1$, $g$, and $f_2$ (below). The three colors (blue, green, pink) indicate the three repelling directions of the parabolic cycle $\beta_1$, and the corresponding dynamics for $f_1$ and $f_2$. The black heavy lines indicate part of the cyclic degenerating arcs $I(\Ohat_{f_1})$ and $I(\Ohat_{f_2})$. The arrows show the actions of $f_1^3$, $g^3$, and $f_2^3$ downstairs.}
\label{fig_rabbit_leaves}
\end{figure}

\parag{Example: The Rabbits.}
For $f_1$($=$ Case (a)) and $g$ , both $\Lam_{f_1}$ and $\Lam_{g}$ consist of $q=q'=3$ leaves which are permuted by the actions of $\fhat_1$ and $\ghat$. 
On the other hand, for $f_2$($=$ Case (b)), $\Lam_{f_2}$ is just one $(=q)$ leaf which is invariant under the action of $\fhat_2$. 
However, we have obvious correspondence between the dynamics on $\Lam_{f_1}$, $\Lam_g$, and $\Lam_{f_2}$. (See \figref{fig_rabbit_leaves}.)

\parag{Cyclic degenerating arcs.}
As the parabolic cycles $O_g$ and $\Ohat_g$ are important elements of the dynamics of $g$ and $\ghat$, the corresponding objects
$$
I(O_f) \dee  h^{-1}(O_g)
~~\text{and}~~
I(\Ohat_f) \dee  \hhat^{-1}(\Ohat_g)
$$
play very important roles in the dynamics of $f$ and $\fhat$. We call them the {\it cyclic degenerating arcs} of $f$. Note that $I(\Ohat_f)$ upstairs contains irregular points $\Ohat_f$.

By \thmref{thm_natural_extension} and \propref{prop_L}, it is not difficult to check the following:

\begin{prop}[Cyclic degenerating arcs]\label{prop_cyclic_deg_arcs}
The cyclic degeneration arcs have the following properties:
\begin{enumerate}[\rm (1)]
\item Both $I(\Ohat_f)$ and $I(O_f)$ are decomposed into $l'$ connected components as follows:
$$
I(\Ohat_f) \ee  
\bigsqcup_{j'=1}^{l'} I(\alhat_{j'})~~\text{and}~~  
I(O_f) \ee \pi_f(I(\Ohat_f)) \ee \bigsqcup_{j'=1}^{l'} I(\al_{j'}),
$$
where $I(\alhat_{j'}) = \hhat^{-1}(\{\betahat_{j'}\})$ and $I(\al_{j'}) = h^{-1}(\skakko{\beta_{j'}})$. 
\item The set $I(\Ohat_f)$ (resp. $I(O_f)$) is the union of all the cyclic $\Ihat_f$-components (resp. $I_f$-components).
\item For each $j'$ modulo $l'$, we have  
$$
\Lam(\alhat_{j'}) \cap I(\hat{O}_f) \ee 
\Lam(\alhat_{j'}) \cap I(\alhat_{j'}) \ee 
I(\alhat_{j'})-\hat{O}_f.
$$
Thus we also have 
$$
\Lam_f-I(\Ohat_f) \ee 
\bigsqcup_{j'=1}^{l'} (\Lam(\alhat_{j'})-I(\alhat_{j'})).
$$ 
\end{enumerate}

\end{prop} 

\parag{Notation for the affine laminations.}
For later use, we also need to translate the sets defined above to those of the affine laminations in the universal setting.  
It is simply done by sending them by $\iota_f:\Afn \to \Af$ or $\iota_g:\Agn \to \Ag$. We set as follows: 
\begin{align*}
\text{Degenerating arc system:~}&\II_f \dee \iota_f(\Afn \cap \Ihat_f),~~\II_g \dee \iota_g(\Agn \cap \Ihat_g).\\
\text{Cyclic degenerating arcs:~}&\IO_f \dee \iota_f(\Afn \cap I(\Ohat_f)). \\
\text{Principal leaves:~}& 
\Lam^\a_f \dee \iota_f(\Lam_f),~~ 
\Lam^\a(\alhat_{j'}) \dee \iota_f(\Lam(\alhat_{j'})),~~\\
 &\Lam^\a_g \dee \iota_g(\Lam_g),~~ 
\Lam^\a(\betahat_{j'}) \dee \iota_g(\Lam(\betahat_{j'})).
\end{align*}
Note that we need to cut the irregular points when we translate $\Ihat_f$, $\Ihat_g$ and $I(\Ohat_f)$. 

Later we will use the fact that $\IO_f$ is the union of all cyclic components of the degenerating arc system $\II_f$ by \propref{prop_cyclic_deg_arcs}(2) above.  

\subsection{Degeneration of the affine parts/laminations.}

Now we start describing the degeneration of $\Afn$ and $\Af$. 
Let us check:

\begin{prop}[Semiconjugacy onto the affine part $\bs\Agn$.]
\label{prop_semiconj_onto_Agn}
The preimage of the affine part $\Agn$ by $\hhat:\Nf \to \Ng$ is $\Afn -I(\Ohat_f)$. In particular, the restriction $\hhat: \Afn - I(\Ohat_f) \to \Agn$ is a semiconjugacy. Hence we also have a semiconjugacy $\hhat_\a: \Af - \IO_f \to \Ag$ in the affine laminations.
\end{prop}

\begin{pf}
By \thmref{thm_natural_extension} we have $\hhat^{-1}(\skakko{\hat{\infty}})=\skakko{\hat{\infty}}$ and $\hhat^{-1}(\Ohat_g)=I(\Ohat_f)$. 
By \propref{prop_Afn_Agn} we have $\Agn \ee \Ng-\hat{O}_g \cup \{\hat{\infty}\}$ and thus $\hhat^{-1}(\Agn)=\Nf- I(\Ohat_f) \cup \{\hat{\infty}\} =\Afn -I(\Ohat)$.
Now the restriction $\hhat: \Afn - I(\Ohat_f) \to \Agn$ is a continuous and surjective map. Since the affine parts and $I(\Ohat_f)$ are invariant under the dynamics in the natural extensions, the restriction is a semiconjugacy.

Under the identification of the affine parts as the affine laminations of $(f \to g)$ (\propref{prop_Afn_Agn}), we have the corresponding semiconjugacy $\hhatA: \Af-\IO_f \to \Ag$.
\QED\end{pf}

\parag{Meekness.}
As we will see, the set $\Afn - I(\Ohat_f)$ in the affine part is an $\fhat$-invariant Riemann surface leaf space. (Recall the definition of \textit{leaf spaces} given in \secref{sec_03}.)
Thus $\hhat: \Afn - I(\Ohat_f) \to \Agn$ gives a semiconjugacy on the leaf spaces. 

To describe more precise properties of this semiconjugacy, we introduce some terminologies to classify surjective and continuous maps between leaf spaces.  
Let $\LL_1$ and $\LL_2$ be leaf spaces whose leaves are all Riemann surfaces or all hyperbolic 3-manifolds (which may be incomplete and may have boundaries). We say
\begin{itemize}
\item a map $H: \LL_1 \to \LL_2$ is {\it leafy} if $H$ is surjective, continuous, and for any leaf $L'$ in $\LL_2$, there exits a unique leaf $L$ in $\LL_1$ with $L=H^{-1}(L')$; 
\item a map $H: \LL_1 \to \LL_2$ is {\it meek} if $H$ is leafy, and for any leaf $L'$ in $\LL_2$, its preimage leaf $L=H^{-1}(L')$ is homeomorphic to $L'$; and 
\item a meek map $H: \LL_1 \to \LL_2$ is {\it very meek} if for any leaf $L'$ in $\LL_2$, there exists a leafwise quasiconformal or quasi-isometric homeomorphism between $L'$ and its preimage leaf $L=H^{-1}(L')$.
\end{itemize}
A leafy map gives one-to-one correspondence between the leaves in two leaf spaces. A meek map $H$ does not change the leafwise topology, but $H$ itself need not be a leafwise homeomorphism. 
For example, consider the simply connected Riemann surfaces $\Chat$, $\D$, and $\C$, which are all trivial leaf spaces with one leaf. 
Then the pinching semiconjugacy $h:\Chat \to \Chat$ given in \secref{sec_05} is very meek. 
Any homeomorphism $H:\D \to \C$ is meek, but not very meek since $H$ cannot be quasiconformal. 

\begin{figure}[htbp]
\begin{center}
\includegraphics[width=.55\textwidth]{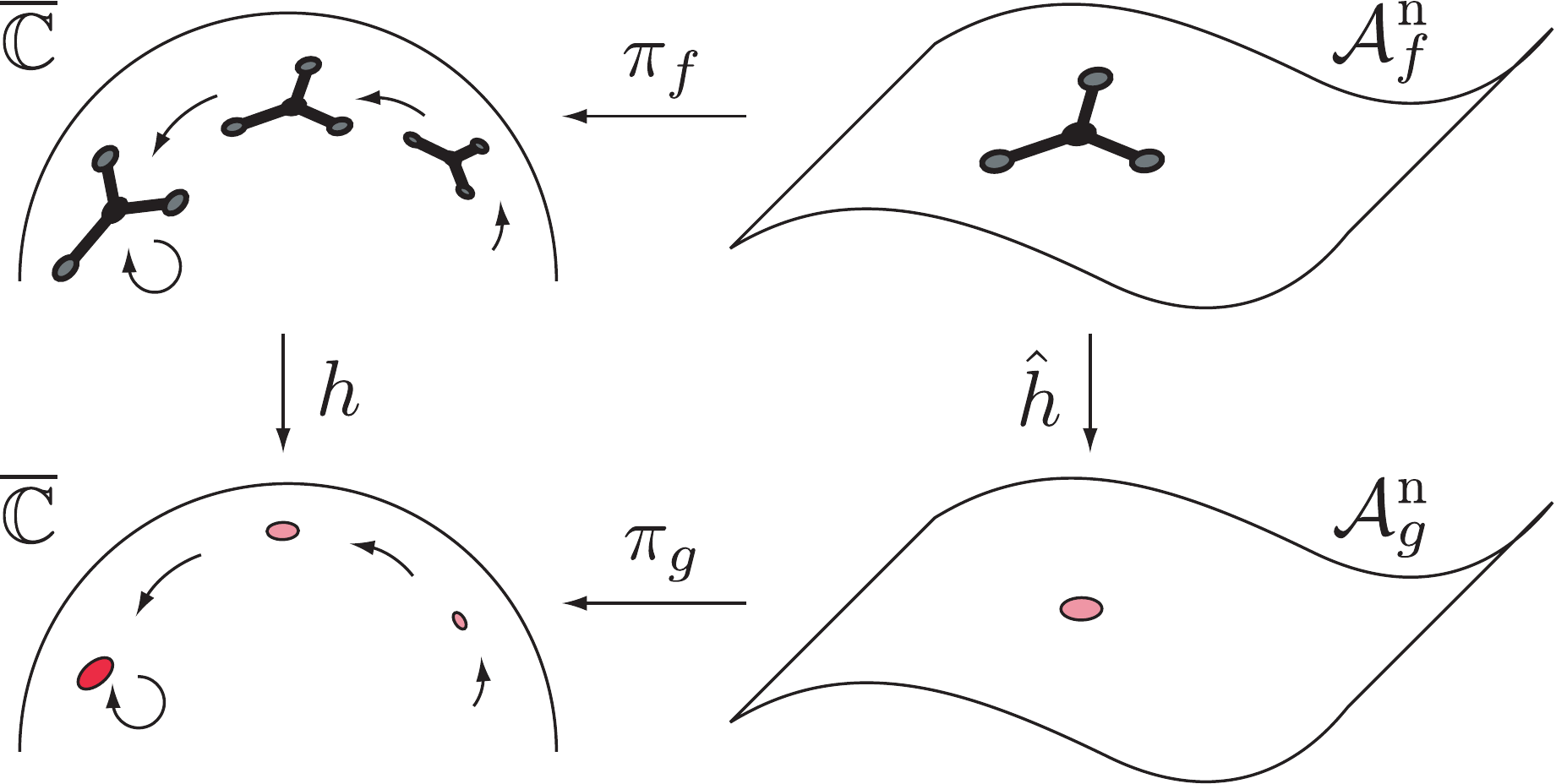}
\end{center}
\caption{Degeneration of a non-cyclic $\Ihat_f$-component.}\label{fig_deg_Af}
\end{figure}

Here is the main theorem in this section:
\begin{thm}[Affine part/laminations]\label{thm_affine_part}
For any degeneration pair $(f \to g)$, the set $\Afn-I(\Ohat_f)$ is a leaf space and the restriction $\hhat: \Afn-I(\Ohat_f) \to \Agn$ of $\hhat=\hhat_\n:\Nf \to \Ng$ is a meek pinching semiconjugacy. More precisely, 
\begin{enumerate}[{\rm (1)}] 
\item {\bf Non-principal part.} A further restriction $\hhat|~\Afn - \Lam_f \to \Agn - \Lam_g$ is a very meek map between the leaf spaces with leaves isomorphic to $\C$.
\item {\bf Principal part.} The complementary restriction $\hhat|~ \Lam_f - I(\Ohat_f)  \to \Lam_g$ is a meek map between the leaf spaces with finitely many leaves described as follows. For each $j'$ modulo $l'$, let $L_1', \ldots, L_{q'}'$ be the $q'$ principal leaves in $\Lam_g(\betahat_{j'})$. 
Then the set $\Lam(\alhat_{j'})-I(\alhat_{j'})$ has $q'$ path-connected components $S_1, \ldots, S_{q'}$ satisfying:
\begin{itemize}
\item Each $S_{i}~(1\le i \le q')$ is a simply connected Riemann surface that is invariant under $\fhat^{\lbar}$.
\item The preimage $\hhat^{-1}(L_{i}')$ is $S_i$. Thus $\skakko{S_i}$ and $\skakko{L_{i}'}$ has one-to-one correspondence.
\end{itemize}
\end{enumerate}
The corresponding semiconjugacy $\hhatA: \Af-\IO_f \to \Ag$ in the affine laminations also has the same properties.
\end{thm}
Property (1) means that any two non-principal leaves in $\Afn$ are never merged into one leaf by $\hhat$. 
One can easily check that every non-cyclic degenerating arc $I(\zetahat)$ containing $\zetahat \in \alhat_f-\hat{O}_f$ is compactly contained in a leaf of $\Afn$ with respect to the leafwise topology. Then $\hhat$ just pinches $I(\zetahat)$ to one point $\hhat(\zetahat) \in \Ihat_g$ with the same type in a leaf of $\Agn$ (\figref{fig_deg_Af}). Though each non-principal leaf contains infinitely many such $\Ihat_f$-components, the pinching map $\hhat$ preserves the topology of the leaves isomorphic to the complex plane. 
Thus the major topological change in the degeneration is described in (2).

\begin{figure}[htbp]
\centering{
\includegraphics[width=0.9\textwidth]{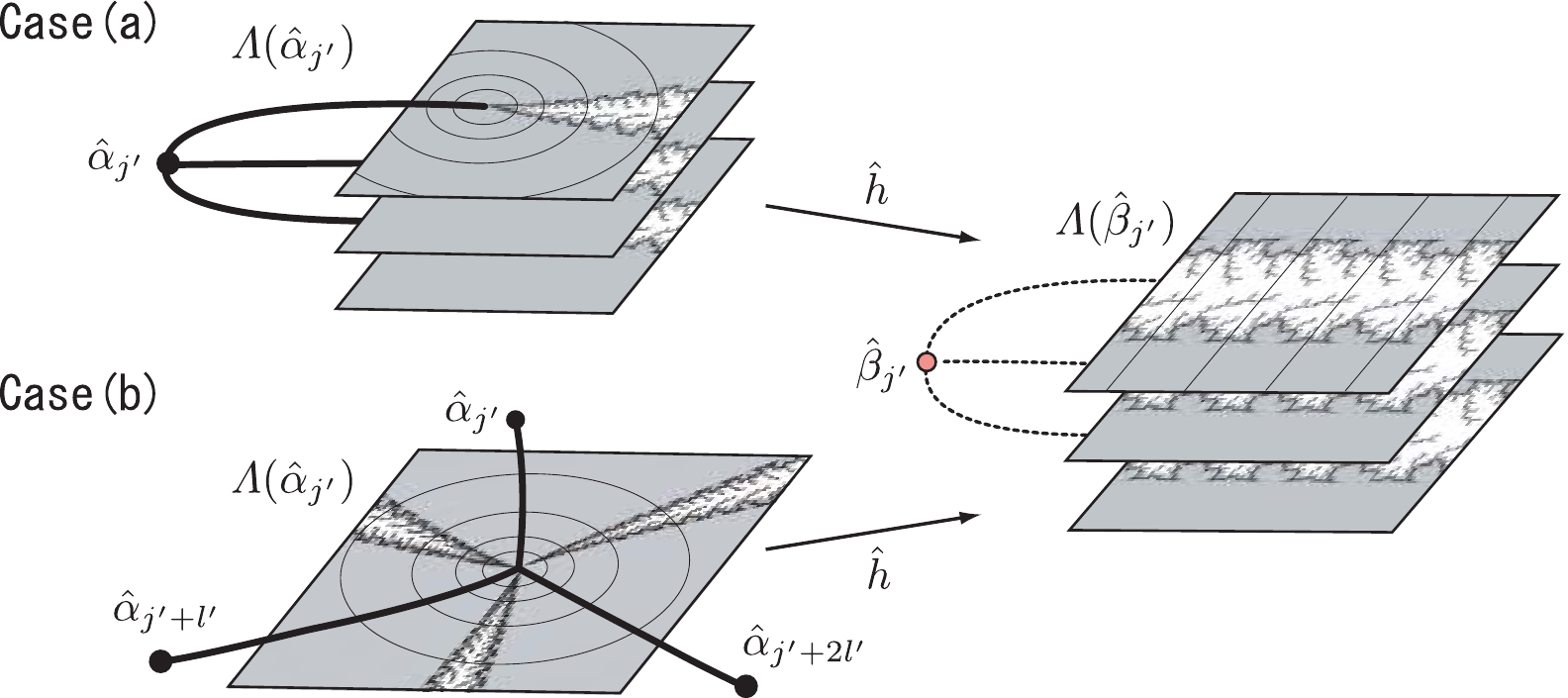}
\caption{A caricature of principal leaves $\Lam(\alhat_{j'})$ and $\Lam(\betahat_{j'})$ for $q'=3$. Here we take subscripts $j'$ modulo $l'$. Recall that $q=q'$ and $l=l'$ in Case (a) and $q=1<q'$ and $l=l'q'$ in Case (b). In Case (b), note that $\Lam(\alhat_{j'})=\Lam(\alhat_{j'+l'})=\Lam(\alhat_{j'+2l'})$. The heavy curves indicate $I(\alhat_{j'})$. The Rabbits are the case of $l'=1$.}\label{fig_principal_leaves}
}
\end{figure}

\figref{fig_principal_leaves} explains the situation of (2). 
As $I(\alhat_{j'})$ is pinched to $\betahat_{j'}$, the $q'$ components of $\Lam(\alhat_{j'})-I(\alhat_{j'})$ become the $q'$ leaves of $\Lam(\betahat_{j'})$. Then the hyperbolic action of $\fhat^{\lbar}$ on each path-connected component of $\Lam(\alhat_{j'})-I(\alhat_{j'})$ is sent to the parabolic action of $\ghat^{\lbar}$ on each leaf of $\Lam(\betahat_{j'})$. 

\parag{Example: The Cauliflowers.}
(See \figref{fig_cauli_leaves} again.) In this case the only cyclic $\Ihat_f$-component $I(\Ohat_f)$ is the invariant lift of the interval $[\al_1, \gam_1]$ by the projection $\pi_f:\Nf \to \Cbar$, where $\gam_1$ is the repelling fixed point. 
On the other hand, the corresponding set $\IO_f$ in the universal setting is the invariant lift of $(\al_1, \gam_1]$ by the projection $\pihat_f: \Af \to \Cbar$.

The set $\Lam_f-I(\Ohat_f)$ (resp. $\Lam_f^\a - \IO_f$) is isomorphic to a plane with a half-line removed.  
There are many other compactly contained $\Ihat_f$-components in $\Lam_f-I(\Ohat_f)$ (resp. $\II_f$-components in $\Lam_f^\a - \IO_f$) but the map $\hhat|~\Lam_f-I(\Ohat_f) \to \Lam_g$ (resp. $\hhat_\a|~\Lam_f^\a-\IO_f \to \Lam_g^\a$) just pinches them to points. Thus $\Lam_f-I(\Ohat_f)$ (resp. $\Lam^\a-\IO_f$) has the same topology as $\Lam_g$ (resp. $\Lam_g^\a$).

\parag{Example: The Rabbits.}
As in Figures \ref{fig_rabbit_leaves} and \ref{fig_principal_leaves}, Case (a) and Caes (b) exhibit different types of degenerations. 
What happens in the Case (a) Rabbits $(f_1 \to g)$ is similar to the Cauliflowers. 
The set $\Lam_{f_1}-I(\Ohat_{f_1})$ is topologically three planes with three half-lines removed, which are mapped onto three principal leaves of $\Lam_g$. 

However, in the Case (b) Rabbits $(f_2 \to g)$ the only cyclic $\Ihat_{f_2}$-component $\Ohat_{f_2}$ is a non-compact trivalent graph which divides the invariant leaf $\Lam_{f_2}$ into three simply connected components. These components are cyclic under $\fhat_2$ and mapped onto the three principal leaves of $\Lam_g$.

\parag{Proof of \thmref{thm_affine_part}.}
Let us start the proof of the theorem. We first show: 
\begin{lem}\label{lem_L'_and_L}
For any leaf $L'$ of $\Agn$, there exists a leaf $L$ of $\Afn$ such that $\hhat^{-1}(L') \subset L$. In particular, $\hhat^{-1}(L')$ is path-connected. 
\end{lem}
\begin{pf}(cf. \cite[Lemma 3.2]{LM}) Fix two distinct points $\zhat'=(z_0', z_{-1}', \ldots)$ and $\what'=(w_0',w_{-1}', \ldots)$ in $L'$. Let $\hat{\eta}'$ be a path in $L'$ joining $\zhat'$ and $\what'$. Then $\eta'_{-n}:=\pi_g \circ \ghat^{-n}(\hat{\eta}')$ is a path joining $z_{-n}'$ and $w_{-n}'$, and $\eta'_{-n}$ with $n \gg 0$ has a neighborhood $U_{-n}$ whose pull-back along $\zhat'$ and $\what'$ is eventually univalent. (That is, the path $\eta'_{-n}$ with $n \gg 0$ does not pass through $O_g$ and $\infty$.)

Choose any $\zhat=(z_0, z_{-1}, \ldots) \in \hhat^{-1}(\zhat')$ and $\what=(w_0, w_{-1}, \ldots) \in \hhat^{-1}(\what')$. 
For any $N \le 0$, even if $\eta'_{-N}$ passes through $I_g$, its preimage $h^{-1}(\eta_{-N}')$ is still a path-connected set by (1) of \thmref{thm_semiconj}. 
Since $z_{-N}$ and $w_{-N}$ are contained in $h^{-1}(\eta_{-N}')$, we can choose a path $\eta_{-N}$ joining $z_{-N}$ and $w_{-N}$. 
Let us fix $N \gg 0$ so that $\eta'_{-N}$ avoids both $O_g$ and $\infty$. 
Then $\eta_{-N}$ also avoids $I(O_f)=\bigsqcup_{j'} I(\al_{j'})$ and $\infty$. 
Now we can take a neighborhood of $\eta_{-N}$ whose pull-back along $\zhat$ and $\what$ is eventually univalent. Since we can lift paths $\skakko{\eta_{-N-n}}$ to a path in $\NN_f$ joining $\zhat$ and $\what$, both $\zhat$ and $\what$ are in the same leaf $L$ of $\Afn$. Thus $\hhat^{-1}(L')$ is a path-connected subset of $L$. 
\QED(\lemref{lem_L'_and_L})
\end{pf}

It is possible that $\hhat^{-1}(L_1') \cup \hhat^{-1}(L_2') \subset L$ for distinct leaves $L_1'$ and $L_2'$ of $\Agn$. In fact, it happens when $L_1'$ and $L_2'$ are principal leaves and $q=1<q'$ (Case (b)). 

As we have seen, the map $\hhat: \Afn - I(\Ohat_f) \to \Agn$ is a semiconjugacy, thus surjective and continuous. 
To show that this is a leafy and meek between two leaf spaces, it is enough to show (1) and (2). 

\parag{Proof of (2).}
We show (2) first. 
In both Case (a) and (b), the set $\Lam(\alhat_{j'})$ consists of $q$ principal leaves by \propref{prop_L}(2) and (3). 
Let $L$ be such a leaf in $\Lam(\alhat_{j'})$, and set $L^-:=L - I(\alhat_{j'})$ for simplicity. 
The degenerating arc $I(\al_{j'})=\pi_f(I(\alhat_{j'}))$ contains $q$ repelling periodic points. 
Let $\gam_L$ be the one of such repelling points so that its periodic lift $\gamhat_L$ is contained in $L$. 
Then there are $\qbar:=l/l'$ ($=1$ or $q'$ according to Case (a) or (b)) degenerating arcs in $I(\al_{j'})$ landing at $\gam_L$ downstairs, and such arcs are all invariant under $f^{\lbar}$. 
Since the quotient $L/\kk{\fhat^{\lbar}}$ is isomorphic to the quotient torus of the dynamics by $f^{\lbar}$ near $\gam_L=\pi_f(\gamhat_L)$, we conclude that the torus $L/\kk{\fhat^{\lbar}}$ is divided into $\qbar$ parts by $\qbar$ simple closed curves $(I(\alhat_{j'})-\skakko{\text{vertices}})/\kk{\fhat^{\lbar}}$. 
This implies that $L^-$ is the union of $\qbar$ path-connected components divided by $I(\alhat_{j'})$, say $S_1, \ldots, S_{\qbar}$. (In Case (a), we have $\qbar=1$ and $L^-$ is topologically a plane minus a half line.) 
In particular, each $S_i$ is simply connected.
Since $lq=l'q'$, we have $q'$ path-connected components in $\Lam(\alhat_{j'})-I(\alhat_{j'})$ in total.

Now take $L' \subset \Agn$ and $L \subset \Afn$ as in the lemma above. Suppose that $\hhat(L)$ contains a irregular point, that is, a point in $\Ohat_g \cup \{\hat{\infty}\}$. 
Since $\hhat^{-1}(\Ohat_g)=I(\Ohat_f)$ and $\hhat^{-1}(\{\hat{\infty}\})=\{\hat{\infty}\}$, it is equivalent to $L \cap I(\Ohat_f) \neq \emptyset$. 
This implies that $L$ is a principal leaf, i.e., $L \subset \Lam(\alhat_{j'})$ for some $j'$ modulo $l'$. 
Since $\Ohat_g$ and $L'$ are disjoint in $\Ng$, we have
$$
\hhat^{-1}(L') \subset L-\hhat^{-1}(\Ohat_g)=
L - I(\Ohat_g) = L - I(\alhat_{j'}) = L^-,
$$
where $L^-$ is the disjoint union of $S_1,\ldots, S_{\qbar}$ as above.
By the lemma above $\hhat^{-1}(L')$ is a path-connected subset of $L^-$. 
Since each $S_i$ are also path-connected, we can find some $i$ with $\hhat^{-1}(L') \subset S_i$. 
On the other hand, the set $\hhat(S_i)$ is also path-connected by continuity of $\hhat$ and thus contained in a leaf of $\Agn$, which must be $L'$. Thus we have $\hhat(S_i)=L'$ and it implies
$$
S_i ~\subset~ \hhat^{-1}(\hhat(S_i))~=~\hhat^{-1}(L') ~\subset~ S_i.
$$
Now let us check that $L'$ is one of $q'$ leaves of $\Lam(\betahat_{j'})$. 
To identify the leaf, we consider lifts of external rays. 
Let $\gamhat_L$ be the repelling periodic point in $L$ and set $\gam_L:=\pi_f(\gamhat_L)$. 
Now there are $\qbar$ external rays $R_f(\theta_1), \ldots, R_f(\theta_{\qbar})$ landing at the repelling periodic point $\gam_L$ downstairs, and they all have ray-period $\lbar$. 
By comparing the quotient tori upstairs and downstairs, we conclude that each of $S_i$ has an invariant lift of such an external ray, say $R_f(\thetahat_i)$.  
(Where $\thetahat_i$ is an cyclic angle in $\Thetahat \subset \TThat$.) 
Since $\hhat$ sends external rays to external rays without changing the angles (\thmref{thm_natural_extension}(2)), the image $R_g(\thetahat_i)=\hhat(R_f(\thetahat_i))$ is an invariant lift of a periodic external ray downstairs which must land on the parabolic point $\beta_{j'}$.  
Hence $L'$ must contain $R_g(\thetahat_i)$ and thus $L'$ is a principal leaf in $\Lam(\beta_{j'})$. 
By setting $L_i':=L'$, we have the desired correspondence between $\skakko{S_i}$ and $\skakko{L'_i}$ with $\hhat^{-1}(L_i')=S_i$. 

Since the same happens for each $L$ in $\Lam(\alhat_{j'})$, we have $\hhat^{-1}(\Lam(\betahat_{j'})) = \Lam(\alhat_{j'})-I(\alhat_{j'})$. By taking the union with respect to $j'$, we have $\hhat^{-1}(\Lam_g) = \Lam_f-I(\Ohat_f)$. 
Since any leaf of $\Afn$ and $\Agn$ is dense (\propref{prop_LM_4.5}(2)), the sets $\Lam_f-I(\Ohat_f)$ and $\Lam_g$ are connected, thus leaf spaces.
Now we conclude (2).

\parag{Proof of (1).}
Next we show (1): Suppose that $\hhat(L)$ avoids irregular points. Then $L$ is not a principal leaf. Since $\hhat(L)$ is a path-connected subset of $\Agn$, there is a leaf of $\Agn$ containing $\hhat(L)$, which must be $L'$. It follows that $\hhat(L)=L'$ and thus 
$$
L ~\subset ~\hhat^{-1}(\hhat(L))~=~\hhat^{-1}(L') ~\subset~ L.
$$
By property (2), the leaf $L'$ is not a principal leaf. 
Thus the preimage of any non-principal leaf is a non-principal leaf. 
Since every leaf in the affine parts are dense, the non-principal parts are also connected, thus leaf spaces. 
Since every leaf in the non-principal parts is isomorphic to $\C$, it follows that $\hhat$ restricted on the non-principal part is very meek. 
This implies property (1).

\parag{Affine laminations.}
Finally, the last part of the statement on the affine laminations is justified under the identification given by \propref{prop_Afn_Agn}, and by defining $\hhatA:=\iota_g \cc \hhat_\n \cc \wp$. \QED

\parag{Remark.}
The affine laminations $\Af$ and $\Ag$ in the universal setting have metrics induced from a ``universal metric" on $\UUhat^\a$. Such a metric is given in \cite[\S 3]{KL}. By using the fact that $\hhat_\n$ tends to identity in the natural extension, we can prove the Gromov-Hausdorff convergence of $\Af-\IO_f$ to $\Ag$ as $f \to g$ with respect to the universal metric. 

\section{Degeneration of the $\Hyp^3$-laminations}\label{sec_07}

In this section we consider degeneration of the $\Hyp^3$-laminations associated with $(f \to g)$.
The main tool is a three dimensional extension of the pinching semiconjugacy $\hhatA:\Af-\IO_f \to \Ag$ in the universal setting.
Unfortunately, we cannot extend $\hhatA$ to the whole vertical fibers of $\Af-\IO_f$ without losing the semiconjugating property. However, we will see that we can extend $\hhatA$ to a fairly large part of $\Hf$ (including the non-principal part of $\Hf$) still as a semiconjugacy onto $\Hg$, and it is enough to explain the degeneration of the dynamics. 
More precisely, we define an $\fhat$-invariant set $\VV_f$ (called \textit{valley}) in the principal leaves $\Lam_f^\h = (\Lam_f^\a)^\h$ and we will have a meek pinching semiconjugacy $\hhat_\h:\Hf -\VV_f \to \Hg$ (Theorems \ref{thm_Hyp^3} and \ref{thm_extension_principal_leaves}).

In particular, we will conclude that $\hhat_\h$ is very meek on the non-principal part, i.e., the geometry of the non-principal part of $\Hf$ is leafwise preserved as $f \to g$. 
In the principal part, we will see that the loxodromic (or hyperbolic) dynamics is degenerate to the parabolic dynamics as $f \to g$. 
However, the restriction of $\hhat_\h$ is still meek.

The semiconjugacy above will play a crucial role when we consider the quotient laminations in the next section.

\parag{Remarks on notation.} 
Here we will work in the universal setting. Let us recall some notation given in \secref{sec_02}. 
\begin{itemize}
\item To denote leaves of the affine {\it laminations}, we will use $L$ or $L'$ for simplicity. (We used them for leaves of the affine {\it parts}.) We also use $\phi:\C \to L$ or $\phi':\C \to L'$ for their uniformizations.
\item For $\phi:\C \to L$ above, the map $\pi_f \cc \wp \cc \phi:\C \to \Cbar$ is an element of $\Kf$. Recall that $\pihat_f:=\pi_f \cc \wp$ on $\Af$. 
Note that $\psi_{-N}:=\pihat_f \cc \fhat_\a^{-N} \cc \phi=\pi_f \cc \fhat_\n^{-N}  \cc \wp \cc \phi$ determines an element $\psihat=(\psi_{-N})_{N \ge 0}$ in $\Kfhat$. 

\item Recall that for $X \subset \Af$, we denote its fibers in $\Hf$ by $X^\h$. In particular, the uniformization of $L^\h$ associated with $\phi:\C \to L$ is denoted by $\phi^\h:\Hyp^3 \to L^\h$. 
\end{itemize}

\subsection{Pinching semiconjugacy on $\Hyp^3$-laminations}
Since the $\Hyp^3$-lamination $\Hf$ is an $\R^+$-bundle of the affine lamination $\Af$, one may imagine that the pinching semiconjugacy $\hhatA:\Af-\IO_f \to \Ag$ is naturally extended to a surjective map from $\Hf-\IO_f^\h \to \Hg$, 
and this would be enough to describe the degeneration of $\Hf$ to $\Hg$.
It is true, however, the map may not be a semiconjugacy; i.e., it may not preserve the dynamics on $\Hf$ and $\Hg$. 
Such a map is useless when we consider their quotient 3-laminations. 

Here we try to have the best possible 3D-extension of $\hhatA$ as a semiconjugacy. 
The main difficulty comes from the universal setting: We need to manage plenty of complicated notation.

\parag{Pinching semiconjugacy between the non-principal parts.}
We naturally have an extension except over the principal leaves. We first establish:
\begin{thm}[Non-principal part]\label{thm_Hyp^3}
The very meek pinching semiconjugacy $\hhatA:\Af-\Lam_f^\a \to \Ag-\Lam_g^\a$ extends to a very meek semiconjugacy $\hhatH:\Hf - \Lam_f^\h \to \Hg - \Lam_g^\h$. In particular, for $\xhat \in \Af-\Lam_f^\a$, the extended semiconjugacy $\hhatH$ sends the fiber of $\xhat$ onto the fiber of $\hhatA(\xhat)$.
\end{thm}

\begin{figure}[htbp]
\begin{center}
\includegraphics[width=.55\textwidth]{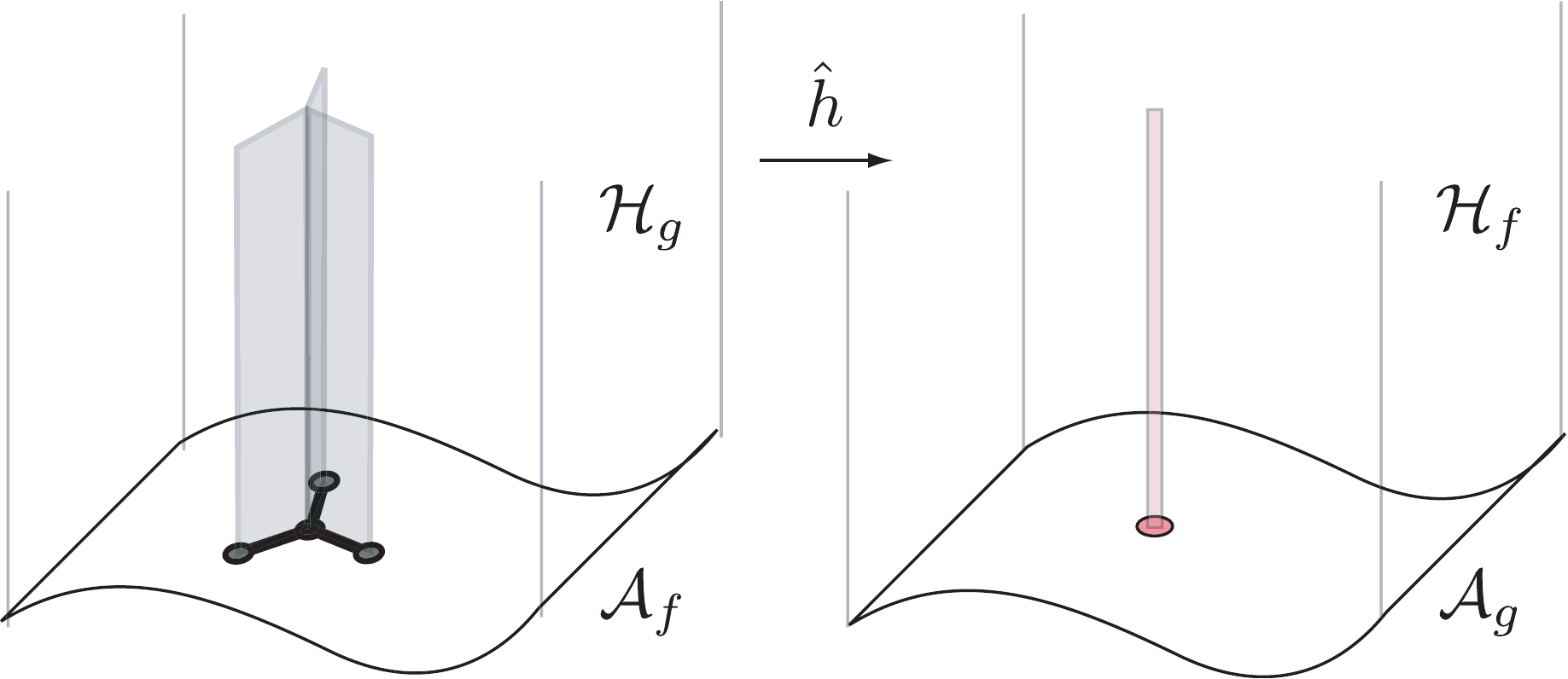}
\end{center}
\caption{Degeneration of a non-periodic $\II_f^\h$-component}
\label{fig_deg_Hf}
\end{figure}

What happens in the non-principal part is very simple: 
Each $\II_f^\h$-component is pinched to be an $\II_g^\h$-component, which is a vertical fiber of an $\II_g$-component (\figref{fig_deg_Hf}).

Before the proof of this theorem, let us start with a lemma which slightly generalize \cite[Lemma 9.2]{LM}:

\begin{lem}[Affinely natural extension]\label{lem_extension}
Let $H:\C \to \C$ be a continuous and surjective map such that for each $w \in \C$ the function 
$$
v_w(t) \ee \max_{|W|=t} |H(w+W)-H(w)|
$$ 
is a non-decreasing function from $\R_{\ge 0}$ onto itself. Then the map $e(H):\Hyp^3\cup\C \to \Hyp^3\cup\C$ defined by 
$$
e(H)(w,t) \dee \kakko{H(w), v_w(t)} 
$$
is a continuous and surjective map which restrict $H$ on $\C$. Moreover, the following is satisfied:
\begin{enumerate}
\item The extension is affinely natural: If $\delta_1$ and $\delta_2$ are complex affine maps of $\C$ then $e(\delta_1)$ and $e(\delta_2)$ are the unique possible similarities of $\Hyp^3$, and
$$
e(\delta_1 \cc H \cc \delta_2) \ee e(\delta_1)\cc e(H)\cc e(\delta_2).
$$
\item $e(H)$ depends continuously on $H$, in the compact-open topology on maps of $\C$ and $\Hyp^3$.
\item The vertical line over $z$ is mapped onto to that of $H(z)$ (but not one-to-one in general).
\end{enumerate}
\end{lem}
The proof is routine and left to the reader.

\parag{Proof of \thmref{thm_Hyp^3}.}
For a leaf $L$ of $\Af -\Lam_f^\a$, the set $L':=\hhatA(L)$ is a leaf of $\Ag-\Lam_g^\a$ (\thmref{thm_affine_part}(1)). Take uniformizations $\phi:\C \to L$ and $\phi':\C \to L'$. Then the map $H:=(\phi')^{-1} \cc \hhatA \cc \phi$ is a surjective and continuous map satisfying the condition of the lemma above. (This part is also routine and left to the reader.) Let us define $\hhatH: L^\h \to L'^\h$ by $\hhatH:=\phi'^\h \cc e(H) \cc (\phi^\h)^{-1}$, where $\phi^\h$ and $\phi'^\h$ are the uniformizations of $\Hyp^3$-leaves $L^\h$ and $L'^\h$ naturally induced by $\phi$ and $\phi'$, respectively. Then one can easily check that this definition does not depend on the choice of $\phi$ and $\phi'$ by affine naturality of the extension. Moreover, the map $\hhatH$ sends a fiber onto a fiber by property 3 of the lemma above. 

Now $\hhatH: \Hf - \Lam_f^\h \to \Hg-\Lam_g^\h$ is defined just leafwise, however, it is not hard to check the continuity of this map by applying property 2 of \lemref{lem_extension} to each local product box of the laminations. (This part is more justified in the proof of \thmref{thm_extension_principal_leaves}.) Property 1 of \lemref{lem_extension} insures that $\hhatH$ semiconjugates $\fhat$ to $\ghat$ in our domain. 

Finally note that $\Hf-\Lam_f^\h$ and $\Hf-\Lam_f^\h$ are leaf spaces by minimality (\cite[\S 7.7]{LM}).
\QED 


\paragraph{}
Next we continuously extend $\hhatH$ above to some part of the principal leaves $\Lam_f^\h$. We start with carving.

\begin{figure}[htbp]
\centering{
\vspace{0cm}\hspace{0cm}
\includegraphics[width=0.6\textwidth]{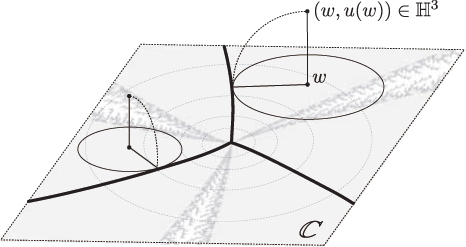}
}
\caption{Carving a valley in the case of $q=1<q'=3$, Case (b).}\label{fig_valley}
\end{figure}

\begin{figure}[htbp]
\begin{center}
\includegraphics[width=0.80\textwidth]{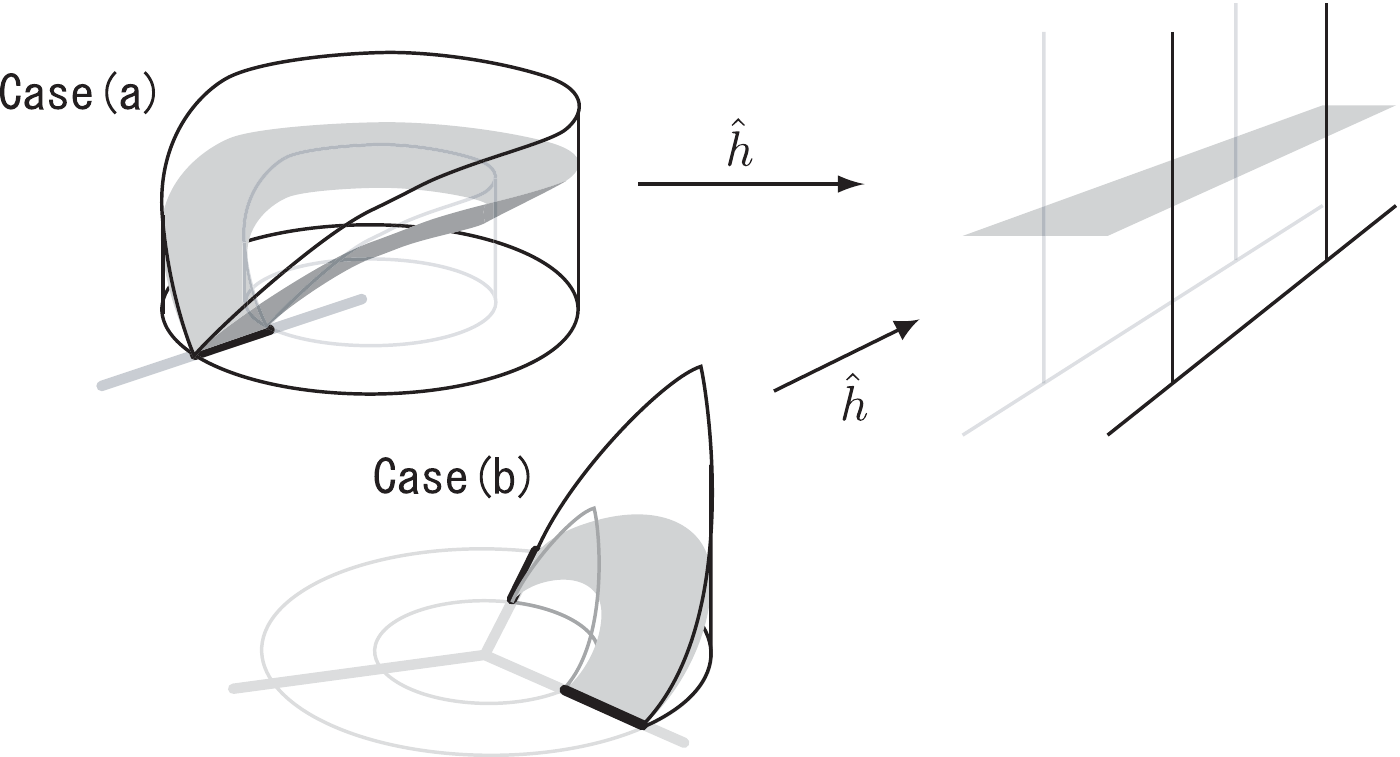}
\end{center}
\caption{A fundamental domain of the action $\fhat^{\lbar}:\Lam_f^\v \to \Lam_f^\v$ is mapped onto that of $\ghat^{\lbar}:\Lam_g^\h \to \Lam_g^\h$}\label{fig_cusp}
\end{figure}

\parag{Valleys carved by degenerating arcs.}
For an $\alhat_{j'} \in \Ohat_f$ with $j'$ modulo $l'$, take a leaf $L$ of $\Af$ contained in $\Lam^\a(\alhat_{j'})$. Since any leaf is isomorphic to $\C$, choose a uniformization $\phi: \C \to L$, and consider a continuous function $u=u_\phi:\C \to [0, \infty)$ defined by $u(w):=\mathrm{dist}_\C(w, \phi^{-1}( L \cap \IO_f   ))$, where the distance is measured in Euclidean metric on $\C$. Set 
$$
L^\v~:=~\phi^\h \kakko{ \skakko{ (w,t) \in \Hyp^3: t < u(w) } } 
~\subset~ ( L  - \IO_f )^\h,
$$
where $\phi^\h: \Hyp^3 \to L^\h$ is the natural uniformization induced by $\phi$. (See \figref{fig_valley}.) 
It is not difficult to show that the definition of $L^\v$ does not depend on the choice of $\phi$. In fact, for another uniformization $\phi_1:\C \to L$, there exists a unique affine map $\delta_1(w)=a + bw$ such that $\phi=\phi_1 \cc \delta_1$. Thus $u_{\phi}(w)=u_{\phi_1}(w)/|b|$ and the difference is canceled when we consider $\phi^\h=\phi_1^\h \cc e(\delta_1)$.  
We call this well-defined set $L^\v$ a \textit{carved principal leaf}, and the complement $L^\h-L^\v$ the \textit{valley} of $L^\h$.  
(This is just like a valley carved by a degenerating ``stream" $\IO_f$ with a source which is the lift of a repelling point. See \figref{fig_cusp}.)
Note that $L^\v$ has $\qbar=l/l'$ path-connected components. 
Take a union of $L^\v$ over all leaves in $\Lam(\alhat_{j'})$, and denote it by $\Lam(\alhat_{j'})^\v$. 
Then the carved principal leaves $\Lam_f^\v:=\bigsqcup \Lam(\alhat_{j'})^\v$ and the removed part $\VV_f:=\bigsqcup (\Lam(\alhat_{j'})^\h-\Lam(\alhat_{j'})^\v)$ are completely invariant under the action of $\fhat:\Hf \to \Hf$. 
Note that both $\Lam_f^\v$ and $\VV_f$ have $\lbar(=lq=l'q')$ path-connected components. 
We also say $\VV_f$ is the \textit{valley} of $\Hf$.
Again by minimality of leaves, $\Hf-\VV_f$ and $\Lam_f^\v$ are connected, thus leaf spaces. 

\parag{Examples: The Cauliflowers.}
For the Cauliflowers $(f \to g)$, the ``stream" $\IO_f$ is the half-line on the principal leaf $\Lam_f^\a$. 
The carved principal leaf $\Lam_f^\v=\Lam_f^\h-\VV_f$ has a fundamental domain as shown in \figref{fig_cusp}, Case (a).

\parag{Examples: The Rabbits.}
For the Case (a) Rabbits $(f_1 \to g)$ the carved principal leaf $\Lam_{f_1}^\v$ is leafwise similar to the Cauliflower case. 
In particular, it has three simply-connected components cyclic under $\fhat_{1\h}$.
For the Case (b) Rabbits $(f_2 \to g)$, the ``stream" $\IO_{f_2}$ has three directions but the carved principal leaf $\Lam_{f_2}^\v$ also has three simply connected components cyclic under $\fhat_{2\h}$. See \figref{fig_cusp}, Case (b).

\parag{}
Now we are ready to state the theorem:
\begin{thm}[Extention to the principal part]\label{thm_extension_principal_leaves}
There exists a continuous extension of $\hhatH: \Hf-\Lam_f^\h \to \Hg-\Lam_g^\h$ above to a meek pinching semiconjugacy $\hhatH :\Hf-\VV_f \to \Hg$ with the following properties:
\begin{enumerate}[\rm (1)]
\item For each $j'$ modulo $l'$, there are $q'$ path-connected components $S_1^\v, \ldots, S_{q'}^\v$ in $\Lam(\alhat_{j'})^\v$ such that $S_{i}^\v=\hhat_\h^{-1}(L_{i}'^\h)$ for the $q'$ principal leaves $L_1'^\h, \ldots, L_{q'}'^\h$ in $\Lam(\betahat_{j'})^\h$. In particular, $\hhat_\h^{-1}(\Lam_g^\h)=\Lam_f^\v=\Lam_f^\h-\VV_f$.
\item For $\xhat \in \Lam_f^\a$, the semiconjugacy $\hhatH$ sends $\skakko{\xhat}^\h \cap \Lam_f^\v$ onto the fiber of $\hhatA(\xhat)$.
\end{enumerate}
\end{thm}
Thus we can understand the dynamics of $\ghat:\Hg \to \Hg$ by pinching $\fhat:\Hf-\VV_f \to \Hf$. The proof breaks into two main steps: First we define a semiconjugacy $\hhat_\h:\Lam_f^\v \to \Lam_g^\h$. Then we show that it is a continuous extension of the pinching map $\hhat_\h$ in \thmref{thm_Hyp^3}. To show this continuity, we have to go back to the universal setting.

\parag{Proof of \thmref{thm_extension_principal_leaves}: Construction.}
We first describe $L^\v$ of $L \subset \Lam(\alhat_{j'})$ in detail. Let $S_1, \ldots, S_{\qbar}$ be the connected components of $L-\IO_f$. Each of them are sent onto leaves in $\Lam_g^\a$, say $L'_1, \ldots, L'_{\qbar}$ (\thmref{thm_affine_part}(2)). Let $\phi$ and $\phi'_i$ be uniformizations of $L$ and $L'_i$ for $i=1, \ldots, \qbar$. We denote $\phi^{-1}(S_i)$ by $\mathbb{S}_i \subset \C$. Then we can define a surjective and continuous map $H_i:\mathbb{S}_i \to \C$ by $H_i:= (\phi_i')^{-1} \cc \hhatA \cc \phi$. For each $w \in \mathbb{S}_i$, set $u(w)=\mathrm{dist}_{\C} (w, \partial \mathbb{S}_i)$ and we define $v_w: [0, u(w)) \to [0, \infty) $ by
$$
v_w(t) \ee \max_{|W|=t} |H_i(w+W)-H_i(w)|,
$$
which is a non-decreasing function from $[0, u(w))$ onto $[0, \infty)$. Set
\begin{align*}
\mathbb{S}_i^\v~&:=~\skakko{(w,t) \in \Hyp^3: w \in \mathbb{S}_i,~t < u(w)}
~~\text{and}~~ \\
S_i^\v~&:=~\phi^\h(\mathbb{S}_i^\v).
\end{align*}
Note that $L^\v$ is the disjoint union of $S_i^\v$ for $i=1, \ldots, \qbar$. Now it is not difficult to check that 
there exists a continuous extension of $H_i:\mathbb{S}_i \to \C$ to $e(H_i):\mathbb{S}_i^\v \to \Hyp^3$. Moreover, this extension is affinely natural: That is, property 1 of \lemref{lem_extension} holds over $\delta_2^{-1}(\mathbb{S}_i)$.
We define $\hhatH:S_i^\v \to L_i'^\h$ by $\hhatH:=\phi_i'^\h \cc e(H_i) \cc (\phi^\h)^{-1}$. By affine naturality, one can check that this map does not depend on the choice of the uniformizations $\phi$ and $\phi_i'$; and that $\hhatH$ semiconjugates the actions of $\fhat^{l'q'}: S_i^\v \to S_i^\v$ onto $\ghat^{l'q'}: L_i'^\h \to L_i'^\h$. Now we have $\hhatH:\Lam_f^\v \to \Lam_g^\h$ that is defined leafwise. 

\parag{Continuity.} 
Next we prove that the map $\hhatH: \Hf -\VV_f \to \Hg$ is continuous; that is, \textit{for any sequence $\xtil_n \to \xtil$ within $\Hf-\VV_f$, the sequence $\ytil_n:=\hhatH(\xtil_n)$ converges to $\ytil:=\hhatH(\xtil)$.}

\parag{Convergence in $\Hf$.}
To characterize the convergence in $\Hf-\VV_f$, we need to start with a lot of notation. For each $\xtil_n \in \Hf-\VV_f$, take a representative $(\psi_{n,-N})_{N \ge 0}$ in $\Kfhat$. Set $\xhat_n:=\pr(\xtil_n) \in \Af$ and $\zhat_n=(z_{n,-N})_{N \ge 0}:=\wp(\xhat_n) \in \Afn$. Let $L_n$ denote the leaf $\iota_f(L(\zhat_n))$ in $\Af$ containing $\xhat_n$. Then there exists a uniformization $\phi_n:\C \to L_n$ such that $\psi_{n,-N}=\pihat_f \cc \fhat_\a^{-N} \cc \phi_n$. For $\xhat_n \in L_n$, let $S_n$ be the connected component of $L_n-\IO_f$ containing $\xhat_n$. It may be a sector when $\xtil_n \in \Lam_f^\v$, or whole a leaf $L_n$ when $\xtil_n \notin \Lam_f^\v$. Set $\mathbb{S}_n:=\phi_n^{-1}(S_n)$. We also define the function $u_{\phi_n}: \C \to [0, +\infty]$ by $u_{\phi_n}(w):=\dist_\C(w, \partial \mathbb{S}_n)$ if $\xtil_n \in \Lam_f^\v$, and by $u_{\phi_n}(w) \equiv +\infty$ otherwise.

For $\xtil$, we also define the following sets just by removing $n$'s from the definitions above: $(\psi_{-N})_{N \ge 0}$ in $\Kfhat$, $\xhat$ in $\Af$, a leaf $L$, uniformization $\phi:\C \to L$, $S \subset L$ containing $\xhat$, $\mathbb{S}=\phi^{-1}(S)$, and the function $u_{\phi}:\C \to [0, +\infty]$.

Since $\xtil_n \to \xtil$, we may assume that $\psi_{n,-N} \to \psi_{-N}$ uniformly on compact sets for any fixed $N$. Note that for the uniformizations $\phi_n^\h:\Hyp^3 \to L_n^\h$ and $\phi^\h:\Hyp^3 \to L^\h$, both $\xtil_n$ and $\xtil$ have coordinate $(0,1)$ with respect to $\psi_n$ and $\psi$: i.e., $(\phi_n^\h)^{-1}(\xtil_n)=(\phi^\h)^{-1}(\xtil)=(0,1)$. Moreover, both $u_{\phi_n}(0)$ and $u_{\phi}(0)$ are greater than $1$ since $\xtil_n, \xtil \in \Hf-\VV_f$. 

Here is an important lemma for the continuity:
\begin{lem}\label{lem_lower_semiconti}
Set $C_n:=u_{\phi_n}(0)$ and $C:=u_{\phi}(0)$. Then  
$$
\liminf_{n \to \infty} C_n ~\ge~ C.
$$
\end{lem}

\paragraph{Proof of \lemref{lem_lower_semiconti}.}
First we consider the case of $\xtil \in \Lam_f^\v$. Since $C_n=+\infty$ when $\xtil_n \notin \Lam_f^\v$, it is enough to show the case of $\xtil_n \in \Lam_f^\v$ for all $n \gg 0$. 

Take a disk $\tilde{D}:=\D_{\tilde{C}}$ with $\e:=C/\tilde{C} \ll 1$. If $N \gg 0$, the map $\psi_{-N}:\tilde{D} \to \Cbar$ is univalent. 
Moreover, the sequence $\psi_{n,-N}|_{\tilde{D}}$ converges uniformly to  $\psi_{-N}|_{\tilde{D}}$ and thus is univalent for all $n \gg 0$. 
Set $B:=\psi_{-N}(\D_{2 C}) \subset \C$. 
(Since $f$ is a polynomial, infinity is the exceptional value for all elements of $\Kf$.) 
Then $B \cap I(O_f) \neq \emptyset$ since $\D_{2 C} \cap \partial \mathbb{S}_n \neq \emptyset$. 
In particular, there is a point in $B \cap I(O_f)$ that attains $\dist_\C(z_{-N}, B \cap I(O_f))$.
Since $z_{n,-N}=\psi_{n,-N}(0)$ tends to $z_{-N}=\psi_{-N}(0) \in B$, we have 
$$
\dist_\C(z_{n,-N}, B \cap I(O_f)) ~\to~ \dist_\C(z_{-N}, B \cap I(O_f)).
$$
By the Koebe distortion theorem, 
\begin{align*}
C_n &\ee |(\psi_{n,-N})'(0)|^{-1} \, \dist_\C(z_{n,-N}, B \cap I(O_f))
\,(1+O(\e))\\
C &\ee |(\psi_{-N})'(0)|^{-1} \, \dist_\C(z_{-N}, B \cap I(O_f))
\,(1+O(\e)).
\end{align*}
This implies that $C_n$ is arbitrarily close to $C$ as $n \to \infty$ if we take $\tilde{C},~N \gg 0$ and $\e \ll 1$ as above. Hence we have $\lim C_n = C$. 

Next we consider the case of $\xtil \notin \Lam_f^\v$, that is, $C=+\infty$. Suppose that $\liminf C_n=M< +\infty$. Passing through a subsequence, we may assume that $C_n \to M$ and $\xtil_n \in \Lam_f^\v$ for all $n \gg 0$. By a similar argument as above, we have 
\begin{align*}
\dist_\C(z_{-N}, \psi_{-N}(\D_{2C}) \cap I(O_f)) &~\sim~ \dist_\C(z_{n,-N}, B \cap I(O_f)) \\
  &~\sim~ M|(\psi_{n,-N})'(0)| \\
  &~\sim~ M|(\psi_{-N})'(0)|
\end{align*}
for all $N,n \gg 0$. Since $(\psi_{-N-m})'(0)=(\psi_{-N})'(0)/(f^{m})'(z_{-N-m})$ for all $m>0$, we have 
$$
\dist(z_{-N-m}, \psi_{-N-m}(\D_{2C}) \cap I(O_f)) ~\asymp~
|(\psi_{-N-m})'(0)| ~\to~ 0~~(m \to \infty)
$$
by hyperbolicity of $f$. This implies $\zhat=(z_{-N})_{N \ge 0}$ accumulates on $I(O_f)$ but it is possible only if $\xhat \in \Lam_f^\a$. It contradicts $\xtil \notin \Lam_f^\v$.\QED{\small(\lemref{lem_lower_semiconti})}

\parag{Convergence in $\Hg$.}
For $\ytil_n$ and $\ytil$ in $\Hg$, we also define the following sets as we did for $\xtil$: $(\psi_{n,-N}')_{N \ge 0}$ and $(\psi_{-N}')_{N \ge 0}$ in $\Kghat$, $\yhat_n$ and $\yhat$ in $\Ag$, leaves $L_n'$ and $L'$ of $\Ag$ containing $\yhat_n$ and $\yhat$, and their uniformizations $\phi_n':\C \to L_n'$ and $\phi':\C \to L'$. 
One can easily check the following properties:
\begin{enumerate}[(1)]
\item $\yhat_n=\hhatA(\xhat_n)$ and $\yhat=\hhatA(\xhat)$. Thus $\yhat_n \to \yhat$ in $\Ag$ by continuity of $\hhatA$.
\item $L_n'=\hhatA(S_n)$ and $L'=\hhatA(S)$. 
\item $\psi_{n,-N}'=\pihat_g \cc \ghat_\a^{-N} \cc \phi_n'$ and $\psi_{-N}'=\pihat_g \cc \ghat_\a^{-N} \cc \phi'$.
\end{enumerate}
Moreover, property (1) insures that we can choose the representatives satisfying $\psi_{n,-N}' \to \psi_{-N}'$ for any fixed $N$ as $n \to \infty$. Suppose that $\ytil_n$ and $\ytil$ have coordinates $(0, \rho_n)$ and $(0, \rho)$ with respect to $\psi'_n$ and $\psi'$ respectively. Now it is enough to show that $\rho_n \to \rho$ for the continuity of $\hhatH$.

Let $H_n:\mathbb{S}_n \to \C$ and $H:\mathbb{S} \to \C$ be surjective maps defined by $H_n:=(\phi_n')^{-1} \cc \hhatA \cc \phi_n$ and $H:=(\phi')^{-1} \cc \hhatA \cc \phi$. Now we have the following diagram:
$$
\begin{CD}
 ~\mathbb{S}_n~ @>{\phi_n} >> ~S_n~ @>{\fhat_\a^{-N}} >> ~\fhat_\a^{-N}(S_n)~ @>{\pihat_f} >> ~\Cbar~ \\
 @V{H_n}VV       @VV{\hhatA}V   @VV{\hhatA}V   @VV{h}V   \\
 ~\C~          @>{\phi'_n} >> ~L_n'~ @>{\ghat_\a^{-N}} >>~\ghat_\a^{-N}(L_n')~ @>{\pihat_g} >> ~\Cbar. 
\end{CD}
$$
By property (3) and this diagram, we have 
$$
\psi_{n,-N}'\cc H_n=h \cc \psi_{n,-N}~~\text{and}~~
\psi_{-N}'\cc H=h \cc \psi_{-N}.
$$ 
Take an $\e \ll 1$ and an $N \gg 0$ such that $1+\e < C=u_\phi(0)$ and $\psi_{-N}:\D_{1+\e} \to \C$ is univalent. By \lemref{lem_lower_semiconti}, the disk $\D_{1+\e}$ is contained in $\mathbb{S}_n$ for all $n \gg 0$. Moreover, the sequence $\psi_{n,-N}$ converges  to $\psi_{-N}$ uniformly on $\D_{1+\e}$ and they are all univalent on it. By the dynamics of $f$ and $g$ related by $h$, this implies $\psi_{n,-N}'$ and $\psi_{-N}'$ are also univalent on $H_n(\D_{1+\e})$ and $H(\D_{1+\e})$ respectively. Thus we have 
$$
H_n=(\psi_{n,-N}')^{-1} \cc h \cc \psi_{n,-N}~~\longrightarrow~~
H=(\psi_{-N}')^{-1} \cc h \cc \psi_{-N}
$$ 
uniformly on $\D_{1+\e}$. By definition of $\hhatH$, we have
$$
\rho_n=\max_{|W|=1 \atop w=0}|H_n(w+W)-H_n(w)|~~\text{and}~~
\rho=\max_{|W|=1 \atop w=0}|H(w+W)-H(w)|,
$$
and thus $\rho_n \to \rho$. This completes the proof of \thmref{thm_extension_principal_leaves}. \QED

\parag{Remark.}
Again $\Hf$ and $\Hg$ have metrics induced from a universal metric on $\UUhat^\h$. The Gromov-Hausdorff convergence of $\Hf -\VV_f$ to $\Hg$ can be proved by using a similar idea as the proof above and the fact that $h$ tends to identity.

\section{Degeneration of the quotient 3-laminations}\label{sec_08}
From this section we will mainly deal with the quotient laminations and their ends. (Recall the setting of \secref{sec_04}.)

For our degeneration pair $(f \to g)$, the meek pinching semiconjugacy $\hhatH: \Hf-\VV_f \to \Hg$ in the previous section provides us with a natural meek pinching map $\hhatM: \Mf-\vv_f \to \Mg$, where $\vv_f=\VV_f/\fhat$ is the quotient valley (\thmref{thm_quotient}). 
This gives a reasonable picture of the degeneration of the quotient 3-laminations. 
We conclude that the geometric difference of $\Mf$ and $\Mg$ appears only in the principal parts. 
Indeed, the unique principal leaf of $\Mf$ is a solid torus, but that of $\Mg$ is a rank one cusp. 
They are both homeomorphic to $\D \times \T$, thus we will see that the leafwise topology of $\Mf$ is preserved as $f \to g$. 

In addition, we check that the map $\hhatM: \Mf-\vv_f \to \Mg$ extends to their conformal boundaries. 

\subsection{Pinching map on the quotient 3-laminations}
For our degeneration pair $(f \to g)$, we consider the quotient 3-laminations $\Mf=\Hf/\fhat$ and $\Mg=\Hg/\ghat$, and give a description of $\Mg$ by translating the topological/geometric properties of $\Mf$.

For any $\tilde{x} \in \Hf$ and any integer $m \in \Z$, the quotient by $\fhat$ identifies $\fhat^m(\tilde{x})$ and $\tilde{x}$. We denote this identified class $x=[\xtil]$. We will generally use $x$ or $y$ for points of $\Mf$, and $\ell$ for a leaf of $\Mf$.

\parag{Solid tori and rank one cusps.}
Let us briefly recall two kinds of hyperbolic 3-manifolds with $\pi_1$ isomorphic to $\Z$. 
We regard the hyperbolic 3-space $\Hyp^3$ as the set of all $(z, t) \in \C \times \R^+$ with hyperbolic metric $ds^2=(|dz|^2+dt^2)/t^2$.
A \textit{solid torus} (or \textit{handle body}) is a hyperbolic 3-manifold isomorphic to $\Hyp^3/\kk{(z,t) \mapsto (\lam z, |\lam| t)}$, where $|\lam|>1$. 
The \textit{core curve} of a solid torus is the only closed geodesic corresponding to the $t$-axis.  
A \textit{rank one cusp} is a hyperbolic 3-manifold isomorphic to $\Hyp^3/\kk{(z,t) \mapsto (z+1, t)}$. 
A rank one cusp has non-trivial closed curves with arbitrarily short length. 
Note that both solid tori and rank one cusps have the same topology homeomorphic to $\D \times \T$. 

\parag{Principal leaves and valley.}
The quotient laminations can be divided into the principal parts and the non-principal parts. Let us characterize the leaves in these parts:
\begin{prop}[Quotient leaves]\label{prop_qpl}
For any degeneration pair $(f \to g)$, the leaves in the quotient laminations are classified as follows:
\begin{itemize}
\item {\bf Non-principal part.} For any non-principal leaf $L^\h$ of $\Hf$ (thus in $\Hf-\Lam_f^\h$), its quotient leaf is isomorphic to either a solid torus or $\Hyp^3$ according to $L^\h$ is periodic or non-periodic. The same is true for the non-principal leaves of $\Hg$. 

\item {\bf Principal part.} The quotients $\ell_f := \Lam_f^\h/\fhat_\h$ and $\ell_g := \Lam_g^\h/\ghat_\h$ are a solid torus and a rank one cusp respectively. 
\end{itemize}
\end{prop}

Note that even if $f$ (or $g$) has two or more principal leaves, the quotient is only one leaf. 

\begin{pf}
The statement about the non-principal parts is straightforward since periodic points in the non-principal parts are the lifts of repelling periodic points downstairs.

Any principal leaf $L^\h$ in $\Lam_f^\h$ is invariant under the iteration of $\fhat^\lbar=\fhat_\h^\lbar$ where $\lbar=lq=l'q'$. Since $\fhat^\lbar:L^\h \to L^\h$ is a hyperbolic action (\propref{prop_L}(5)) and every point in $\Lam_f^\h$ cyclically lands on $L^\h$, the quotient $\ell_f$ of $\Lam_f^\h$ is isomorphic to a solid torus with leafwise ideal boundary isomorphic to the torus $L/\fhat_\a^\lbar$. (Note that the \textit{conformal boundary} of $\ell_f$ is different. See the remark below.) 

Similarly, any leaf $L'^\h$ in $\Lam_g^\h$ has period $\lbar$ too and the action $\ghat^\lbar:L'^\h \to L'^\h$ is conjugate to a translation in $\Hyp^3$. Thus the quotient $\ell_g$ of $\Lam_g^\h$ is isomorphic to a rank one cusp whose \textit{leafwise} ideal boundary is isomorphic to $L'/\ghat_\a^\lbar$, or the cylinder $\C/\Z$.
\QED
\end{pf}

We also say $\ell_f$ and $\ell_g$ the \textit{principal leaves} of $\Mf$ and $\Mg$ respectively. 
Since any component of $\VV_f \subset \Lam_f^\h$ is invariant under $\fhat^\lbar$, its quotient $\vv_f:=\VV_f/\fhat$ is well-defined. We say $\vv_f$ is the \textit{valley} of $\ell_f$.

\begin{figure}[htbp]
\begin{center}
\includegraphics[width=.75\textwidth]{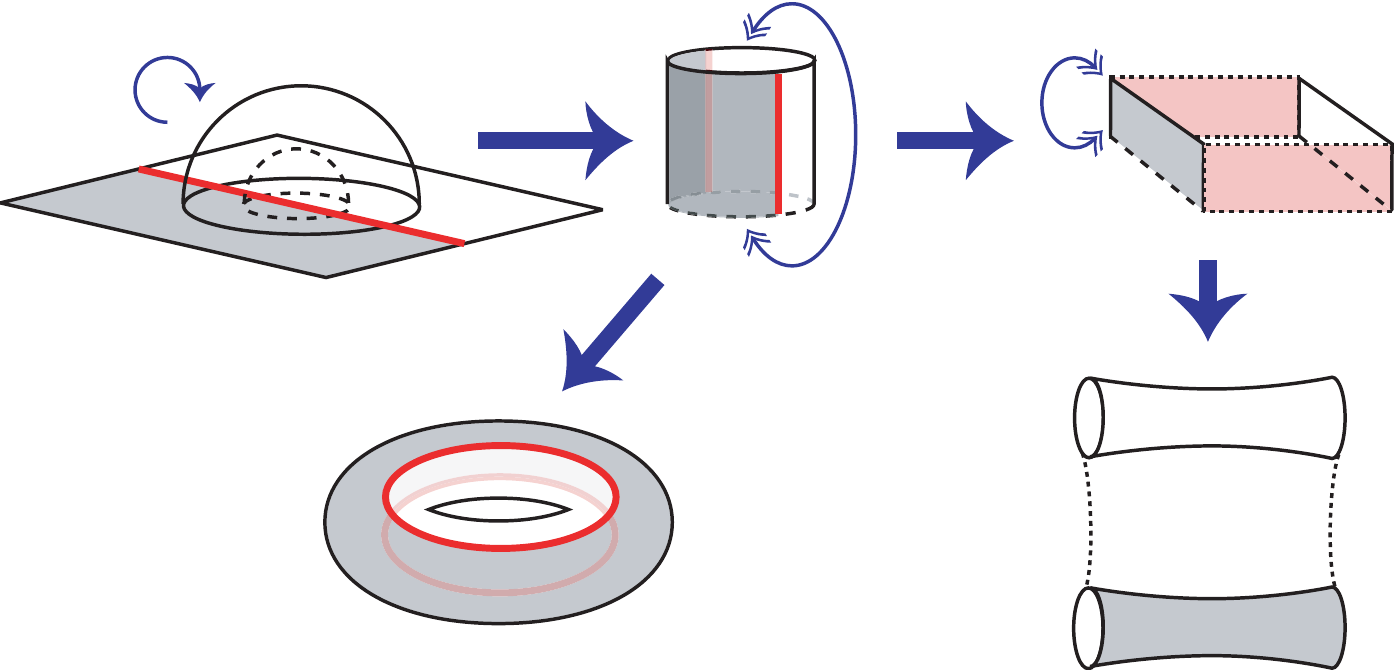}
\end{center}
\caption{The leafwise ideal boundary of $\ell_0$ (at the lower left) and its conformal boundary in $\Bar{\MM_{f_0}}$ (at the lower right). The difference comes from how we represent a solid torus: As $\D \times \T$, or as (an annulus)$\times (0,1)$.
The red curves indicate the lifted Julia set.}
\label{fig_boundaries}
\end{figure}

\parag{Remark on the boundaries of leaves.}
For any leaf $\ell$ in $\Mf$, its leafwise ideal boundary and global conformal boundary are different. 
For example, consider the repelling fixed point $z=1$ of $f_0(z)=z^2$ and the invariant leaf $L_0$ of $\AAA_{f_0}$ corresponding to the backward orbits converging to $z=1$. 
Then the action $\fhat_0: L_0 \to L_0$ is conformally conjugate to $w \mapsto 2w$. 
By taking the quotient of its 3D extension $\fhat_{0}=\fhat_{0\h}: L_0^\h \to L_0^\h$, we have a leaf $\ell_0$ of $\MM_{f_0}$ isomorphic to a solid torus whose ideal boundary is conformally isomorphic to the torus $\Cstar/\kk{w \mapsto 2w}$. 
But this ideal boundary contains points in the Julia set $\JJ_{f_0}$. 
The global conformal boundary $\partial \ell_0$ which is actually visible from $\MM_{f_0}$ is two annuli isomorphic to $\Hyp/\kk{w \mapsto 2w}$ and $\Hyp^-/\kk{w \mapsto 2w}$.
According to the product structure described as in \secref{sec_04}, it is natural to draw $\Bar{\ell_0}:=\partial \ell_0 \cup \ell_0$ as (an annulus)$\times [0,1]$ like \figref{fig_boundaries}.

\parag{Pinching map on the quotient 3-laminations.}
Now we are ready to state the main theorem in this section:
\begin{thm}[Quotient 3-laminations]\label{thm_quotient}
There exists a meek map $\hhatM:\Mf - \vv_f \to \Mg$ such that
\begin{enumerate}[\rm (1)]
\item {\bf Non-principal part:} The restriction $\hhatM|~\Mf- \ell_f \to \Mg-\ell_g$ is very meek. In particular, for any non-principal leaf $\ell$ of $\Mf-\ell_f$, $\ell$ and $\hhatM(\ell)$ are hyperbolic 3-manifolds of the same kind.
\item {\bf Principal part:} The complementary restriction $\hhatM|~\ell_f - \vv_f \to \ell_g$ is meek. In particular, $\ell_f, \ell_f-\vv_f$ and $\ell_g$ are all homeomorphic. 
\end{enumerate}
\end{thm}
Thus the inner topology of $\Mf$ is leafwise preserved as $f \to g$. 
(Note that this theorem is only about the inner structure of $\Mf$; not about $\overline{\Mf}$.) 
By (1), the global laminar structure is preserved at least in the non-principal part. 

As to be expected by the observation of \secref{sec_04}, we would find the major topological difference of $\Mfbar$ and $\Mgbar$ in the conformal boundaries, in particular, in the lower ends. This is the main topic of the next two sections.

\parag{Proof of \thmref{thm_quotient}: Existence.}
We construct such a map $\hhatM$ from the pinching semiconjugacy $\hhatH: \Hf-\VV_f \to \Mg$.
Take $\xtil$ and $\ytil$ in $\Hf-\VV_f$ with $\ytil=\fhat^n(\xtil)$ for some $n$. 
Then $\hhatH(\ytil)=\hhatH \cc \fhat^n(\xtil)=g^n \cc \hhatH (\xtil)$ and thus $[\hhat_\h(\xtil)]=[\hhat_\h(\ytil)]$ in the quotient $\Mg$. This implies that a map $\hhatM: \Mf-\vv_f \to \Mg$ is induced from $\hhatH:\Hf-\VV_f \to \Hg$ by defining $\hhatM([\xtil]):=[\hhatH(\xtil)]$. 

Note that both $\Mf-\vv_f$ and $\Mg$ are leaf spaces since $\HH_f-\VV_f$ and $\Hg$ are. Thus the map $\hhatM: \Mf-\vv_f \to \Mg$ is a surjective and continuous map between leaf spaces. 
To show that $\hhatM$ is meek, we divide the map into the non-principal part $\hhatM|_{\Mf- \ell_f}$ and the principal part $\hhatM|_{\ell_f-\vv_f}$, and check their meekness separately.
(Note also that $\Mf- \ell_f$ is a leaf space and $\ell-\vv_f$ is a leaf space with one path-connected component.)

\parag{(1): Non-principal part.}
Now we can describe the topology of $\Mg$ by pinching $\Mf$. 
Take any non-principal leaf $\ell'$ of $\Mg$. 
By \thmref{thm_Hyp^3}, its preimage $\ell:=\hhatM^{-1}(\ell')$ must be also a non-principal leaf, and by definition, the map $\hhatM:\ell \to \ell'$ just pinches $\II_f^\h/\fhat$-components of $\ell$ to $\II_g^\h/\ghat$-components of $\ell'$.
Since $\II_f^\h$ has no periodic component in the non-principal leaves, each $\II_f^\h/\fhat$-component in $\ell$ is topologically a star-like graph times $\R_+$ where $\R_+$ corresponds to the direction of the vertical extension from $\Af$ to $\Hf$.
(cf. \figref{fig_deg_Hf}.) 
When $\ell$ is isomorphic to $\Hyp^3$, so is $\ell'$ since $\hhatM:\ell \to \ell'$ sends the vertical fibers to the vertical fibers. 
When $\ell$ is a solid torus, it is the quotient of a non-principal periodic leaf $L^\h$ in $\Hf$. 
Since $\hhat_\h$ is a semiconjugacy, $\ell'$ is also the quotient of non-principal periodic leaf $\hhat_\h(L^\h)$, thus a solid tori.

\begin{figure}[htbp]
\begin{center}
\includegraphics[width=.60\textwidth]{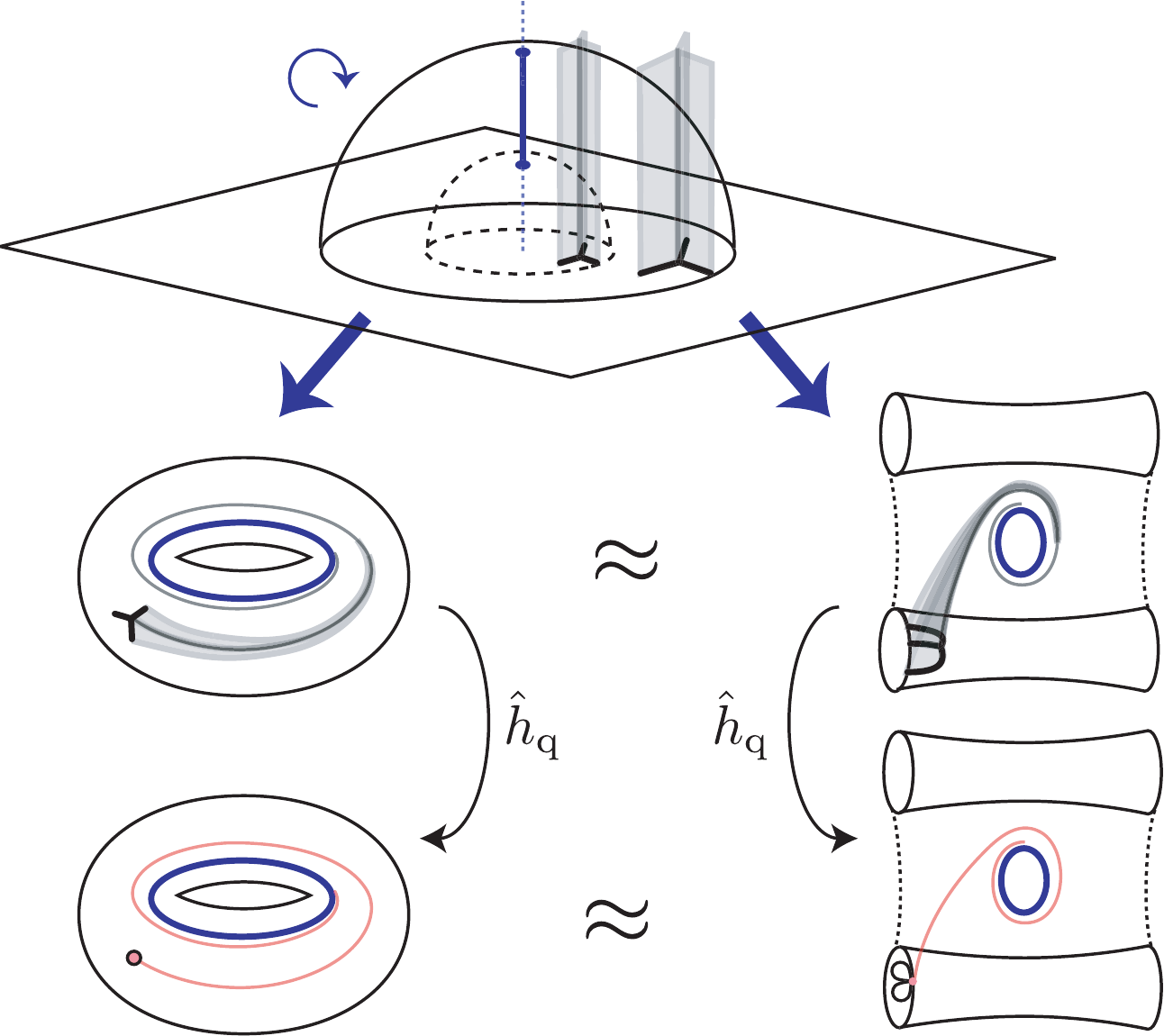}
\end{center}
\caption{A non-principal solid torus in $\Mf$ is pinched to a solid torus in $\Mg$. The blue curves show the core curves in the solid tori. }
\label{fig_solid_tori}
\end{figure}

\parag{Remark.}
This argument would be more justified by considering the \textit{leafwise homotopy}:
When $\ell'$ is a non-principal leaf in $\Mg$, the leaf $\ell=\hhatM^{-1}(\ell')$ is homotopy equivalent to $\ell'$. Thus $\ell$ and $\ell'$ must be the same type $\simeq \Hyp^3$ or solid tori.

When $\ell$ is a solid torus, one can see the vertical fibers of $L^\h$ in $\ell$ as a one-dimensional foliation with fibers coming from the ideal boundary and coiling around the axis of $\ell$. 
(See \figref{fig_solid_tori}.)
Even in this case, the pinching map $\hhatM:\ell \to \ell'$ does not change the topology and $\ell'$ is still a solid torus.

\begin{figure}[htbp]
\begin{center}
\includegraphics[width=.70\textwidth]{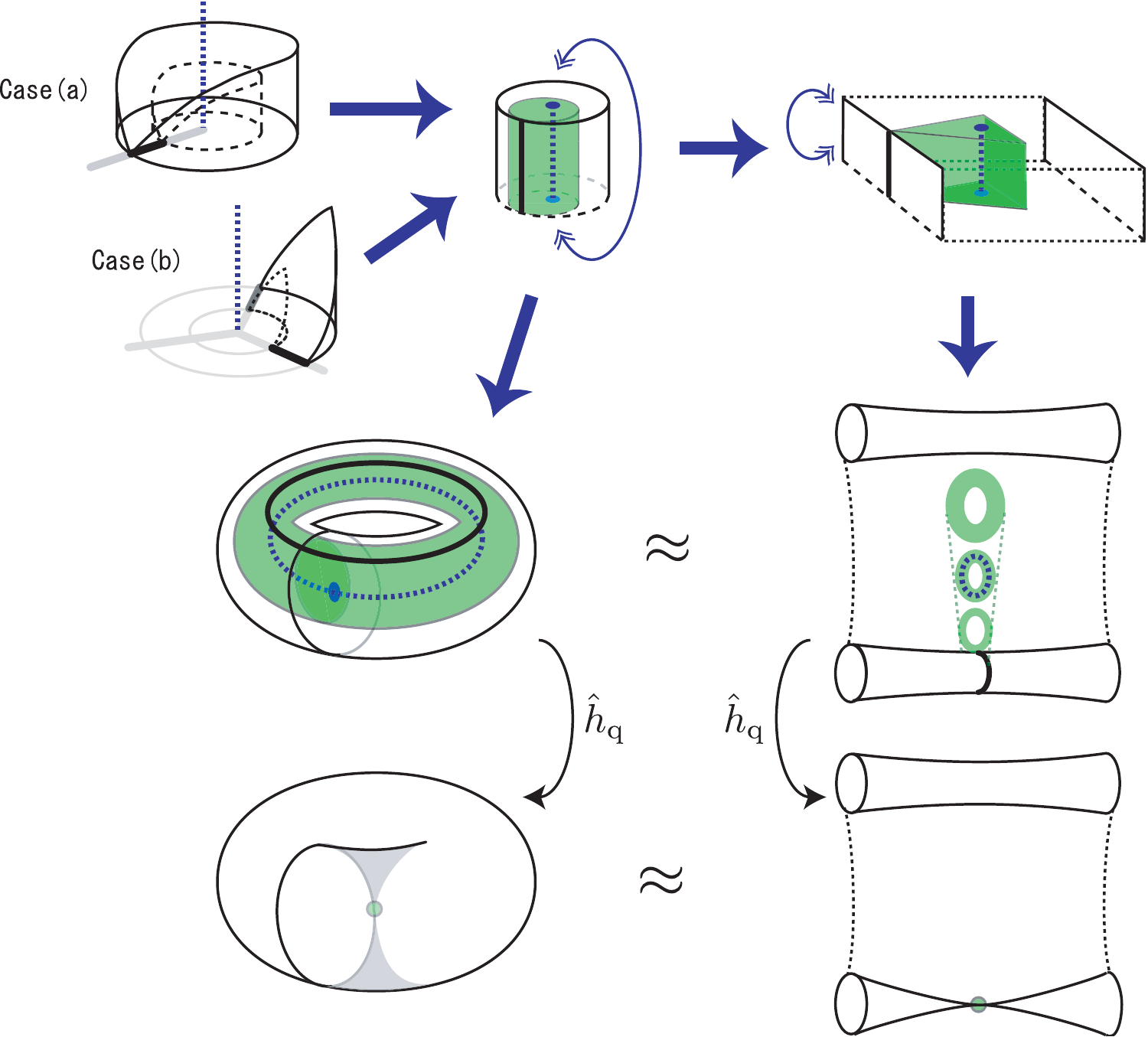}
\end{center}
\caption{A topological picture of the degeneration of the principal leaf $\ell_f$ to $\ell_g$ by $\hhatM$. The valleys are shown in green, and the blue dotted curves show the core curves. }
\label{fig_deg_valley}
\end{figure}

\parag{(2): Principal part.} 
By \thmref{thm_extension_principal_leaves}(1), we have $\hhatM^{-1}(\ell_g)=\ell_f-\vv_f$. 
As \figref{fig_deg_valley} indicates, the valley $\vv_f$ and $\ell_f$ are topological solid tori sharing a circle at their boundary. In particular, $\ell_f-\vv_f$ is again topologically a solid torus $\D \times \T$. 
(Note that the ideal/conformal boundaries are not included in this argument.)
The principal leaf $\ell_g$ is a rank one cusp, but it is also homeomorphic to $\D \times \T$. 
\QED

\parag{Remarks on \figref{fig_deg_valley}: Generation of the cuspidal part.}
On the left of \figref{fig_deg_valley}, a leafwise view with ideal boundaries, the valley $\vv_f$ is shown as a solid torus. On the right, quasi-Fuchsian-like view with conformal boundaries, $\vv_f$ is shown as (an annulus)$\times$(an open interval). 
By the pinching map $\hhatM$, the valley $\vv_f$ in $\ell_f$ is pushed away to the cuspidal part of the lower end.

The lifted Julia set in $\partial \ell_f$ or $\partial \ell_f$

The $\II_f^\h/\fhat$-components in $\ell_f$ should be also drawn as in \figref{fig_solid_tori}. (We abbreviate it for simplicity.)  
On the other hand, the $\II_g^\h/\ghat$-components in $\ell_g$ are a little different from those in \figref{fig_solid_tori}, 
which are coiling around the core curve.
Since $\ell_g$ do not have the core curve, 
any $\II_g^\h/\ghat$-component in $\ell_g$ escapes to the cuspidal part of the lower end. 
We may consider that the coiling $\II_f^\h/\fhat$-components are pushed away with the valley to the cuspidal part of $\Bar{\ell_g}$ by $\hhatM$.

Note also that the lower end may be much more complicated in general. (cf. The Rabbits, see \figref{fig_ell_4_rabbits}.) 

\subsection{Extending the pinching map to the Kleinian laminations}

Next we extend the meek pinching map $\hhatM: \Mf-\vv_f \to \Mg$ to their conformal boundaries. 
We need some more notation:

\parag{Degenerating arcs in the conformal boundary.}
To describe the degeneration of the conformal boundary $\partial \Mf=\Ff/\fhat$, the set $\II_f$ of degenerating arcs again plays an essential role.

Since we have to take a quotient in the Fatou set $\FF_f$, we need to cut out the part of $\II_f$ in the Julia set. We set
$$
\ZZ_f \dee \II_f \cap \Ff.
$$
Since $\Bf$ does not intersect with $\II_f$, the set $\ZZ_f$ is a subset of $\Df$. 
(Recall that $\Ff=\Bf\sqcup\Df$, defined in \secref{sec_04}.)
Moreover, $\fhat_\a$ acts on $\ZZ_f$ homeomorphically and properly discontinuously. Thus the quotient 
$$
\vz_f \dee \ZZ_f/\fhat_\a
$$
lives in the lower end $\Sf=\Df/\fhat_\a$.
We also call the path-connected components of $\ZZ_f$ or $\vz_f$ (``$\ZZ_f$-components" or ``$\vz_f$-components" for short) the \textit{degenerating arcs} in $\Df$ or $\Sf$. 

Since the restriction $\hhatA|~ \Ff-\ZZ_f \to \Fg$ of $\hhatA: \Af-\IO_f \to \Ag$ is a conjugacy by Theorems \ref{thm_natural_extension} and \ref{thm_affine_part}, we have a homeomorphism $\hhatM:\partial \Mf - \vz_f \to \partial \Mg$. Hence:

\begin{thm}[Pinching map on the Kleinian laminations]\label{thm_Kleinian}
The meek pinching map $\hhatM:\Mf - \vv_f \to \Mg$ extends to a meek map $\hhatM:\Bar{\Mf} - \vv_f  \cup \vz_f \to \Bar{\Mg}$ between two leaf spaces whose leaves are 3-manifolds with boundaries. 
\end{thm}

\section{Degeneration of the lower ends}\label{sec_09}
In the following two sections, we observe topological changes of the Kleinian lamination as the parameter $c$ of $f_c$ moves from one hyperbolic component to another, via parabolic $g=f_\sigma$. 
This motion corresponds to a pair of degeneration and bifurcation processes of hyperbolic quadratic maps like the Rabbits. 

The main tool is again tessellations associated with a degeneration pair $(f \to g)$.

\parag{Abstract of this section.}
This section is devoted for the degeneration process. 
Our aim is to complete a topological and combinatorial description of the Kleinian laminations associated with parabolic quadratic maps. 
By \thmref{thm_quotient}, it is enough to study the structure of the conformal boundaries. In particular, we may focus on the lower ends by \propref{prop_Bf_Df}.

First we stand on the fact that any parabolic quadratic map $g = f_\s$ has a Case (a) degeneration pair $(f \to g)$. 
To describe the structures of the lower ends $\Sf=\Df/\fhat_\a$ and $\Sg=\Dg/\ghat_\a$ for this $(f \to g)$, we will realize fundamental regions $\QP_f$ and $\QP_g$ of the actions $\fhat_\a:\Df \to \Df$ and $\ghat_\a:\Dg \to \Dg$ by means of panels in the tessellations. 
Then we will obtain detailed topological descriptions of the lower ends. 
For example, we show that $\Sf$ is compact but $\Sg$ is non-compact (\propref{thm_Sf_Sg}).
As a corollary, the Kleinian lamination $\Bar{\Mf}$ cannot be homeomorphic to $\Bar{\Mg}$.

In addition, we define the \textit{quotient tessellations} $\Tess^\q(f)$ and $\Tess^\q(g)$ supported on $\Sf$ and $\Sg$ respectively. 
Then we will check that the combinatorics of the tiles in the quotient tessellations are the same (\propref{prop_quotient_tiles}). 

Finally, as an example we apply these results to the Cauliflowers.

\parag{Note.}
The structure of $\Sf$ can be described in a different way by applying the argument of \cite[\S\S 3.5]{LM}, but we would need an extra care for $\Sg$, the parabolic case. 
In this section we give a method that can deal with both $\Sf$ and $\Sg$. 


\subsection{Kleinian laminations for parabolic quadratic maps}
So far we have seen that the topologies of the laminations of general parabolic $g$ mostly inherit those of $f$, where $f$ is a perturbation of $g$ such that $(f \to g)$ is a degeneration pair. 
In particular, the leafwise topology of the quotient 3-lamination $\Mf$ is preserved as $f \to g$ by \thmref{thm_quotient}. 
According to the analogy of the quasi-Fuchsian example (\secref{sec_04}), we turn to the study of the Kleinian laminations $\Bar{\Mf}$ and $\Bar{\Mg}$ to detect the difference of ${\Mf}$ and ${\Mg}$.
Here we make a thorough study of the combinatorial and topological properties of the lower ends for such a degeneration pair. 

\parag{Assumption (a).}
To simplify the arguments, we make a technical assumption. 
Here is a direct corollary of \cite[Lemma 4.4]{Mi} (cf. \cite[Theorem 1.1]{Ka2}):
\begin{prop}\label{prop_assumption_a}
For any parabolic quadratic map $g$, there exists a Case (a) degeneration pair $(f \to g)$, that is, $q=q'$ and $l=l'$.
\end{prop}
In the rest of this section we stand on this fact and assume that $(f \to g)$ is Case (a). (We will refer this assumption as ``Assumption (a)".) Note that in the tessellation $\Tess(f)$ the degenerating edge of $\Tf(\theta,m, +)$ is shared by $\Tf(\theta,m, -)$ under this assumption. (See \cite[Proposition 3.2]{Ka3} and \figref{fig_rabbit_tess}.) 

The following argument also works for Case (b) (the case of $q=1<q'$) with a slight modification. But we will deal with Case (b)  in the next section as a bifurcation from $g$. 

\subsection{Fundamental region of the lower ends}

\parag{Tessellations over $\bs{\Df}$ and $\bs{\Dg}$.}
To describe the topology of $\Sf$ and $\Sg$ in detail, we construct fundamental regions of them in $\Df$ and $\Dg$ by means of tessellation. 

Let us check some basic facts on the tessellations $\Tess^\a(f)$ and $\Tess^\a(g)$ on the affine laminations. 
A \textit{tile} in $\Tess^\a(f)$ is defined by 
$$
T^\a_f(\thetahat, m, \ast) \dee \iota_f\kakko{\Tf(\thetahat,m , \ast)}
$$ 
where $T_f(\thetahat,m , \ast)$ is a tile in $\Tess^\n(f)$ and $\iota_f:\Afn \to \Af$ sends the objects in the natural extension to the universal setting (\secref{sec_06}). 
By $\Pif^\a(\thetahat, \ast)$ we denote the union of all possible tiles of the form $\Tf^\a(\thetahat, m , \ast)$ with $m \in \Z$. 
Now the set $\Df-\ZZ_f$ is tessellated by the tiles and the panels of $\Tess^\a(f)$, since 
$$
\pihat_f^{-1}(\Kfc-I_f) \ee \Df-\II_f \ee \Df-\ZZ_f.
$$
We apply similar notations for the lifted tessellation of $\Dg = \pihat_g^{-1}(\Kgc)$.

Recall that the set $\Thetahat=\Thetahat_{(f \to g)}$ is an invariant subset of the lifted angle doubling map $\delhat: \TThat \to \TThat$. 
For any $\thetahat \in \Thetahat$, tiles and panels satisfy
$$
\fhat_\a(\Tf^\a(\thetahat, m , \ast))\ee
\Tf^\a(\delhat(\thetahat), m+1 , \ast)
~~~\text{and}~~~
\fhat_\a(\Pif^\a(\thetahat, \ast))\ee
\Pif^\a(\delhat(\thetahat), \ast).
$$
There are only finitely many cyclic panels in the tessellation, which are associated with cyclic angles in $\Thetahat$.

\parag{Cyclic panels vs. non-cyclic panels.}
The fundamental regions we will construct are based on the dynamics of the panels. 
We will divide the fundamental region into two parts: The \textit{cyclic part} $\Qf(\pm)$ and the \textit{non-cyclic part} $\Pf(\pm)$, which correspond to the cyclic panels and the non-cyclic panels respectively.
To construct them, let us start with some definitions on special angles.

\parag{A cyclic angle $\boldsymbol{\theta_0^+}$ (downstairs).} 
We first fix a cyclic angle $\theta_0^+$ in $\Theta$ so that we have a basepoint of the eventually periodic dynamics $\delta: \Theta \to \Theta$. 
We choose $\theta_0^+$ as follows: 
 Let $\beta_0$ be the parabolic periodic point in $O_g$ that is on the boundary of the Fatou component containing the critical point $0$. In the type $\Theta(\beta_0)$, there exist two angles $\theta_0^+$ and $\theta_0^-$ with representatives $\theta_0^+<\theta_0^- \le \theta_0^+ +1$ such that $\delta^{\lbar} (\theta_0^\pm)=\theta_0^\pm$ and $\partial \Pi_g(\theta_0^\pm,\pm)$ contains $0$. (Such angles are characterized by  the \textit{critical sector} in \cite[\S 2.3]{Ka3}). 

\parag{Examples.}
For the Cauliflowers, we set $\theta_0^+ = 0$ and $\theta_0^- = 1 \equiv 0$. 
For the Rabbits, we set $\theta_0^+ = 4/7$ and $\theta_0^- = 8/7 \equiv 1/7$. 
For the Airplains, we set $\theta_0^+ = 5/7$ and $\theta_0^- = 9/7 \equiv 2/7$.

Now we have:

\begin{prop}\label{prop_panels}
Every panel in $\Tess(f)$ or $\Tess(g)$ is uniquely parameterized by an angle that eventually lands on $\theta_0^+$ and signature $\ast=\pm$. In addition, under Assumption (a), the cyclic panel of angle $\theta_0^+$ and signature $\ast$ consists of the tiles with addresses of the form $(\theta_0^+, \mu l, *)~(\mu \in \Z)$. 
\end{prop}

\begin{pf}
Since the combinatorics of tessellations are the same for $f$ and $g$ (\secref{sec_05}), it is enough to show the statement for $g$. 

If $q=q'=1$ and $l>1$ (like the Airplains), the angles $\theta_0^+$ and $\theta_0^-$ belong to the different cycles under angle doubling but $\Pig(\theta_0^+, \ast)=\Pig(\theta_0^-,\ast)$ by definition. Indeed, each panel has two angles in $\Theta = \Theta_{(f \to g)}$ that eventually land on either $\theta_0^+$ or $\theta_0^-$. (See ``Remark on angles and levels" in \cite[\S\S 3.2]{Ka3}.)
Otherwise $\theta_0^+$ and $\theta_0^-$ belong to the same cycle. Each panel has a unique angle in $\Theta$ that eventually lands on $\theta_0^+$.
Hence in both cases, every panel is uniquely parameterized by an angle that eventually lands on $\theta_0^+$ and signatures $\ast=\pm$.

By definition of tiles in \cite[\S 3]{Ka3} and Assumption (a), the panel $\Pi_g(\theta_0^+,*)$ consists of the tiles with addresses of the form $(\theta_0^+, \mu l, *)~(\mu \in \Z)$. 
\QED
\end{pf}

\parag{Fundamental angles (upstairs).}
Next we consider the ``fundamental angles" of the dynamics $\delhat: \Thetahat \to \Thetahat$ upstairs. Let $\thetahat^+=(\theta_{-n}^+)_{n \ge 0 } \in \Thetahat$ denote the cyclic lift of $\theta_0^+$; that is, $\theta_{-k\lbar}^+=\theta_0^+$ for all $k \ge 0$. Then every cyclic panel in $\Tess^\a(g)$ lands on either $\Pig^\a(\thetahat^+, +)$ or $\Pig^\a(\thetahat^+, -)$ .

We say $\hat{\omega}=\bon{\omega}$ is an element of $\Omega \subset \TThat$ if $\omega_0=\theta_0^+$ and $\omega_{-n} \neq \theta_0^+$ for any $n \ge 1$. Note that $\Omega$ is a compact subset of $\Thetahat$. Let $\thetahat = \bon{\theta} \in \Thetahat$ be a non-cyclic angle such that the forward iteration $\skakko{\delta^m(\theta_0)}_{m>0}$ eventually lands on $\theta_0^+$. Then there exist unique $n \in \Z$ and $\hat{\omega} \in \Omega$ such that $\delhat^n(\thetahat)= \hat{\omega}$ and $\ghat^n(\Pig^\a(\thetahat,\ast))= \Pig^\a(\hat{\omega},\ast)$. Since the panels in $\Tess^\a(f)$ also has the same properties, we have:

\begin{prop}\label{prop_panels_upstairs}
The quotient of each non-cyclic panel is uniquely represented by a panel with an angle in $\Omega$. 
If $\thetahat= \bon{\theta}$ is $\thetahat^+$ or in $\Omega$ as above, then $\theta_0=\theta_0^+$. Thus the panel of angle $\thetahat$ and signature $\ast$ upstairs is projected onto that of angle $\theta^+_0$ and the signature $\ast$ downstairs.  
Moreover, this panel upstairs consists of tiles with address of the form $(\thetahat, \mu l, \ast)~(\mu \in \Z)$. 

\end{prop}

\begin{figure}[htbp]
\centering{
\vspace{0cm}
\includegraphics[width=0.8\textwidth]{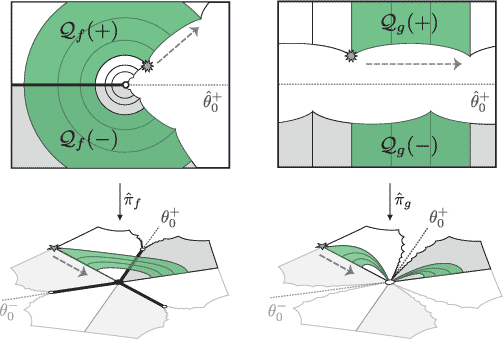}
}
\caption{A caricature of $\Qf(\ast)$ and $\Qg(\ast)$ (in green) for $q'=3$ under Assumption (a). The dashed arrow show the actions of $f^{\lbar}$ and  $g^{\lbar}$ downstairs, and those of $\fhat^{\lbar}$ and $\ghat^{\lbar}$ upstairs. The stars show the critical points $z=0$ downstairs, and their lifted images upstairs. The blue doted lines are the external rays.}
\label{fig_Qf_Qg}
\end{figure}

\parag{Fundamental regions: Cyclic parts.}
Let us define fundamental regions of the cyclic panels in $\Df$ and $\Dg$. One can find such regions in $\Pif^\a(\thetahat^+, \ast)$ and $\Pig^\a(\thetahat^+, \ast)$ as follows: Set 
$$
\Qf(*) \dee \bigcup_{\mu =0}^{q-1} 
\overline{\Tf^\a(\thetahat^+,  \mu l, *)}  
~~~\text{and}~~~
\Qg(*) \dee \bigcup_{\mu =0}^{q-1} \Tg^\a(\thetahat^+, \mu l, *)  
$$
for $\ast = +$ or $-$. Here the closures of tiles in $\Qf(*)$ are just tiles plus their degenerating edges on $\IO_f$. (That is, we take the leafwise closures.) 
By Assumption (a), $\Qf(+)$ and $\Qf(-)$ share their degenerating edges as in \figref{fig_Qf_Qg} (upper left). 
Note that for $f$ (resp. $g$) the principal leaf containing the lifted external ray of angle $\thetahat^+$ also contains $\Qf(\ast)$ (resp.  $\Qg(\ast)$). 

By properties of the panels of angle $\thetahat^+$, one can easily check that every orbit in the cyclic panels passes through those set just once, except some orbits landing on the boundaries. 


\begin{figure}[htbp]
\centering{
\vspace{0cm}
\includegraphics[width=0.8\textwidth]{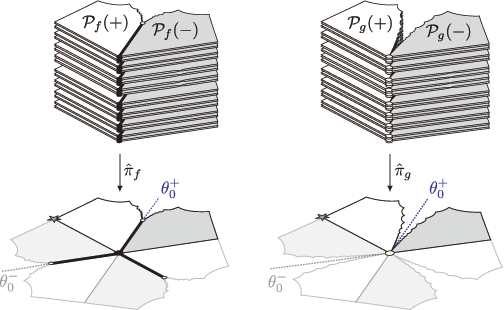}
}
\caption{A caricature of $\Pf(\ast)$ and $\Pg(\ast)$ for $q' = 3$. Note that they are projected onto $\overline{\Pif(\theta_0^+,\ast)}-\skakko{\gam_f(\theta_0^+)}$ and $\Pig(\theta_0^+,\ast)$ downstairs. The stars show the critical point $z=0$. The doted lines are the external rays of angles $\theta^\pm_0$.}\label{fig_Pf_Pg}
\end{figure}

\parag{Fundamental regions: Non-cyclic parts.}
Next one can choose fundamental regions of the non-cyclic parts in $\Df$ and $\Dg$ as follows: For each $\ast=\pm$, set 
$$ 
\Pf(\ast) \dee \bigcup_{\omehat \in \Omega} 
\kakko{\overline{\Pif^\a(\omehat, \ast)}-\{\gam_f^\a(\omehat)\}}
~~\text{and}~~
\Pg(*) \dee \bigcup_{\omehat \in \Omega} \Pig^\a(\omehat,\ast),
$$
where $\gam_f^\a(\omehat):=\iota_f(\gam_f(\omehat)) \in \Jf$. 
By Assumption (a), panels in $\Pf(+)$ share the degenerating edges with those in $\Pf(-)$ (\figref{fig_Pf_Pg}, upper left).
Note that each $\overline{\Pif^\a(\omehat, \ast)}$ contains one point of $\pihat_f^{-1}(\al_0)$, since $\pihat_f:\Af \to \Chat$ projects $\Pi_f^\a(\omehat, \ast)$ over the panel $\Pi_f(\theta^+_0, \ast)$ downstairs, and this panel contains $\al_0$ on its boundary.
Note also that for each $\ast=\pm$, topologically the interior of $\Pf(*)$ is the product of the disk $\Pif(\theta_0^+,\ast)^\cc$ downstairs and the Cantor set $\Omega$. (The same holds for the interior of $\Pg(*)$.) 

By properties of the panels with angles in $\Omega$, every orbit in non-cyclic panels passes through those set just once except the orbits landing on the boundaries of panels. Indeed, such exceptional orbits land on these boundaries at most finite times. 

By those construction of $\Qf(\ast)$ and $\Pf(\ast)$, we conclude:
\begin{prop}[Fundamental regions]\label{prop_fundamental_region}
The sets 
$$
\QQ\PP_f \dee \bigcup_{\ast=\pm} \kakko{\, \Qf(\ast) \, \cup \, \Pf(\ast) \, }
~~\text{and}~~ 
\QQ\PP_g \dee \bigcup_{\ast=\pm} \kakko{\, \Qg(\ast) \, \cup \, \Pg(\ast) \, }
$$
are fundamental regions of $\Sf=\Df/\fhat_\a$ and $\Sg=\Dg/\ghat_\a$ in the dynamics on $\Df$ and $\Dg$. 
In particular, $\hhatA$ sends $\QQ\PP_f - \ZZ_f$ homeomorphically onto $\QQ\PP_g$.
\end{prop}

\subsection{Topology of the lower ends}
Let us describe the topology of $\Sf=\Df/\fhat$ and $\Sg=\Dg/\ghat$ by means of the fundamental regions $\QQ\PP_f$ and $\QQ\PP_g$.
For $\ast = \pm$, we set $\Qf^\ast:=\Qf(\ast)/\fhat$, $\Qg^\ast:=\Qg(\ast)/\ghat$, $\Pf^\ast:=\Pf(\ast)/\fhat$ and $\Pg^\ast:=\Pg(\ast)/\ghat$ for simplicity. Then it is clear by construction that 
$$
\Sf \ee \QQ\PP_f/\fhat \ee \bigcup_{\ast=\pm} \kakko{\, \Qf^\ast \, \cup \, \Pf^\ast \, }
~~\text{and}~~ 
\Sg \ee \QQ\PP_g/\ghat \ee \bigcup_{\ast=\pm} \kakko{\, \Qg^\ast \, \cup \, \Pg^\ast \, }.
$$

\parag{Relations between $\boldsymbol{\PP}$ and $\boldsymbol{\QQ}$.}
The lower ends $\Sf$ and $\Sg$ have both differences and common properties. 
Let us first check some common properties to $\Sf$ and $\Sg$:
\begin{prop}[$\PP$ winds around $\QQ$]\label{thm_Qf_Pf}
In the lower ends $\Sf$ and $\Sg$, the subsets $\Qf^\ast$, $\Pf^\ast$ and $\Qg^\ast$, $\Pg^\ast$ $(\ast=\pm)$ satisfy the following:
\begin{enumerate}[\rm (1)]
\item For each $\ast=\pm$, 
$$
\Qf^\ast \cup \Pf^\ast \ee \overline{\Pf^\ast}
~~~\text{and}~~~ 
\Qg^\ast \cup \Pg^\ast \ee \overline{\Pg^\ast}.
$$
\item $\overline{\Pf^+} \cup \overline{\Pf^-} = \Sf$ and  $\overline{\Pg^+} \cup \overline{\Pg^-} = \Sg$.
\end{enumerate}
\end{prop}

\paragraph{Proof.}
Fix a signature $\ast=\pm$ and take any point in $\xhat \in \Qf(\ast) \subset \Df$. We first show that for every $\omehat \in \Omega$ there exists a sequence $\xhat_n \in \Pif^\a(\omehat, \ast)$ with $[\xhat_n] \to [\xhat]$ in $\Sf$, thus $\Qf^\ast \subset \overline{\Pf^\ast}$. Let us divide $\Pif^\a(\omehat,\ast)$ into the following subsets parameterized by $n \in \Z$:
$$
Q_{-n}~\dee \bigcup \, \skakko{ \,
\overline{\Tf^\a(\omehat, \mu l  - n \lbar, \ast)}
: 0 \le \mu \le q-1
\,}
$$
so that $\pihat_f(\Qf(\ast))=\pihat_f(\fhat^{n\lbar}(Q_{-n}))$. Since $\delhat^{n \lbar}(\omehat) \to \thetahat^+$ as $n \to +\infty$, we can approximate any point $\xhat \in \Qf(\ast)$ by points in $\fhat^{n \lbar}(Q_{-n})$. Thus we can take such $\xhat_n \in Q_{-n}$ as desired. Note that $[Q_{-n}]$ (tiles with deeper levels) are folding to $\Qf^\ast$ as $n \to \infty$.  

Next we take any point $x$ in $\overline{\Pf^\ast}-\Pf^\ast$. Then there exists a sequence $\xhat_n \in \Pf(\ast)$ which has no limit point in $\Pf(\ast)$ but $[\xhat_n] \to x$ in the quotient. By definition of $\Pf(\ast)$, such sequence must accumulate on $\{\gam_f^\a(\omehat):\omehat \in \Omega \} \subset \Jf$. It implies that $\xhat_n$ are contained in tiles with deeper levels as $n \to \infty$. Thus $[\xhat_n]$ must accumulate on $\Qf^\ast$ and it follows that $x \in \Qf^\ast$. This completes the proof of (1) for $f$. 

Since $\Qf(\ast) \cup \Pf(\ast)$ is a fundamental region for $\Sf$, (2) immediately follows from (1). These arguments clearly works if we replace $f$ by $g$. \QED

\begin{figure}[htbp]
\centering{
\includegraphics[width=0.9\textwidth]
{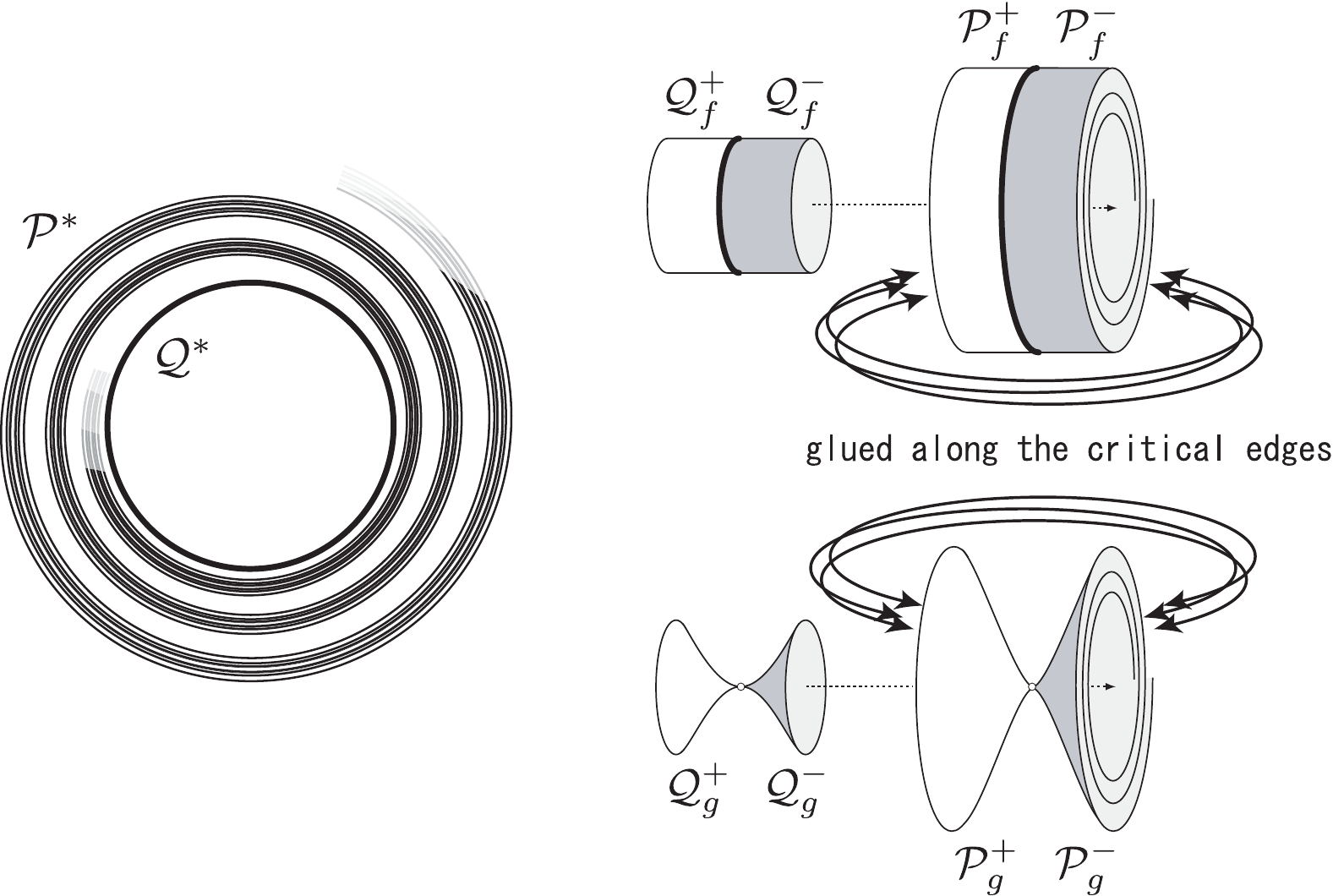}
}
\caption{Laminar structures of $\Sf=\bigcup_{\ast = \pm} (\Pf^\ast \cup \Qf^\ast)$ (or $\Sg=\bigcup_{\ast = \pm}(\Pg^\ast \cup \Qg^\ast)$).}
\label{fig_winding}
\end{figure}

\paragraph{}
By this proof it follows that any quotient panel $[\Pif^\a(\omehat, \ast)]$ in $\Pf^\ast$ accumulates on $\Qf^\ast$ in $\Sf$. 
On the left of \figref{fig_winding} a one-dimension down caricature of this situation is shown. 
The central thick circle is $\Qf^\ast$ with laminated quotient panels in $\Pf^\ast$ winding around. 
The same is true for $g$. 

\parag{Compactness of the lower ends.}
A significant difference between $\Sf$ and $\Sg$ is the following:
\begin{prop}[Compactness]\label{thm_Sf_Sg}
The lower end $\Sf$ is compact but $\Sg$ is non-compact.
\end{prop}
\begin{pf}
Take any sequence $\xhat_n$ in $\Pf(\ast)$ escaping from all compact subsets of $\Pf(\ast)$. By construction of $\Pf(\ast)$, the sequence $[\xhat_n]$ must accumulate on $\Qf^\ast$ by \propref{thm_Qf_Pf}. Since $\Qf^\ast$ is compact, we have the compactness of $\Sf$. On the other hand, similar sequence $\xhat_n$ in $\Pg(\ast)$ may get closer and closer to $\Qg^\ast$ in the quotient, but it may not accumulate since $\Qg^\ast$ is non-compact. In fact, we have such an escaping sequence $[\xhat_n] \in \Sg$ by taking $\xhat_n \in \Pg(\ast)$ so that the heights of $z_n=\pihat_g(\xhat_n)$ are unbounded with respect to a fixed Fatou coordinate of $\Kgc =D_g$. \QED
\end{pf}

Here is a direct corollary:
\begin{cor}\label{cor_Kleinian_not_homeo}
The Kleinian lamination $\Bar{\Mf}$ is not homeomorphic to $\Bar{\Mg}$. More generally, there is no pair of hyperbolic and parabolic quadratic maps with topologically the same Kleinian laminations.
\end{cor}

\subsection{Quotient tessellations.}
Next we focus on how the boundaries of the fundamental regions $\QQ\PP_f$ and $\QQ\PP_g$ are glued by taking quotient. 
Since it is determined by the connection of the panels in the quotient, we introduce the \textit{quotient tessellations} supported on $\Sf$ and $\Sg$. 
(The quotient tessellations will also play crucial roles in the next section.)
A \textit{quotient tile} in the lower ends is an equivalent class $[T]$ of a tile $T$ in $\Tess^\a(f)$ or $\Tess^\a(g)$. We set 
\begin{align*}
\Tess^\q(f) &\dee \skakko{[T] ~:~ T \in \Tess^\a(f)} \\
\Tess^\q(g) &\dee \skakko{[T] ~:~ T \in \Tess^\a(g)}. 
\end{align*}
Recall that the boundaries of tiles in $\Tess^\a(f)$ and $\Tess^\a(g)$ consist of equipotential, critical, and degenerating edges. 
For a quotient tile $[T]$, we define its equipotential, critical and degenerating edges by the equivalent classes of the corresponding edges of $T$. (Of course there is no degenerating edges in tiles of $\Tess^\q(g)$.)

Note that quotient panels also make sense. 
We first study how the quotient tiles/panels in $\Sf$ and $\Sg$ share their edges. We claim:

\begin{prop}[Quotient tiles and their addresses]\label{prop_quotient_tiles}
The quotient tiles have the following properties:
\begin{enumerate}[\rm (1)]
\item The tiles in $\Tess^\q(f)$ and $\Tess^\q(g)$ tessellates $\Sf-\vz_f$ and $\Sg$ respectively. 
\item The set of degenerating arcs $\vz_f$ is the union of degenerating edges for all quotient tiles in $\Tess^\q(f)$.
\item Each tile $[T] \subset \Sf$ (resp. $\Sg$) is uniquely represented by a tile $T'$ in the fundamental region $\QPf$ (resp. $\QPg$) respectively. Thus $[T]$ has a well-defined address given by that of $T'$. 
\item The map $\hhatM: \Sf-\vz_f \to \Sg$ sends the quotient tiles to the quotient tiles without changing the addresses.
\item Any tiles in $\Tess^\q(f)$ share their equipotential or critical edges iff so do the tiles in $\Tess^\q(g)$ with the same addresses.
\end{enumerate}
\end{prop}

\begin{pf}
Property (1) comes from the fact that the tiles in $\Tess^\a(f)$ and $\Tess^\a(g)$ tessellates $\Df-\ZZ_f$ and $\Dg$ respectively. 
Since $\ZZ_f$ in $\Df$ is the union of degenerating edges in the tessellation $\Tess^\a(f)$, we have Property (2) by taking quotient.
Property (3) holds since $\QPf$ and $\QPg$ are fundamental regions that are based on the union of panels with angles $\{\thetahat^+\} \cup \Omega$. 
To check properties (4) and (5), recall that $\hhatA: \Df -\ZZ_f \to \Df$ conjugates $\fhat_\a$ and $\ghat_\a$ and preserve the address of tiles in $\Tess^\a(f)$ and $\Tess^\a(g)$. Then (4) and (5) are natural conclusions. 
\QED
\end{pf}

\parag{Connection of the edges on $\bs{\QQ}$ and $\bs \PP$.} 
Let us remark a little more about the connection of edges. 
The connection of the equipotential edges are trivial, since they are uniquely glued in the panels. 
On the other hand, as indicated on the right of \figref{fig_winding}, the critical edges in the fundamental regions $\QQ\PP_f$ and $\QQ\PP_g$ are glued in a complicated way according to the dynamics. 
However, by \propref{prop_quotient_tiles}(4), the connection of the equipotential or critical edges are preserved as $f \to g$. 

Note that the topological difference between $\Sf$ and $\Sg$ is given by the existence of the degenerating arcs $\vz_f$ (the black thick curves in \figref{fig_winding}, upper right) for $\Sf$. 
The degenerating arcs will also play an important role in bifurcation processes, which we will deal with in the next section.

\subsection{Example: The Cauliflowers}\label{subsec_cauli}
Here we give a brief summery of our results applied to the Kleinian laminations of the Cauliflowers $(f \to g)$. 
(The Airplanes are described in a similar way.)

Recall that the Cauliflowers represent the motion of $f_c$ starting with $f_0(z)=z^2$ and moving toward the parabolic $g(z)=z^2+{1/4}$. 
For $f=f_c$ with $0<c<1/4$, \corref{cor_hyp_center} and \thmref{thm_prod_str_of_M} imply:
\textit{
\begin{itemize}
\item[\bf (C1)] The Kleinian lamination $\overline{\MM_{f_0}}$ is quasi-isometrically equivalent to $\overline{\MM_{f}}$. 
In particular, $\overline{\MM_{f}}$ has a product structure $\approx \SSS_0 \times [0,1]$.
\end{itemize}
}
Thus it is enough to investigate the difference between $\overline{\MM_{f}}$ and $\overline{\MM_{g}}$.
According to the notation in \secref{sec_08}, we set the principal leaves $\ell_f:=\Lam_f^\h/\fhat$ (a solid torus) and $\ell_g:=\Lam_g^\h/\ghat$ (a rank one cusp). By \thmref{thm_quotient},
\textit{
\begin{itemize}
\item[\bf (C2)] There exists a very meek map from the non-principal parts $\Mf-\ell_f$ onto $\Mg-\ell_g$. 
\item[\bf (C3)] The principal part $\ell_f$ and $\ell_g$ has the same topology homeomorphic to $\D \times \T$.
\end{itemize}
}
Thus an essential difference appears in the conformal boundaries. 
Before we deal with the conformal boundaries, let us remark that we can upgrade (C2) as follows:
\textit{
\begin{itemize}
\item[\bf (C2)'] The non-principal parts $\Bar{\Mf}-\Bar{\ell_f}$ and $\Bar{\Mg}-\Bar{\ell_g}$ of the Kleinian laminations are homeomorphic. 
\end{itemize}
}
\parag{Sketch of the proof of (C2)'.}
In the dynamics downstairs, every non-cyclic panel eventually lands on the panel of angle $1/2$ and then mapped to the invariant panel of angle $0$. It is not difficult to construct a homeomorphism $H: \bigcup_{\ast=\pm}\Pif(1/2, \ast) \to  \bigcup_{\ast=\pm} \Pig(1/2, \ast)$ which extends to the boundaries and agrees with the semiconjugacy $h$. By using dynamics and B\" ottcher coordinate on/outside the Julia set, we can extend the map $H$ to the sphere except the union of panels of angle $0$ so that $H$ is a conjugacy.
This lifts to a conjugacy between the non-principal part $\hat{H}_\a: \Af-\Lam_f^\a \to \Ag-\Lam_g^\a$.
Then $\hat{H}_\a$ also extends to a conjugacy between $\hat{H}_\h: \Hf-\Lam_f^\h \to \Hg-\Lam_g^\h$ by \cite[Lemma 9.2]{LM} or \lemref{lem_extension}. 
By taking quotient we can easily yield a homeomorphism between ${\Mf}-{\ell_f}$ and ${\Mg}-{\ell_g}$. 
It extends to the conformal boundaries since the conjugacy $\hat{H}_\a$ preserves the Fatou set and the Julia set.
\QED

\parag{Conformal boundaries.}
Let $\SSS^f$ and $\SSS^g$ denote the upper ends of $\Mf$ and $\Mg$, and $\SSS_f$ and $\SSS_g$ their lower ends.
By \propref{prop_Bf_Df} and \propref{thm_Sf_Sg}, we have
\textit{
\begin{itemize}
\item[\bf (C4)] The upper ends $\SSS^f$ and $\SSS^g$ are conformally equivalent to Sullivan's solenoidal Riemann surface lamination $\SSS_0$. 
\item[\bf (C5)] The lower end $\Sf$ is compact, but $\Sg$ is non-compact.
\end{itemize}
}
\parag{Principal leaves.}
Let us check the difference between $\Bar{\ell_f}$ and $\Bar{\ell_g}$ in detail. 
By (C4) the leafwise upper ends of $\ell_f$ and $\ell_g$ are conformally the same annulus, say $s_0$. 
Let $s_f$ and $s_g$ denote their leafwise lower ends in $\ell_f$ and $\ell_g$. 
Then $s_f$ is an annulus containing the cyclic part $\QQ_f^+ \cup \QQ_f^-$ of the fundamental region.
On the other hand, $s_g$ consists of two punctured disks which has a cuspidal part $\QQ_g^+ \cup \QQ_g^-$.

\begin{figure}[htbp]
\centering{
\vspace{0cm}
\includegraphics[width=0.55\textwidth]{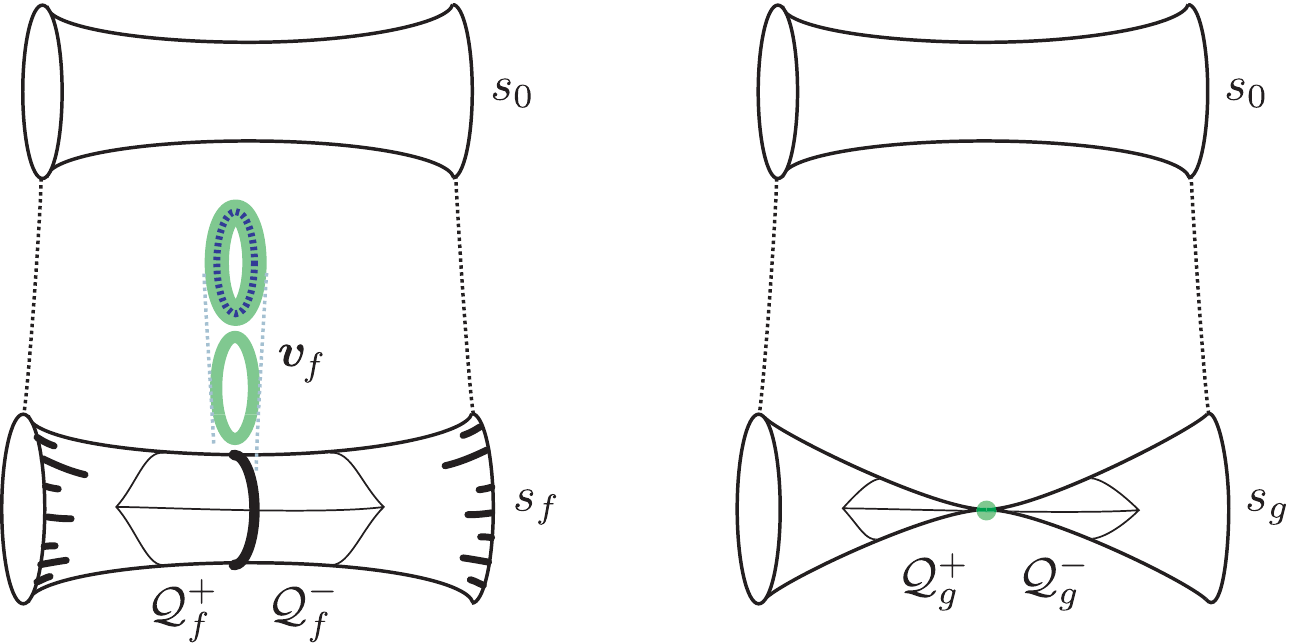}
}
\caption{Caricatures of $\Bar{\ell_f}$ and $\Bar{\ell_g}$. Note that the core curve (the dotted circle) of $\ell_f$ is contained in the valley $\vv_f=\VV_f/\fhat$. There are $\II_f^\h/\fhat$-components (not drawn) coiling around the core curve which are sent by $\hhat_\q$ onto geodesics escaping to the cusp (cf. \figref{fig_solid_tori}).}\label{fig_ell_4_cauliflower}
\end{figure}

In the sense of the product structure given in \thmref{thm_prod_str_of_M}, we may regard $s_0 \cup \ell_f \cup s_f$ as a product $s_0 \times [0,1]$. 
However, the product structure obviously breaks for $\ell_g$ (\figref{fig_ell_4_cauliflower}).
Thus we conclude:
\textit{
\begin{itemize}
\item[\bf (C6)] $\overline{\Mf}=\Mf \cup \partial \Mf$ has a product structure $\approx \SSS_0 \times [0,1]$, but $\overline{\Mg}$ does not. In particular, The Kleinian laminations $\overline{\MM_{f_0}}$ and $\overline{\Mg}$ are not homeomorphic.
\end{itemize}
}

\parag{Remark.}
By the same argument as in \cite[\S 6.4]{LM}, one can show that $\Mg - \ell_g$ has the product structure $\Bar{\Mg} - \Bar{\ell_g} \simeq (\SSS_0-s_0) \times [0,1]$.

\section{Twisting the lower ends}\label{sec_10}
This section is devoted for the bifurcation processes after Case (a) degenerations as in the preceding section.

Suppose that $f=f_c$ varies continuously from one hyperbolic component to another via parabolic $g=f_\s$. 
This degeneration and bifurcation process is represented by a couple of Case (a) degeneration pair $(f_1 \to g_1)$ and Case (b) degeneration pair $(f_2 \to g_2)$ with $g_1=g_2=g$. 
(This cannot happens when $g$ is the root of a primitive copy of the Mandelbrot set, like the Cauliflowers and the Airplanes.)

The aim of this section is to compare the Kleinian laminations for $f_1$ and $f_2$. 
Again by an analogy to the quasi-Fuchsian example (\secref{sec_04}), we continue to study the topological and combinatorial structure of the lower ends. 
Here is a brief summary of this section:
\begin{itemize}
\item 
First we check that the topology of $\MM_{f_c}$ is leafwise preserved as $c$ moves as above (\propref{prop_M1_M2}). So we will focus on describing the difference of the lower ends for $f_1$ and $f_2$.
\item 
We compare the tessellations $\Tess(f_1)$ and $\Tess(f_2)$. Then we introduce the {\it subdivided tessellation} $\Tess'(f_2)$ of $\Tess(f_2)$ so that there is a natural one-to-one correspondence between the tiles in $\Tess(f_1)$ and $\Tess'(f_2)$.
\item 
We define a fundamental region $\QP_{f_2}$ of the action $\fhat_{2\a}:\DD_{f_2} \to \DD_{f_2} $ for the lower end $\SSS_{f_2}=\DD_{f_2}/\fhat_{2\a}$ in the same way as we did for $f_1$. 
Then we check some common properties of $\SSS_{f_1}$ and $\SSS_{f_2}$ (\propref{prop_common_prop}).
\item 
By investigating the connections of the degenerating edges of (subdivided) tiles in the fundamental regions, we will show that the quotient tessellations in $\SSS_{f_1}$ and $\SSS_{f_2}$ have different combinatorics, given by \textit{twisting} and \textit{dislocation} (Theorems \ref{thm_twisting} and \ref{thm_twisting2}).
\item 
Then we return to the case when $f_1$ is in the main cardioid of the Mandelbrot set and give a proof of \thmref{thm_introduction_2}. 
\item 
Finally as an example, we apply these results to the Rabbits.
\end{itemize}

\subsection{Trans-component partial conjugacies.}
Let $f_1, g$ and $f_2$ as above.

\parag{Assumption and notation.} 
Let $p', ~q',~ l'$ be as usual for $g=g_1=g_2$. To distinguish $p, ~q,~ l$'s of $f_1$ and $f_2$, we denote them by $p_i,~ q_i, ~l_i$ for $f_i~(i=1,2)$. 
Since $g$ is a shared boundary point of two distinct hyperbolic components, we assume that $q' \ge 2$. 
By \propref{prop_case_a_case_b} we may assume that the Case (a) degeneration pair $(f_1 \to g_1)$ and the Case (b) $(f_2 \to g_2)$ satisfy the following:
 $$
\left\{ \begin{array}{rl}
q_1&=~q' ~=:~ q  \\
l_1&=~l' ~=:~ l 
          \end{array}
          \right.
~~\text{and}~~         
\left\{ \begin{array}{rl}
q_2&=~1  \\
l_2&=~q'l'~=~ql \ee \lbar. 
          \end{array}
          \right.
$$
Note that we may regard $(f_1 \to g_1)$ as $(f \to g)$ in the preceding section satisfying Assumption (a).

\parag{Examples.}
The Rabbits are typical examples. 
More generally, we can take $f_1$ in the main cardioid and $f_2$ in the center of the $p/q$-satellite limb. 

\parag{}
Note that we always have $\Theta_{f_1}=\Theta_{f_2}=\Theta_g$. 
For simplicity we denote them by $\Theta$ and their common invariant lift to $\TThat$ by $\Thetahat$. 

To simplify the notation, objects of the form $X_{f_i}~(i=1$ or $2)$ will be denoted by $X_i$ if it does not lead confusion. For example, we write $\AAA_1=\AAA_{f_1}$, $I_2=I_{f_2}$, etc.

\parag{Difference in the affine and $\Hyp^3$-laminations.}
Let us roughly compare the topologies of $\AAA_1$ and $\AAA_2$.
By \thmref{thm_affine_part}(1), every non-principal leaf of $\AAA_1$ is pinched to be a non-principal leaf of $\Ag$, then again plumped to be a non-principal leaf of $\AAA_2$. In particular, there is a one-to-one correspondence between the non-principal leaves of $\AAA_1$ and $\AAA_2$, and one can find no significant topological difference. 
On the other hand, each one of the principal leaves in $\AAA_2$ originally comes from $q'$ principal leaves of $\AAA_1$. (cf. \figref{fig_principal_leaves}.)  
Thus the leafwise topology of the affine lamination changes at the principal part when $f_c$ moves continuously from one hyperbolic component to another via a parabolic map. 

One can easily see that the same happens for the $\Hyp^3$-laminations. However, if we turn to the quotient laminations, the leafwise topology is preserved in the following sense.

\parag{Leafwise topology of the quotient 3-laminations.}
Let us compare $\MM_1$ and $\MM_2$, the interiors of the Kleinian laminations. 
By \thmref{thm_quotient}, we have meek pinching maps $\hhat_{1\rm{q}}:\MM_1 - \vv_1 \to \Mg$ and $\hhat_{2\rm{q}}:\MM_2 - \vv_2 \to \Mg$. This immediately implies:
\begin{prop}[Leafwise topology]\label{prop_M1_M2}
For any non-principal leaf $\ell'$ in $\MM_g-\ell_g$, let $\ell^1=\hhat_{1\rm{q}}^{-1}(\ell')$ and $\ell^2=\hhat_{2\rm{q}}^{-1}(\ell')$ be its preimages, which are non-principal leaves. Then $\ell'$, $\ell^1$, and $\ell^2$ are either all isomorphic to solid tori or all isomorphic to $\Hyp^3$. 

The principal leaves $\ell_1=\ell_{f_1}$, $\ell_2=\ell_{f_2}$ and $\ell_g$ are all homeomorphic. Moreover, $\ell_g$ is homeomorphic to the preimages $\ell_{1}-\vv_1=\hhat_{1\rm{q}}^{-1}(\ell_g)$ and $\ell_{2}-\vv_2=\hhat_{2\rm{q}}^{-1}(\ell_g)$.
\end{prop}

Hence the inner topology is leafwise preserved when $f_c$ moves as continuously as above. 

\parag{Remark.}
During the motion of $f_c$, the leafwise homotopies are also preserved. It seems that $f_c$ gives a family of faithful representions of the ``laminated fundamental group" of $\MM_{1}$. 

\parag{}
Now intriguing topological differences would be observed in the lower ends. 
Let us compare the lower ends $\SSS_1$ and $\SSS_2$ of $\Bar{\MM_1}$ and $\Bar{\MM_2}$ as in the previous section.

\parag{Partial conjugacies.}
To investigate the combinatorial change between the lower ends of the quotient 3-laminations, we need to compare the tessellations $\Tess^\a(f_1)$ and $\Tess^\a(f_2)$. 
We start with defining a map which relates the tiles in $\Tess^\a(f_1)$ to those of $\Tess^\a(f_2)$.

Let $h_i:\Chat \to \Chat$ be the semiconjugacy from $f_i~(i=1,2)$ to $g$ given by \thmref{thm_semiconj}.
Then the restriction $h_i|~\Cbar-I_{i} \to \Cbar-I_{g_i}$ is a topological conjugacy. 
Since $I_{g_1}=I_{g_2}$, we have a partial conjugacy $\kappa=h_2^{-1}\cc h_1:\Cbar-I_{1} \to \Cbar-I_{2}$ from $f_1$ to $f_2$ (\cite[Corollary 4.2]{Ka3}). One can easily check that $\kappa:\Cbar-I_1 \to \Cbar-I_2$ is lifted to the partial conjugacies
\begin{align*}
\kaphat_\n&: \NN_1-\Ihat_1 ~\to~ \NN_2-\Ihat_2 \\
\kaphat_\n&~|~ \AAA_1^\n -\Ihat_1 ~\to~ \AAA_2^\n  -\Ihat_2 \\
\kaphat_\a&: \AAA_1 -\II_1 ~\to~ \AAA_2  -\II_2. 
\end{align*} 
In particular, the restriction $\kaphat_\a |~ \DD_1 - \ZZ_1 ~\to~ \DD_2  -\ZZ_2$ is a conjugacy between $\fhat_{1\a}|_{\DD_1 - \ZZ_1}$ and $\fhat_{2\a}|_{\DD_2 - \ZZ_2}$ too. Thus by taking quotient we have a homeomorphism
$$
\kaphat_\q : \SSS_1-\vz_1 \to \SSS_2 -\vz_2.
$$

\parag{Subdivision of tessellation.} 
By comparing the definitions of $\Tess(g_1)$ and $\Tess(g_2)$ (\cite[\S 3]{Ka3}), one can easily check that 
$$
T_{g_2}(\theta, m, \ast) \ee 
\bigcup_{\mu =0}^{q-1}T_{g_1}(\theta,~ m + \mu l,~ \ast)
$$
for any $T_{g_2}(\theta, m, \ast) \in \Tess(g_2)$. Thus $\Tess(g_1)$ is just a subdivision of $\Tess(g_2)$. 

Take a tile $T_1(\theta, m , \ast) \in \Tess(f_1)$. Then there is a homeomorphic image $T_{2}'(\theta, m , \ast):=\kappa(T_{1}(\theta, m , \ast))$ in $K_{2}^\cc=K_{f_2}^\cc$. We say the family
$$
\Tess'(f_2) \dee \skakko{\kappa(T) ~:~ T \in \Tess(f_1)}
$$
is the \textit{subdivided tessellation} of $K_{2}^\cc-I_{2}$. Since $\Tess(f_1)$ and $\Tess(f_2)$ have the same combinatorics as $\Tess(g_1)$ and $\Tess(g_2)$ respectively,
$$
T_{f_2}(\theta, m, \ast) \ee 
\bigcup_{\mu =0}^{q-1}T_{f_2}'(\theta, ~m+\mu l,~ \ast)
$$
for any $T_{f_2}(\theta, \mu, \ast) \in \Tess(f_2)$. Now we have a natural tile-to-tile correspondence between $\Tess(f_1)$, $\Tess(g_1)$ and $\Tess'(f_2)$.

By $\Tess'^\a(f_2)$ we denote the natural subdivisions of $\Tess^\a(f_2)$ associated with $\Tess'(f_2)$. 
More precisely, the {\it subdivided tessellation} $\Tess'^\a(f_2)$ is the family of the {\it subdivided tiles} 
$$
T_2'^\a(\thetahat, m, \ast) \dee \kaphat_\a\kakko{T_1^\a(\thetahat, m, \ast)}, 
$$
where $T_1^\a(\thetahat, m, \ast)$ is a tile of $\Tess^\a(f_1)$. 
Finally we define the subdivided quotient tessellation $\Tess'^\q(f_2)$ as the family of $\kaphat_\q([T])$ with $[T] \in \Tess^\q(f_1)$.

\begin{figure}[htbp]
\centering{\vspace{0cm}
\includegraphics[height=0.19\textheight]{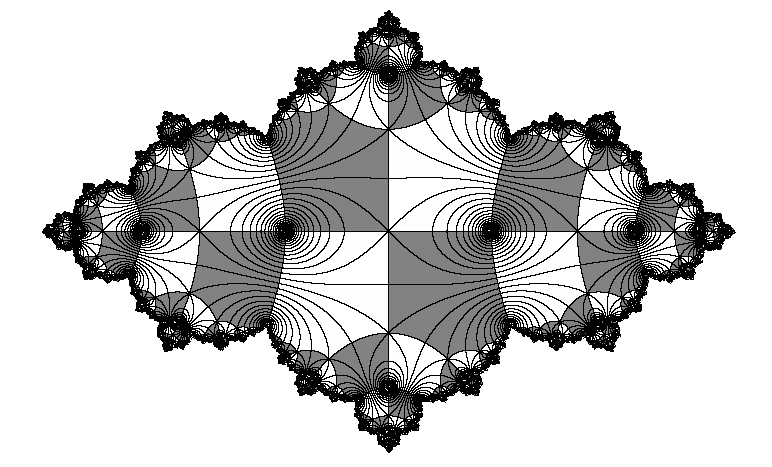} \hspace{-0.5cm}
\includegraphics[height=0.19\textheight]{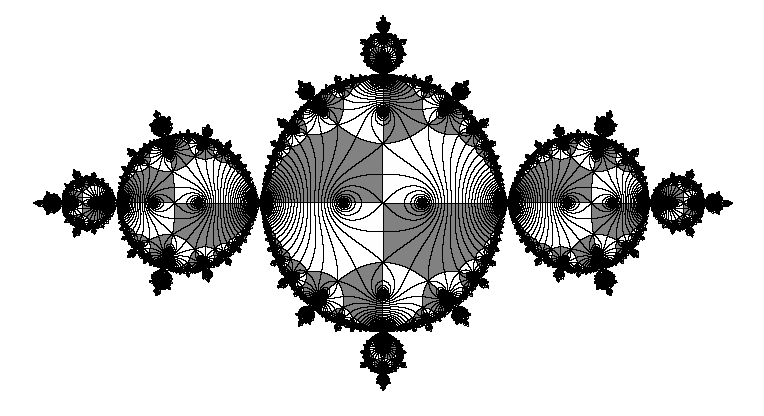}
}
\caption{$\Tess(f_1)$ and $\Tess'(f_2)$ for $g(z)=z^2-3/4$. Compare with Figures 2 and 7 in \cite{Ka3}.}
\label{fig_sub_tess}
\end{figure}

\subsection{Fundamental regions and common properties}
Now we are ready to observe the common properties and the combinatorial difference between the lower ends $\SSS_1=\SSS_{f_1}$ and $\SSS_2=\SSS_{f_2}$ of the quotient 3-laminations by studying combinatorics of the quotient tiles. 

\parag{Fundamental regions.} 
Recall that we may regard $f_1$ as $f$ in the preceding section (``Assumption (a)"). To compare $\DD_1=\DD_{f_1}=\Df$ and $\DD_2=\DD_{f_2}$, let us define the sets in $\DD_2$ corresponding to $\Qf(\ast)=\QQ_1(\ast)$ and $\Pf(\ast)=\PP_1(\ast)$ in the same way. For $\ast=\pm$, we set 
$$
\QQ_2(\ast) \dee \bigcup_{\mu =0}^{q-1} 
\overline{T_2'^\a(\thetahat^+,  \mu l, *)} 
 ~~\text{and}~~
\PP_2(\ast) \dee 
\bigsqcup_{\omehat \in \Omega} 
\kakko{\overline{\Pi_2^\a(\omehat, \ast)}-\{\gam_2^\a(\omehat)\}}
$$
where we set $\Pi_2^\a(\omehat, \ast):=\kaphat_\a(\Pi_1^\a(\omehat, \ast))$, the union of all subdivided tiles of the form $T_2'^\a(\omehat, \mu l, \ast)$ with $\mu \in \Z$. One can check that the union
$$
\QP_{f_2}\ee
\QP_2 \dee \bigcup_{\ast = \pm} \PP_2(\ast) \cup \QQ_2(\ast)
$$
is a fundamental region of $\SSS_2=\DD_2/\fhat_{2\a}$.

\parag{Common properties.} 
For simplicity, we denote the subsets $\QQ_i(\ast)/\fhat_{i\a}$ and $\PP_i(\ast)/\fhat_{i\a}$ in $\SSS_i$ by $\QQ_i^\ast$ and $\PP_i^\ast$ respectively for $i=1$ and $2$. 
We first check the properties that are common to $f_1$ and $f_2$:

\begin{prop}[Common properties of $\bs{\SSS_1}$ and $\bs{\SSS_2}$]
\label{prop_common_prop}
The lower ends $\SSS_1=\DD_1/\fhat_{1\a}$ and $\SSS_2=\DD_2/\fhat_{2\a}$ of $\Bar{\MM_1}$ and $\Bar{\MM_2}$ have the following common properties:
\begin{enumerate}[\rm (1)]
\item 
The (subdivided) quotient tiles in $\Tess^\q(f_1)$ and $\Tess'^\q(f_2)$ tessellate  $\SSS_1 -\vz_1$ and $\SSS_2 -\vz_2$ respectively. 
In particular, the sets $\vz_1$ and $\vz_2$ of degeneration arcs are the unions of the degeneration edges in these (subdivided) quotient tessellations. 

\item The critical and equipotential edges of the quotient tiles in $\SSS_1=\DD_1/\fhat_{1\a}$ have the same combinatorics as those of the subdivided quotient tiles in $\SSS_2=\DD_2/\fhat_{2\a}$. 
\item For $i=1, 2$ and $\ast=\pm$, we have $\QQ_i^\ast \cup \PP_i^\ast=\overline{\PP_i^\ast}$ and $\overline{\PP_i^+} \cup \overline{\PP_i^-}=\SSS_i$. In particular, every quotient panel $[\Pi_i^\a(\omehat, \ast)]$ in $\PP_i^\ast$ accumulates on $\QQ_i^\ast$. 
\item Both $\SSS_1$ and $\SSS_2$ are compact. 
\end{enumerate}
\end{prop}
\begin{pf}
Property (1) is obvious by the definition of the quotient tessellations. 
Property (2) holds since the critical and equipotential edges of tiles in $\Tess'^\a(f_2)$ have the same combinatorics as those of $\Tess^\a(g_1)$ or $\Tess^\a(f_1)$. For (3) and (4), one can just apply the same argument as Propositions \ref{thm_Qf_Pf} and \ref{thm_Sf_Sg} to $f_2$. Details are left to the reader.
\QED
\end{pf}

So we may observe the main difference between $\SSS_1$ and $\SSS_2$ in the connection of the degenerating edges. 

\subsection{Combinatorial difference along degenerating edges}

Let us investigate how the quotient tiles are glued along the degenerating edges in the lower ends $\SSS_1$ and $\SSS_2$.

\parag{Degenerating arcs in the lower ends.}
By \propref{prop_common_prop}(1), the degenerating edges lie on the set of degenerating arcs $\vz_i=\ZZ_i/\fhat_{i\a}~(i=1,2)$. 
As we did for the fundamental domains of the lower ends, it is convenient to introduce the cyclic and non-cyclic decompositions of the degenerating arcs $\vz_1$ and $\vz_2$. 
We claim:
\begin{prop}[Decomposition $\vz_i=\vx_i \sqcup \vy_i$]\label{prop_z_x_y}
For each $i=1$ or $2$, there exists only one closed $\vz_i$-component given by
$$
\vx_i \dee (\DD_i \cap \IO_i)/\fhat_{i\a}
$$
which is contained in $\QQ_i^+ \cap \QQ_i^-$.
The complement $\vy_i:=\vz_i-\vx_i$ is contained in $\PP_i^+ \cap \PP_i^-$.
Moreover, each path-connected component of $\vy_i$ ($\vy_i$-component for short) accumulates on $\vx_i$.
\end{prop}

\begin{pf}
Let us to go back to the dynamics $\fhat_{i\a}: \DD_i \to \DD_i$. 
Recall that $\ZZ_i=\DD_i \cap \II_i$. 
Now the union of cyclic $\ZZ_i$-components is given by $\XX_i:=\DD_i \cap \IO_i$.
By definition of the closed set $\QQ_i(+)$ (resp. $\QQ_i(-)$), any point in $\XX_i$ lands on the degenerating edges of $\QQ_i(+)$ (resp. $\QQ_i(-)$) just once (or twice on the endpoints) by the dynamics of $\fhat_{i\a}$. 
Thus the quotient $\vx_i:=\XX_i/\fhat_{i\a}$ is contained in the boundaries of both $\QQ_i^+$ and $\QQ_i^-$. 
In particular, it is a simple closed curve.  

Since the set $\YY_i:=\ZZ_i -\XX_i$ is the union of non-cyclic $\ZZ_i$-components, its quotient $\vy_i:=\vz_i-\vx_i$ consists of $\vz_i$-components that are homeomorphic to a path-connected component of $D_{f_i} \cap I_{f_i}$ downstairs. 
Hence any $\vz_i$-components except $\vx_i$ cannot be a closed curve. 
Since non-cyclic $\ZZ_i$-components are contained in non-cyclic panels in $\DD_i$, the set $\vy_i$ is contained in the non-cyclic part $\PP_i^+ \cap \PP_i^-$ in the lower end.

Since any panel in $\PP_i^\ast$ accumulates on $\QQ_i^\ast$ (\propref{thm_Qf_Pf}), any $\vy_i$-component accumulates on $\vx_i$. 

\QED
\end{pf}

In the following we consider the connection of the quotient tiles along $\vx_i$ in $\QQ_i^+ \cap \QQ_i^-$ and $\vy_i$ in $\PP_i^+ \cap \PP_i^-$ separately. 
The question is: 
\begin{quote}
\textit{For each $i=1$ or $2$, when do (subdivided) tiles $[T^+] \in \QQ_i^+$ (resp. $\PP_i^+$) and $[T^-] \in \QQ_i^-$ (resp. $\PP_i^-$) share their degenerating edges?} 
\end{quote}
The answer will give us a particular difference between $\SSS_1$ and $\SSS_2$.

\parag{Difference in the cyclic parts $\boldsymbol{\QQ_i^\ast}$.}
By definition, the sets $\QQ_1(\ast)$ and $\QQ_2(\ast)$ consist of the (subdivided) tiles $T_1^\a(\thetahat^+, \mu l, \ast)$ and $T_2'^\a(\thetahat^+, \mu l, \ast)$ with $0 \le \mu < q$. 
We first claim:

\begin{thm}[Twisting $\boldsymbol{\QQ_i^\ast}$]\label{thm_twisting}
In $\SSS_i$, the tiles in $\QQ_i^\ast$ are glued along the degenerating edge $\vx_i:=\XX_i/\fhat_i$ with the following combinatorics:
\begin{enumerate}[\rm (a)]
\item The degenerating edge of the quotient tile $[T_1^\a(\thetahat^+, \mu l, +)]$ in $\QQ_1^+$ is glued to that of $[T_1^\a(\thetahat^+, \mu' l, -)]$ in $\QQ_1^-$ iff $\mu=\mu'$.
\item The degenerating edge of the subdivided quotient tile $[T_2'^\a(\thetahat^+, \mu l, +)]$ in $\QQ_2^+$ is glued to that of $[T_2'^\a(\thetahat^+, \mu' l, -)]$ in $\QQ_2^-$ iff $\mu \equiv \mu'+\tilde{p}$ for some $\tilde{p}$ with $\tilde{p} \, p \equiv 1$ modulo $q$.
\end{enumerate}
\end{thm}

\begin{figure}[htbp]
\centering{
\includegraphics[width=0.7\textwidth]
{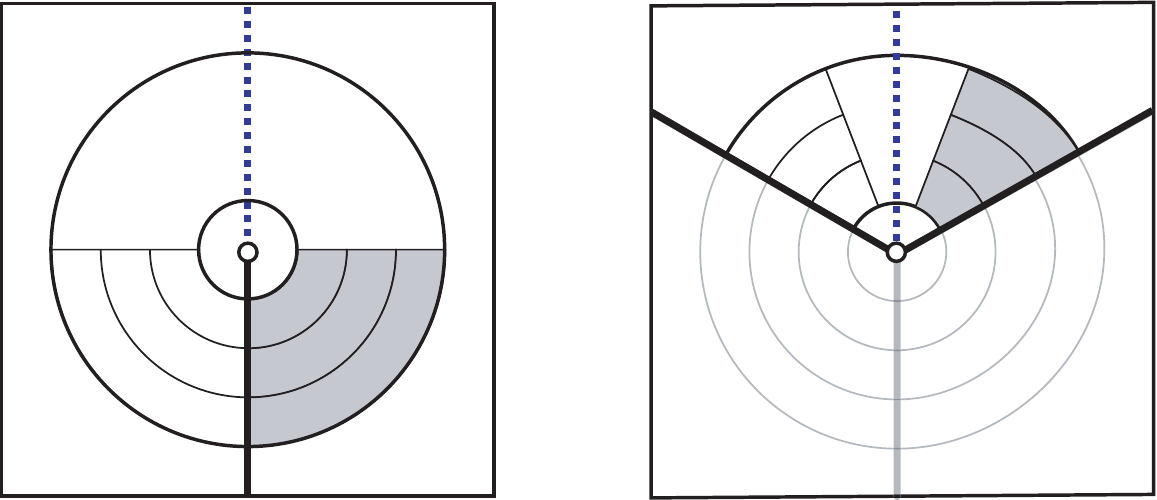}
}
\caption{For $i=1$ (left) and $i=2$ (right), the sets $\QQ_i(+)$ and $\QQ_i(-)$ are subsets of the principal leaves containing the external ray of angle $\thetahat^+$ (the blue dotted lines). 
The thickened lines indicate the degenerating arcs.}
\label{fig_twist}
\end{figure}

\begin{pf}
The statement for $i=1$ is clear by definition of $\QQ_1(\ast)$. Let $L_2$ be the principal leaf that contains the external ray of angle $\thetahat^+$. Then $L_2$ also contains both $\QQ_2(+)$ and $\QQ_2(-)$. Now $L_2 \cap \XX_2$ is given by the non-compact star-like graph $L_2 \cap \IO_2$ with its central joint contained in $\JJ_2$ removed. Thus it has $q$ connected components $I_0, \ldots , I_{q-1}$ in counterclockwise order. Now $L_2 \cap \IO_2$ divide $L_2$ in $q$ sectors and the tiles in both $\QQ_2(+)$ and $\QQ_2(-)$ are contained in one of these sectors that contains the external ray of angle $\thetahat^+$. We may assume that $I_0$ and $I_1$ bounds this sector. Then the degenerating edge of tiles in $\QQ_2(+)$ and $\QQ_2(-)$ are contained in $I_1$ and $I_0$ respectively. (The right square of \figref{fig_twist} shows this situation in $L_2$ for $p/q=1/3$. The left square shows the corresponding picture for $i=1$.) 

By the action of $\fhat_2^{l}:L_2 \to L_2$, the $q$ periodic external rays in $L_2$ and $I_0, \ldots , I_{q-1}$ are permuted with rotation number $p/q$. (cf. \cite[Propositions 2.2 and 2.3]{Ka3}.) In particular, we have $\fhat_2^{lk}(I_{j})=I_{j+kp}$ taking subscript modulo $q$. To glue the degenerating edges of $\QQ_2(+)$ and $\QQ_2(-)$, it is necessary for $I_0$ to be mapped to $I_1$ by some $\fhat_2^{lk}$. Such a $k$ is given by the equation $0+kp \equiv 1$ modulo $q$. This $k$ is our $\tilde{p}$ in the statement.  
\QED
\end{pf}

\parag{Dehn twist along $\boldsymbol{\vx_1}$.} 
The theorem above indicates what is the difference between the cyclic parts $\QQ_1^+ \cup \QQ_1^-$ and $\QQ_2^+ \cup \QQ_2^-$: First cut $\QQ_1^+ \cup \QQ_1^-$ along $\vx_1$, and then twist $\QQ_1^-$ by (combinatorial) angle $\tilde{p}/q \in \T$ and glue them again. Then we have combinatorially the same object as $\QQ_2^+ \cup \QQ_2^-$ with the seam $\vx_2$. We say this operation a combinatorial \textit{$\tilde{p}/q$-Dehn twist} of $\QQ_1^+\cup \QQ_1^-$. In the next paragraph we will see that $\tilde{p}$ is more specified when we look at this change in $\PP_i^+\cup \PP_i^-$.

\parag{Difference in the non-cyclic parts $\boldsymbol{\PP_i^\ast}$.}
Next we consider $\PP_i^\pm$ in the lower end $\SSS_i$ for $i=1$ and $2$. 
Let $T^+ \in \PP_i(+)$ and $T^- \in \PP_i(-)$ be the tiles of addresses $(\omehat, \mu l, +)$ and $(\omehat', \mu' l, -)$ respectively. (If $i=2$ addresses are taken as tiles in $\Tess'^\a(f_2)$.) 
As in the question above, we assume that the quotient tile $[T^+] \in \PP_i^+$ and $[T^-]\in \PP_i^-$ share their degenerating edges. 
Then there exists a unique (subdivided) tile $U^- \subset \DD_i$ of address $(\thetahat, \mu l, -)$ that actually shares its degenerating edge with $T^+$ in the affine lamination before taking the quotient. 
For fixed $i=1$ or $2$, the quotient tiles $[T^+]$ and $[T^-]$ share its degenerating edge in $\SSS_i$ iff there exists a unique integer $\tilde{p}$ such that $\fhat_i^{\tilde{p} l}(T^-)=U^-$. The following theorem says that such a $\tilde{p}$ depends only on $i$ and $p/q$:

\begin{thm}[Dislocating $\boldsymbol{\PP_i}$]\label{thm_twisting2}
For fixed $i=1$ or $2$, let $T^+$, $T^-$ and $U^-$ be the (subdivided) tiles in $\PP_i(+) \cup \PP_i(-)$ defined as above. Then the integer $\tilde{p}$ with $\fhat^{\tilde{p} l}(T^-)=U^-$ is determined as follows:
\begin{enumerate}[\rm (a)]
\item If $i=1$, then $T^-=U^-$ and thus $\tilde{p}=0$. In other words, the quotient $[T^+]$ shares its degenerating edge with $[T^-]$ iff $\omehat=\thetahat=\omehat'$ and $\mu=\mu'$.
\item If $i=2$, then $\tilde{p}$ is uniquely determined so that $|\tilde{p}|<q/2$ and $\tilde{p} \, p \equiv 1$ modulo $q$. In other words, the quotient tile $[T^+]$ shares its degenerating edge with $[T^-]$ iff $\thetahat=\delhat^{\tilde{p}l}(\omehat')$ and $\mu=\mu'+\tilde{p}$.
\end{enumerate}
In particular, $\tilde{p}$ does not depend on the choice of $\omehat$ and $\mu$.
\end{thm}

Note that this theorem is consistent with $\tilde{p}/q$-Dehn twist of $\QQ_1^+ \cup \QQ_1^-$ described in \thmref{thm_twisting}. 

\begin{pf}
By $\Theta_0$ we denote the type of $\gam_2(\theta_0^+)$ (or $\beta_0$); that is, the set of $q$ angles which is cyclic under the iteration of $\delta^l$ and contains $\theta_0^\pm$. For $\omehat =\bon{\omega} \in \Omega$, there exists the maximal $N=N(\omehat)$ such that $\omega_{-Nl} \in \Theta_0$. Since $\omega_0=\theta_0^+$ and $\omega_{-n} \neq \theta_0^+$ for all $n >0$ , we have $0 \le N < q$. In particular, the $q$ angles 
$$
\omega_{-Nl},~\omega_{(-N+1)l}, ~\ldots, ~\omega_0=\theta_0^+,~\ldots,~\omega_{(-N+q-1)l}
$$ 
(where $\omega_k:=\delta^k(\omega_0)$ when $k >0$) coincide with the angles in $\Theta_0$.

Now $T^+ \in \PP_i(+)$ and $T^- \in \PP_i(-)$ have addresses $(\omehat, \mu l, +)$ and $(\omehat', \mu' l, -)$ respectively. 
Set $N:=N(\omehat)$ and $N':=N(\omehat')$. 
Since $U^-$ shares its degeneration edge with $T^+$ in $\DD_i$ before taking the quotient, the address of $U^-$ is of the form $(\thetahat, \mu l, -)$ with non-cyclic $\thetahat=\bon{\thetahat} \in \Thetahat$. 
Note that by the dynamics downstairs tiles of addresses $(\omega_{-n}, \mu l-n,+)$ and $(\theta_{-n}, \mu l-n,-)$ share their degenerating edges for all $n \ge 0$. Moreover, we have 
\begin{itemize}
\item $\theta_{-Nl} \in \Theta_0$ and that $\theta_{-n} \notin \Theta_0$ for all $n > Nl$; and 
\item $\theta_0=\theta_0^+$ or $\theta_0=\theta_0^-$ according to $i=1$ or $2$.
\end{itemize}
In the case of $i=1$, the tiles of addresses $(\omega_{-n}, \mu l-n,+)$ and $(\omega_{-n}, \mu l-n,-)$ share their degeneration edges  so we have $\omega_{-n} = \theta_{-n}$ for all $n \ge 0$. Hence the statement is clear for $i = 1$.

Suppose that $i=2$. Since $\fhat_2^{\tilde{p} l}(T^-)=U^-$ implies $\delhat^{\tilde{p} l}(\omehat')=\thetahat$, we have $\delta^{\tilde{p}l}(\omega'_0)=\theta_{0}$ and $-N'+\tilde{p}=-N$. Recall that $\omega'_0=\theta_0^+$ and $\theta_0=\theta_0^-$. By the same argument as \thmref{thm_twisting}, such a $\tilde{p}$ must satisfy $\tilde{p}\, p \equiv 1 \mod q$. Since $q$ angles
$$
\omega'_{-N'l}~,\ldots, ~\omega'_0=\theta_0^+~, \ldots, ~\omega'_{(-N'+q-1)l}
$$ 
coincide with $\Theta_0$, they contain just one angle equal to $\theta_0^-$. Hence we have $-N'\le \tilde{p} \le -N'+q-1$. Since $N'=\tilde{p}+N$ and $0 \le N<q$, we have $-q/2 < \tilde{p} <q/2$. Now $\tilde{p}$ is uniquely determined and depends neither on $\omehat$ nor $\mu$. \QED

\end{pf}

\parag{The combinatorial $\bs{\tilde{p}/q}$-Dehn twist of the lower ends.}
\propref{prop_common_prop}(1), Theorems \ref{thm_twisting} and \ref{thm_twisting2} above show that the combinatorial difference between $\SSS_1$ and $\SSS_2$ is caused by the connection of the degenerating edges. 
Let us observe it in a continuous motion of the parameter $c$ of $f=f_c$.
Let $c$ move from one hyperbolic component $H_1$ to one of its satellite component $H_2$ via the parabolic parameter $\sigma$. (That is, $\sigma$ is the root of $H_2$.) 
Our degeneration pair $(f_i \to g=f_\sigma)$ represents the motion within $H_i$ for each $i=1$ and $2$. 
\begin{figure}[htbp]
\begin{center}
\includegraphics[width=.75\textwidth]{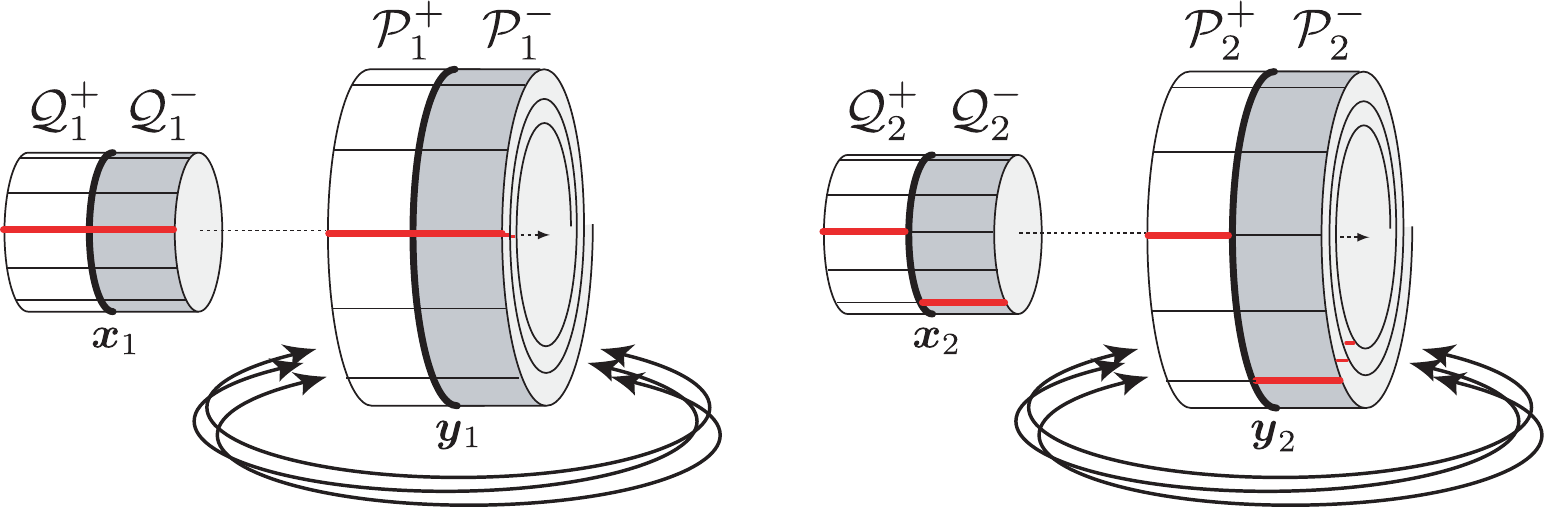}
\end{center}
\caption{Dislocation of $\PP_i$ along twisted $\QQ_i$. 
This caricature only visualize the fact that the panels in $\PP_2^\ast$ are slided by $\tilde{p}l$ tiles. They are also shuffled, thus different pairs of panels are glued along $\vy_1$ and $\vy_2$.
}
\label{fig_twisting_lower_ends}
\end{figure}

As we have seen in \propref{thm_Qf_Pf} and \propref{prop_common_prop}, the sets $\Pg^\ast$ and $\PP_i^\ast$ accumulate on $\Qg^\ast$ and $\QQ_i^\ast$ for all $i=1, 2$ and $\ast=\pm$ respectively. 
Their boundaries are glued in exactly the same way on the critical edges, but there is no edges in $\Qg^\ast \cup \Pg^\ast$ corresponding to the the degenerating edges in $\Qg^\ast \cup \Pg^\ast$ (\figref{fig_winding}). 
More precisely, when $c$ is in $H_1$ first, the degenerating edges $\vz_1$ form a ``laminated path" in $\SSS_1$ which has only one closed path $\vx_1$ with the other $\vy_1$-components winding around (\propref{prop_z_x_y}). 
They are pinched and pushed to ``infinity" as $f_1 \to g~(c \to \sigma)$ then $\SSS_1$ loose compactness at $\Sg$ (\propref{thm_Sf_Sg}). 

Next when $c$ moves into $H_2$, the degenerating edges $\vz_2 = \ZZ_2/\fhat_2$ are again plumped up but the connection of the degenerating edges are changed as Theorems \ref{thm_twisting} and \ref{thm_twisting2} above indicate. 
The difference is the combinatorial $\tilde{p}/q$-Dehn twist along $\vx_1$, and the combinatorial dislocation by $\tilde{p}l$ (subdivided) tiles along $\vy_1$-components. 
In particular, the combination of the glued panels along $\vy_1$ are shuffled  according to the relation $\thetahat=\delhat^{\tilde{p}l}(\omega')$ in \thmref{thm_twisting2}(b).
Such $\tilde{p}$ is uniquely determined by $p/q$ and these two operations are totally consistent with each other. 
So we may consider that the union of these operation is a natural twisting operation of the lamination along $\vz_1$ when $c$ moves from $H_1$ to $H_2$. 
We call this operation is the \textit{combinatorial $\tilde{p}/q$-Dehn twist of the lower ends}.

\subsection{Topological difference of the Kleinian lamination and the lower ends} 
The phenomenon above looks similar to the example of quasi-Fuchsian groups in \secref{sec_04}, but there is a particular difference also: The total topology of the Kleinian lamination may change. 

In the quasi-Fuchsian example the Kleinian manifold $\Bar{M_t}$ has the same product structure $\approx S \times [0,1]$ for all $0<|t|<1$. (See \cite[Theorem 4.34]{MT} for example). This fact corresponds to:

\begin{thm}[Product structure]\label{thm_prod_str_vs_main_cardioid}
Let $f = f_c$ be a quadratic map. Then the Kleinian lamination $\Bar{\MM_{f}}$ is homeomorphic to the product $\SSS_0 \times [0,1]$ if and only if $c$ is in the main cardioid of the Mandelbrot set. 
\end{thm}

\begin{pf}
By \thmref{thm_prod_str_of_M}, the Kleinian lamination $\Bar{\MM_{f}}$ has a product structure $\SSS_0 \times [0,1]$ when $c$ is in the main cardioid. 

Conversely, suppose that there exists a topological homeomorphism $H:\Bar{\MM_{f}} \to \SSS_0 \times [0,1]$.  
Then the conformal boundary $\partial \MM_{f}$ must consist of two indecomposable connected components.

If the Fatou set has only one connected component ($ = B_f$, the basin at infinity) downstairs, then $\partial \MM_{f} = \Ff/\fhat$ cannot have two components. Thus the interior $D_f$ of the filled Julia set $K_f$ is not empty.

Let $\gam = \gam_{f}(0)$ be the $\beta$-fixed point of $f$. 
Then the leaf $\ell = \ell_\gam$ of $\MM_{f}$ associated with this repelling fixed point is a solid torus. 
The leafwise conformal boundary $\partial \ell$ has two indecomposable components corresponding to $B_f$ and $D_f$, the upper and lower ends. 
Since we have the homeomorphism $H$, the image $H(\ell)$ must be a solid torus and the boundary of $H(\ell)$ must be two annuli. This can happen only if $c$ is in the main cardioid, otherwise the lower end of $\ell$ has countably many components or $c=1/4$.
When $c = 1/4$, the lower end $S_f$ is non-compact by \thmref{thm_Sf_Sg}. This contradicts the fact that $\SSS_0$ is compact, and that $H$ sends the boundary component to the boundary component.
\QED
\end{pf}
This yields \thmref{thm_introduction_2} in the introduction.
By using this theorem, we immediately have the following example:
\begin{cor}\label{cor_lower_ends_of_p/q_limbs}
Let $f_1$ be in the main cardioid, and $f_2$ in the $p/q$-limb of the Mandelbrot set. Then the Kleinian lamination $\Bar{\MM_1}$ is not homeomorphic to $\Bar{\MM_2}$. 
\end{cor}

Now it is natural to ask:

\parag{Question.}
{\it For any pair of Case (a) degeneration pair $(f_1 \to g)$ and Case (b) $(f_2 \to g)$, can the Kleinian lamination $\Bar{\MM_1}$ be homeomorphic to $\Bar{\MM_2}$? 
}

\parag{}
It suffices to show that the lower end $\SSS_1$ cannot be homeomorphic to $\SSS_2$. 
However, the author do not have the answer even for the following question: 

\parag{Question.}
{\it Can Sullivan's lamination $\SSS_0$ be homeomorphic to the lower end  $\Sf$ of some hyperbolic quadratic map $f = f_c$ with $c$ outside the main cardioid? 
}

\subsection{Example: The Rabbits}\label{subsec_rabbits}
To conclude this paper let us consider a typical degeneration and bifurcation process when the parameter $c$ moves from $0$ to the center of the $1/3$-limb (Douady's rabbit, with the value denoted by $c_{\mathrm{rab}}$) of the Mandelbrot set. (See Figure 1 of Part I \cite{Ka3}.) 

Let $(f_1 \to g)$ be the Case (a) Rabbits and $(f_2 \to g)$ the Case (b) Rabbits. 
Set $f_0(z)=z^2$ and $f_3(z)=z^2+c_{\mathrm{rab}}$. 
By \corref{cor_hyp_center} and \thmref{thm_prod_str_of_M}, we have:
\textit{
\begin{itemize}
\item[\bf (R1)] The Kleinian lamination $\overline{\MM_{f_0}}$ (resp. $\overline{\MM_{f_3}}$) is quasi-isometrically equivalent to $\overline{\MM_{f_1}}$ (resp. $\overline{\MM_{f_2}}$). 
In particular, $\overline{\MM_{f_1}}$ has a product structure $\approx \SSS_0 \times [0,1]$.
\end{itemize}
}
Thus it is enough to check the changes of laminations associated with $(f_1 \to g)$ and $(f_2 \to g)$ to compare the Kleinian laminations of $f_0(z)=z^2$ and Douady's rabbit.

\parag{Quotient 3-laminations.}
The principal leaf $\ell_1=\Lam_1^\h/\fhat_1$ (a solid torus) is pinched to be the principal leaf $\ell_g =\Lam_g^\h/\ghat$ (a rank one cusp), then blown up to become the principal leaf $\ell_2=\Lam_2^\h/\fhat_2$ (again a solid torus). 

By \thmref{thm_quotient} we have:
\textit{
\begin{itemize}
\item[\bf (R2)] For each $i=1$ or $2$, there exists a very meek map $\hhat_{i\mathrm{q}}: \MM_i-\ell_i \to \Mg-\ell_g$.
\item[\bf (R3)] The principal part $\ell_1$, $\ell_g$, and $\ell_2$ have the same topology homeomorphic to $\D \times \T$.
\end{itemize}
}
Thus we may conclude:\textit{
During the degeneration-and-bifurcation process from $f_0$ to $f_3$ via $f_1$, $g$ and $f_2$, all the leaves of the quotient 3-lamination are topologically preserved. 
}

In addition, let us remark that the same argument as the proof of (C2)' yields an upgrade of (R2):
\textit{
\begin{itemize}
\item[\bf (R2)'] The non-principal part $\Bar{\MM_2}-\Bar{\ell_2}$ is homeomorphic to $\Bar{\Mg}-\Bar{\ell_g}$.
\end{itemize}
}

\parag{Conformal boundaries.}
We may find more significant difference in the conformal boundaries. 
An intriguing fact is, though this motion of the parameter passes across the hyperbolic components, the situation is still analogical to the quasi-Fuchsian deformation.

For $i=1$ or $2$, let $\SSS^i$ and $\SSS^g$ (resp. $\SSS_i$ and $\SSS_g$) be the upper (resp. lower) ends of $\MM_i$ and $\Mg$.
By Propositions \ref{prop_Bf_Df}, \ref{thm_Sf_Sg}, and \ref{prop_common_prop}, we have
\textit{
\begin{itemize}
\item[\bf (R4)] The upper ends $\SSS^i$ and $\SSS^g$ are all conformally equivalent to Sullivan's solenoidal Riemann surface lamination $\SSS_0$. 
\item[\bf (R5)] The lower ends $\SSS_1$ and $\SSS_2$ are compact, but $\Sg$ is non-compact.
\end{itemize}
}

\begin{figure}[htbp]
\centering{
\vspace{0cm}
\includegraphics[width=0.8\textwidth]{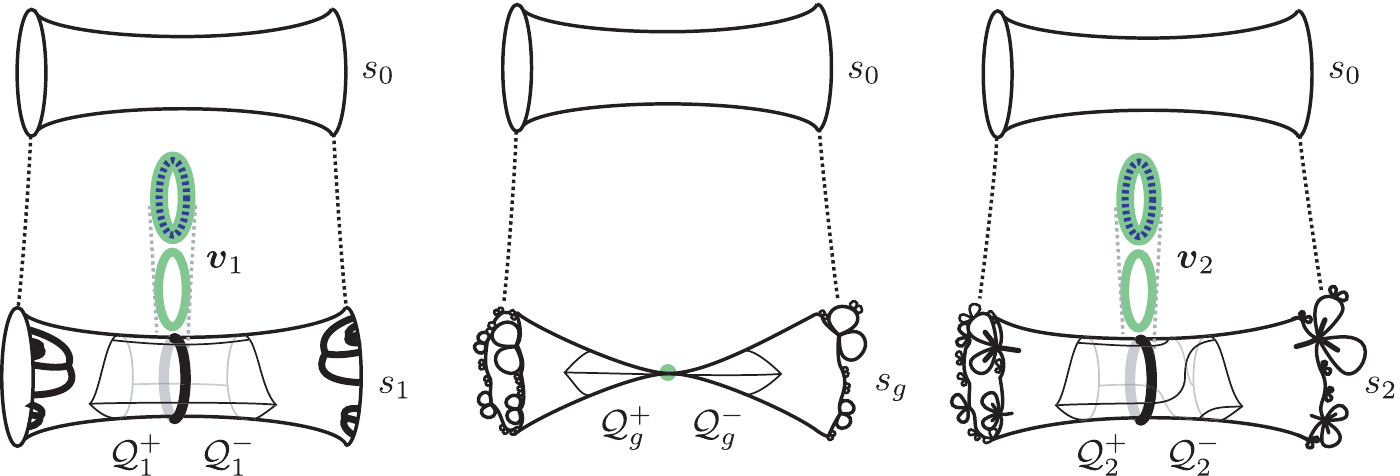}
}
\caption{Caricatures of $\Bar{\ell_1}$, $\Bar{\ell_g}$, and $\Bar{\ell_2}$ from left to right. 
}\label{fig_ell_4_rabbits}
\end{figure}

\parag{Product structure.}
Next we check the product structures of $\overline{\MM_i}~(i=1,2)$ and $\overline{\MM_g}$.
By \thmref{thm_prod_str_vs_main_cardioid} we have:
\textit{
\begin{itemize}
\item[\bf (R6)] The quotient 3-lamination with conformal boundary $\overline{\MM_1}$ has a product structure $\approx \SSS_0\times [0,1]$, but $\overline{\MM_g}$ and $\overline{\MM_2}$ do not. 
\end{itemize}
}
See \figref{fig_ell_4_rabbits}. By summing up (R5) and (R6) (or by Corollaries \ref{cor_Kleinian_not_homeo} and \ref{cor_lower_ends_of_p/q_limbs}), we have:
\textit{
\begin{itemize}
\item[\bf (R7)] The Kleinian laminations $\Bar{\MM_1}$, $\Bar{\Mg}$, and $\Bar{\MM_2}$ are not homeomorphic each other. 
\end{itemize}
}
In particular, we have:
\textit{
\begin{itemize}
\item[\bf (R7)'] The Kleinian laminations of $f_0(z)=z^2$ is not homeomorphic to that of Douady's rabbit. 
\end{itemize}
}

Now we have \thmref{thm_introduction_1} in the introduction by (C6), (R1), and (R7).

\parag{Dehn twist in the lower ends.}
The most intriguing part is the combinatorial difference between $\SSS_1$ and $\SSS_2$ given by the twisting (and sliding) operation along $\vz_i=\ZZ_i/\fhat_i$. 

By the argument of Sections \ref{sec_09} and this section, $\SSS_1$ and $\SSS_2$ have the cyclic and non-cyclic decomposition 
$$
\SSS_i \ee \bigcup_{\ast=\pm} \QQ_i^\ast \cup \bigcup_{\ast=\pm} \PP_i^\ast~~~~(i=1,2).
$$
As $f_1 \to g$, the union of the quotient degenerating arcs $\vz_1$ is pinched and pushed away to ``infinity". 
As $f_c$ moves from $g$ to $f_2$, the quotient degenerating arcs appears again as $\vz_2$ but the connection of (subdivided) tiles along $\vz_1$ and $\vz_2$ are different. 
By applying \thmref{thm_twisting} and \thmref{thm_twisting2} to our case of $p/q = 1/3$ and $l = 1$, we have $\tilde{p} = 1$ and conclude:
\textit{
\begin{itemize}
\item[\bf (R8)] During the degeneration and bifurcation process from $f_1$ to $f_2$ via $g$, $\vz_1$ is pinched and then plumped with combinatorial $\tilde{p}/q = 1/3$-Dehn twist to become $\vz_2$. 
\end{itemize}
}
See Figures \ref{fig_winding} and \ref{fig_twisting_lower_ends} again.

\parag{Principal leaves.}
Let us observe (R8) in the principal leaves $\Bar{\ell_1}$, $\Bar{\ell_g}$, and $\Bar{\ell_2}$ with conformal boundaries.
Their interiors are topologically the same by (R3). 
Moreover, by (R4) their leafwise upper ends are conformally the same annuli, say $s_0$. 
(Note that this $s_0$ is different from $s_0$ in the Cauliflowers: the modulus is three times larger.) 
Hence the topological difference is concentrated in the lower ends. 

Let $s_1$, $s_g$, and $s_2$ denote the leafwise lower ends of $\Bar{\ell_1}$, $\Bar{\ell_g}$, and $\Bar{\ell_2}$ respectively. 
As shown in the left of \figref{fig_ell_4_rabbits}, $s_1$ is an annulus and it has countably many $\vz_1$-components (black thick curves). It also contains the cyclic part $\QQ_1^+ \cup \QQ_1^-$ of the fundamental region $\QQ\PP_1$. 
As $f_1 \to g$, these $\vz_1$-components are pinched and we have non-compact lower end $s_g$, which consists of countably many path-connected components. 
In particular, there are two of them containing $\QQ_g^+$ and $\QQ_g^-$, the cuspidal parts. 
Then $g$ is perturbed to be $f_2$ and the $\vz_2$-components appear. 
The $\vz_1$ and $\vz_2$ consist of the degenerating edges of tiles in $\Tess^\q(f_1)$ and $\Tess'^\q(f_2)$, but the connection of these edges are different as described in (R8).

\parag{Remark: Rabbits in rabbit, etc.} 
The similar change happens when $c$ moves from $c_3 = c_{\mathrm{rab}}$ to another hyperbolic component $H'$ via another parabolic parameter. Mainly the topology of the quotient 3-lamination is leafwise preserved, but the lower end changes its combinatorial structure when it goes into $H'$. The location and rotation number of Dehn twist changes according to the angles of the landing ray at the root of $H'$. 

Now it is natural to consider what may happen after finite or infinitely many bifurcations via parabolic parameters. 
This type of bifurcations are described in terms of \textit{renormalization}. 
At least for the affine laminations, we already have a reasonable picture. See \cite{CK}.

\setcounter{section}{0}
\renewcommand{\thesection}{\Alph{section}}
\section{Appendix: Remarks on H\'enon mappings}\label{sec_11}
For hyperbolic $f_c(z)=z^2+c$ and $a \in \Cstar$, we define the \textit{H\'enon map} $H_{c,a}: \C^2 \to \C^2$ by
$$
H_{c,a}: \, \vz=\Tate{x}{y}~\mapsto~H_{c,a}(\vz)=\Tate{f_c(x)-ay}{x}.
$$
This $H_{c,a}$ has the constant Jacobian $a$ thus invertible. When $|a| \ll 1$, this map is considered as a 2-dimensional perturbation of $f_c(z)=z^2+c$. We say $\vz \in K_\pm(H_{c,a})$ if $\norm{H_{c,a}^{n}(\vz)}$ is bounded as $n \to \pm \infty$. We set $J_\pm(H_{c,a}):= \partial K_\pm(H_{c,a})$, the \textit{forward and backward Julia sets}. They are all invariant under the dynamics of $H_{c,a}$.

In \cite{HO}, Hubbard and Oberste-Vorth showed the following:

\begin{thm}[Theorem 5.1 of \cite{HO}]\label{thm_HO}
If $|a| \ll 1$, there exists a homeomorphism $\Phi_-: \NN_{f_c}-\skakko{\hat{\infty}} \to J_-(H_{c,a})$ which conjugates $\fhat_c$ and $H_{c,a}$.
\end{thm}

Note that $\NN_{f_c}-\skakko{\hat{\infty}}$ is the natural extension of $f:\C \to \C$. Hence $J_-(H_{c,a})$ minus the attracting cycle of $H_{c,a}$ is a (topological) realization of the affine part $\AAA_{f_c}^\n$ and the affine lamination $\AAA_{f_c}$ in $\C^2$. (Recall that $\Afn$ and $\Af$ are identified when $f$ is hyperbolic.) This means that our investigation on natural extensions and affine parts in \secref{sec_06} is also an investigation on $J_-$ of some H\'enon maps. 

For example, take two distinct degeneration pairs $(f_1 \to g)$ and $(f_2 \to g)$ such that $c_i$ of $f_i(z)=z^2+c_i$ are contained in the distinct hyperbolic components for $i=1$ and $2$. Then there exists a sufficiently small $a \neq 0$ such that the H\'enon maps $H_1=H_{c_1, a}$ and $H_2=H_{c_2, a}$ have properties as in the theorem above. 
The maps $H_1$ and $H_2$ can be arbitrarily close in the parameter space, but they have different dynamics on $J_-$ by our investigation on affine parts and affine laminations in \secref{sec_06}. In particular, both $J_-(H_1)$ and $J_-(H_2)$ support tessellation induced by $\Tess^\n(f_1)$ and $\Tess^\n(f_2)$ respectively. 

Here we give a brief application:
\begin{thm}[Distinct $\bs{J^-}$.]\label{thm_HO}
The backawrd Julia sets $J_-(H_1)$ and $J_-(H_2)$ are not homeomorphic.
\end{thm}
\begin{pf}
By the theorem above, $J_-(H_i)$ is homromorphic to $\NN_{f_i}-\skakko{\infty}$. The irregular points of $\NN_{f_i}-\skakko{\infty}$ are the cyclic lift of the attracting cycle of $f_i$. Since $f_1$ and $f_2$ have different number of attracting points, the sets $\NN_{f_1}-\skakko{\infty}$ and $\NN_{f_2}-\skakko{\infty}$ cannot be homromorphic. (Note that an irregular point is characterized as a point that cannot have a laminar neighbourhood with product structure.)
\QED
\end{pf}

\paragraph{Stable and unstable manifolds.}
Conversely, one can translate some notion for the H\'enon mappings to our affine laminations. For example, we can find ``stable and unstable manifolds" in the affine laminations. Take a repelling fixed point $z$ of $f$ and let $\zhat$ be its lift in $\Afn$. Then the action of $\fhat:L(\zhat) \to L(\zhat)$ is conjugate to an affine expanding map on $\C$. We may consider it the ``unstable manifold" of $\zhat$. On the other hand, there is a fiber $\TT_z:=\pi^{-1}(z)$ of $z$ invariant under $\fhat$. One can easily check that the action $\fhat: \TT_z \to \TT_z$ is contracting in the natural topology. In particular, we have $\fhat^n(\zetahat) \to \zhat$ for any $\zetahat \in \TT_z$. We may consider this fiber direction the ``stable manifold" of $\zhat$. 

The notion of ``stable and unstable manifolds" in the affine laminations are not only for quadratic maps. More generally, it may be intriguing to consider the following problem: \textit{For a rational map $f$ of degree $\ge 2$, formulate the stable and unstable manifolds in the universal dynamics $\fhat_\u:\UUhat^\a \to \UUhat^\a$.} See also paragraph 10 of \cite[\S 10]{LM} and Appendix of \cite{BS}.

\end{document}